%=================================================================
%
% New combinatorial formula for modified Hall-Littlewood polynomials
%
% Anatol N. Kirillov, January 1998, last April 1998
%
% PlainTEX, 39p.
% 
%=================================================================
%
% Fill in the data (title, authors name, etc.) of your manuscript.
% The places to be filled in occur immediately after the symbols '!',
% hence you can search for '!' one by one.
\font\tenbf=cmbx10

\font\eightrm=cmr8
\font\eightit=cmti8
\font\germ=eufm10
\def\g{\hbox{\germ g}}
\def\b{\hbox{\germ b}}

\def\sectiontitle#1\par{\vskip0pt plus.1\vsize\penalty-250
\vskip0pt plus-.1\vsize\bigskip\vskip\parskip
\message{#1}\leftline{\tenbf#1}\nobreak\vglue 5pt}
\def\ds{\displaystyle}
\def\wh{\widehat}
\def\wt{\widetilde}
\def\eno{\eqalignno}
\def\ld{\lambda}

\def\al{\alpha}

\def\mod{\hbox{\rm mod~}}

\def\frac#1#2{{#1\over#2}}

\magnification=\magstep1
%==================================================
%IT IS BOX.TEX
%==================================================
\def\m@th{\mathsurround=0pt}

\def\fsquare(#1,#2){
\hbox{\vrule$\hskip-0.4pt\vcenter to #1{\normalbaselines\m@th
\hrule\vfil\hbox to #1{\hfill$\scriptstyle #2$\hfill}\vfil\hrule}$\hskip-0.4pt
\vrule}}

\def\addsquare(#1,#2){\hbox{$
	\dimen1=#1 \advance\dimen1 by -0.8pt
	\vcenter to #1{\hrule height0.4pt depth0.0pt\vss%
	\hbox to #1{\hss{%
	\vbox to \dimen1{\vss%
	\hbox to \dimen1{\hss$\scriptstyle~#2~$\hss}%
	\vss}\hss}%
	\vrule width0.4pt}\vss%
	\hrule height0.4pt depth0.0pt}$}}

\def\Fsquare(#1,#2){
\hbox{\vrule$\hskip-0.4pt\vcenter to #1{\normalbaselines\m@th
\hrule\vfil\hbox to #1{\hfill$#2$\hfill}\vfil\hrule}$\hskip-0.4pt
\vrule}}

\def\Addsquare(#1,#2){\hbox{$
	\dimen1=#1 \advance\dimen1 by -0.8pt
	\vcenter to #1{\hrule height0.4pt depth0.0pt\vss%
	\hbox to #1{\hss{%
	\vbox to \dimen1{\vss%
	\hbox to \dimen1{\hss$~#2~$\hss}%
	\vss}\hss}%
	\vrule width0.4pt}\vss%
	\hrule height0.4pt depth0.0pt}$}}

\def\hfourbox(#1,#2,#3,#4){%
	\fsquare(0.3cm,#1)\addsquare(0.3cm,#2)\addsquare(0.3cm,#3)\addsquare(0.3cm,#4)}

\def\Hfourbox(#1,#2,#3,#4){%
	\Fsquare(0.4cm,#1)\Addsquare(0.4cm,#2)\Addsquare(0.4cm,#3)\Addsquare(0.4cm,#4)}

\def\hthreebox(#1,#2,#3){%
	\fsquare(0.3cm,#1)\addsquare(0.3cm,#2)\addsquare(0.3cm,#3)}

\def\htwobox(#1,#2){%
	\fsquare(0.3cm,#1)\addsquare(0.3cm,#2)}

\def\vfourbox(#1,#2,#3,#4){%
%	\hbox{
	\normalbaselines\m@th\offinterlineskip
	\vtop{\hbox{\fsquare(0.3cm,#1)}
	      \vskip-0.4pt
	      \hbox{\fsquare(0.3cm,#2)}	
	      \vskip-0.4pt
	      \hbox{\fsquare(0.3cm,#3)}	
	      \vskip-0.4pt
	      \hbox{\fsquare(0.3cm,#4)}}}%}

\def\Vfourbox(#1,#2,#3,#4){%
%	\hbox{
	\normalbaselines\m@th\offinterlineskip
	\vtop{\hbox{\Fsquare(0.4cm,#1)}
	      \vskip-0.4pt
	      \hbox{\Fsquare(0.4cm,#2)}	
	      \vskip-0.4pt
	      \hbox{\Fsquare(0.4cm,#3)}	
	      \vskip-0.4pt
	      \hbox{\Fsquare(0.4cm,#4)}}}%}

\def\vthreebox(#1,#2,#3){%
%	\hbox{
	\normalbaselines\m@th\offinterlineskip
	\vtop{\hbox{\fsquare(0.3cm,#1)}
	      \vskip-0.4pt
	      \hbox{\fsquare(0.3cm,#2)}	
	      \vskip-0.4pt
	      \hbox{\fsquare(0.3cm,#3)}}}%}

\def\vtwobox(#1,#2){%
%	\hbox{
	\normalbaselines\m@th\offinterlineskip
	\vtop{\hbox{\fsquare(0.3cm,#1)}
	      \vskip-0.4pt
	      \hbox{\fsquare(0.3cm,#2)}}}%}

\def\Hthreebox(#1,#2,#3){%
	\Fsquare(0.4cm,#1)\Addsquare(0.4cm,#2)\Addsquare(0.4cm,#3)}

\def\Htwobox(#1,#2){%
	\Fsquare(0.4cm,#1)\Addsquare(0.4cm,#2)}

\def\Vthreebox(#1,#2,#3){%
	\normalbaselines\m@th\offinterlineskip
	\vtop{\hbox{\Fsquare(0.4cm,#1)}
	      \vskip-0.4pt
	      \hbox{\Fsquare(0.4cm,#2)}	
	      \vskip-0.4pt
	      \hbox{\Fsquare(0.4cm,#3)}}}

\def\Vtwobox(#1,#2){%
	\normalbaselines\m@th\offinterlineskip
	\vtop{\hbox{\Fsquare(0.4cm,#1)}
	      \vskip-0.4pt
	      \hbox{\Fsquare(0.4cm,#2)}}}

\def\twoone(#1,#2,#3){%
	\normalbaselines\m@th\offinterlineskip
	\vtop{\hbox{\htwobox({#1},{#2})}
	      \vskip-0.4pt
	      \hbox{\fsquare(0.3cm,#3)}}}

\def\Twoone(#1,#2,#3){%
	\hbox{
	\normalbaselines\m@th\offinterlineskip
	\vtop{\hbox{\Htwobox({#1},{#2})}
	      \vskip-0.4pt
	      \hbox{\Fsquare(0.4cm,#3)}}}}

\def\threeone(#1,#2,#3,#4){%
	\normalbaselines\m@th\offinterlineskip
	\vtop{\hbox{\hthreebox({#1},{#2},{#3})}
	      \vskip-0.4pt
	      \hbox{\fsquare(0.3cm,#4)}}}

\def\Threeone(#1,#2,#3,#4){%
	\normalbaselines\m@th\offinterlineskip
	\vtop{\hbox{\Hthreebox({#1},{#2},{#3})}
	      \vskip-0.4pt
	      \hbox{\Fsquare(0.4cm,#4)}}}

\def\Threetwo(#1,#2,#3,#4,#5){%
	\normalbaselines\m@th\offinterlineskip
	\vtop{\hbox{\Hthreebox({#1},{#2},{#3})}
	      \vskip-0.4pt
	      \hbox{\Htwobox({#4},{#5})}}}

\def\threetwo(#1,#2,#3,#4,#5){%
	\normalbaselines\m@th\offinterlineskip
	\vtop{\hbox{\hthreebox({#1},{#2},{#3})}
	      \vskip-0.4pt
	      \hbox{\htwobox({#4},{#5})}}}

\def\twotwo(#1,#2,#3,#4){%
	\normalbaselines\m@th\offinterlineskip
	\vtop{\hbox{\htwobox({#1},{#2})}
	      \vskip-0.4pt
	      \hbox{\htwobox({#3},{#4})}}}

\def\Twotwo(#1,#2,#3,#4){%
	\normalbaselines\m@th\offinterlineskip
	\vtop{\hbox{\Htwobox({#1},{#2})}
	      \vskip-0.4pt
	      \hbox{\Htwobox({#3},{#4})}}}

\def\twooneone(#1,#2,#3,#4){%
	\normalbaselines\m@th\offinterlineskip
	\vtop{\hbox{\htwobox({#1},{#2})}
	      \vskip-0.4pt
	      \hbox{\fsquare(0.3cm,#3)}
	      \vskip-0.4pt
	      \hbox{\fsquare(0.3cm,#4)}}}

\def\Twooneone(#1,#2,#3,#4){%
	\normalbaselines\m@th\offinterlineskip
	\vtop{\hbox{\Htwobox({#1},{#2})}
	      \vskip-0.4pt
	      \hbox{\Fsquare(0.4cm,#3)}
	      \vskip-0.4pt
	      \hbox{\Fsquare(0.4cm,#4)}}}

\def\Twotwoone(#1,#2,#3,#4,#5){%
	\normalbaselines\m@th\offinterlineskip
	\vtop{\hbox{\Htwobox({#1},{#2})}
	      \vskip-0.4pt
	      \hbox{\Htwobox({#3},{#4})}
              \vskip-0.4pt
	      \hbox{\Fsquare(0.4cm,#5)}}}

\def\twotwoone(#1,#2,#3,#4,#5){%
	\normalbaselines\m@th\offinterlineskip
	\vtop{\hbox{\htwobox({#1},{#2})}
	      \vskip-0.4pt
	      \hbox{\htwobox({#3},{#4})}
              \vskip-0.4pt
	      \hbox{\fsquare(0.3cm,#5)}}}

\def\Threeoneone(#1,#2,#3,#4,#5){%
	\normalbaselines\m@th\offinterlineskip
	\vtop{\hbox{\Hthreebox({#1},{#2},{#3})}
	      \vskip-0.4pt
	      \hbox{\Fsquare(0.4cm,#4)}
              \vskip-0.4pt
	      \hbox{\Fsquare(0.4cm,#5)}}}

\def\threeoneone(#1,#2,#3,#4,#5){%
	\normalbaselines\m@th\offinterlineskip
	\vtop{\hbox{\hthreebox({#1},{#2},{#3})}
	      \vskip-0.4pt
	      \hbox{\fsquare(0.3cm,#4)}
              \vskip-0.4pt
	      \hbox{\fsquare(0.3cm,#5)}}}
%------------------------------------------------------------------
% \end box.tex
%-------------------------------------------------------------------
%\parindent=15pt
\nopagenumbers
\vglue 5pc %
\vglue 1pc
\baselineskip=13pt
\headline{\ifnum\pageno=1\hfil\else%
{\ifodd\pageno\rightheadline \else \leftheadline\fi}\fi}
\def\rightheadline{\hfil\eightit 
New combinatorial formula for modified Hall--Littlewood polynomials
\quad\eightrm\folio}
\def\leftheadline{\eightrm\folio\quad 
\eightit 
Anatol N. Kirillov 
%! Running author (even)
\hfil}
\voffset=1\baselineskip%
\centerline{\tenbf 
NEW  \hskip 0.1cm COMBINATORIAL  \hskip 0.1cm FORMULA \hskip 0.1cm FOR}
 
\centerline{\tenbf
MODIFIED \hskip 0.1cm HALL--LITTLEWOOD \hskip 0.1cm POLYNOMIALS}
\vglue 12pt%
%\vglue 2pt
\centerline{\eightrm 
ANATOL N. KIRILLOV
%! AUTHOR: 
%\footnote"$^*$"{\eightrm \baselineskip=10pt    
%! present address
% }
}
\vglue 6pt
%\vglue 2pt
\baselineskip=12pt
\centerline{\eightit
CRM, University of Montreal
%! Use the permanent address (University Name, etc.)
}
\baselineskip=10pt
\centerline{\eightit 
C.P. 6128, Succursale A, Montreal (Quebec) H3C 3J7, Canada 
}
\baselineskip=10pt
\centerline{\eightit
and }
\baselineskip=10pt
\centerline{\eightit 
Steklov Mathematical Institute,
%! Use the permanent address (University Name, etc.)
}
\baselineskip=10pt
\centerline{\eightit 
Fontanka 27, St.Petersburg, 191011, Russia
%! The permanent address (City, State ZIP/Zone, Country)
}
%More than two authors, see the sample prints.
\vglue 30pt
\centerline{\it Dedicated to Richard Askey on the occasion of his 
65 birthday}
\vglue 30pt
\vglue 10pt
\centerline{\eightrm ABSTRACT}
\vglue 10pt
{\rightskip=2.5pc
\leftskip=2.5pc
%{\rightskip=1.5pc
%\leftskip=1.5pc
\eightrm\parindent=1pc 
We obtain new combinatorial formulae for modified Hall--Littlewood polynomials, 
for matrix elements of the transition matrix between the elementary 
symmetric polynomials and Hall-Littlewood's ones, 
and for the number of rational points over the finite field of unipotent 
partial flag variety. The definitions and examples of  
generalized mahonian statistic on the set of transport matrices and 
dual mahonian statistic on the set of transport (0,1)--matrices are given.
%These statistics are a natural generalization of definition of the mahonian 
%statistic on {\it words} introduced by D.~Foata.
We also review known $q$--analogues of Littlewood--Richardson numbers and 
consider their possible generalizations. Some conjectures about 
multinomial fermionic formulae for homogeneous unrestricted one 
dimensional sums and generalized Kostka--Foulkes polynomials are 
formulated. Finally we suggest two parameter deformations of polynomials 
${\cal P}_{\ld\mu}(t)$ and one dimensional sums.
} 
%\vglue12pt%
%\vglue10pt
\vglue 4pt
\overfullrule=1pt
\def\qed{\hfill$\vrule height 2.5mm width 2.5mm depth 0mm$}
\vfill\eject
\baselineskip=13pt
{\bf Contents}
\vskip 0.3cm
 
\S 0. Introduction
\vskip 0.2cm

\S 1. Modified Hall-Littlewood polynomials
\vskip 0.1cm

\hskip 0.3cm 1.1. Definition

\hskip 0.3cm 1.2. Modified Hall--Littlewood polynomials for partition 
$\ld =(1^N)$

\hskip 0.3cm 1.3. Hall-Littlewood polynomials and characters of the 
affine Lie algebra $\wh{sl}(n)$

\hskip 0.3cm 1.4. Modified Hall-Littlewood polynomials and unipotent flag 
varieties

\hskip 0.3cm 1.5. Modified Hall--Littlewood polynomials and Demazure 
characters

\hskip 0.3cm 1.6. Modified Hall--Littlewood polynomials and chains of 
subgroups in a finite 

\hskip 1cm  abelian $p$--group

\vskip 0.2cm
\S 2. Generalized mahonian statistics
\vskip 0.1cm

\hskip 0.3cm 2.1. Mahonian statistics on the set $M(\mu )$

\hskip 0.3cm 2.2. Dual mahonian statistics

\hskip 0.3cm 2.3. Generalized mahonian statistics

\vskip 0.2cm
\S 3. Main results
\vskip 0.1cm

\hskip 0.3cm 3.1. Combinatorial formula for modified Hall--Littlewood 
polynomials

\hskip 0.3cm 3.2. New combinatorial formula for the transition matrix 
$M(e,P)$

\vskip 0.2cm
\S 4. Proofs of Theorems~3.1 and 3.4
\vskip 0.1cm

\hskip 0.3cm 4.1. Proof of Theorem~3.4

\hskip 0.3cm 4.2. Proof of Theorem~3.1

\vskip 0.2cm
\S 5. Polynomials ${\cal P}_{\ld\mu}(t)$ and their interpretations

\vskip 0.2cm
\S 6. Generalizations of polynomials ${\cal P}_{\ld\mu}(t)$ and 
$K_{\ld\mu}(t)$
\vskip 0.1cm

\hskip 0.3cm 6.1. Crystal Kostka polynomials

\hskip 0.3cm 6.2. Fusion Kostka polynomials

\hskip 0.3cm 6.3. Ribbon Kostka polynomials

\hskip 0.3cm 6.4. Generalized Kostka polynomials

\hskip 0.3cm 6.5. Summary

\vskip 0.2cm
\S 7. Fermionic formulae
\vskip 0.1cm

\hskip 0.3cm 7.1. Multinomial fermionic formulae for one dimensional sums

\hskip 0.3cm 7.2. Rigged configurations polynomials

\vskip 0.2cm
\S 8. Two parameter deformation of one dimensional sums

\vskip 0.3cm
References
%}

\vfil\eject

{\bf \S 0. Introduction.}
\vskip 0.2cm

In this paper certain combinatorial and algebraic applications and 
generalizations of 
fermionic formulae for unrestricted one dimensional sums, obtained in 
[HKKOTY], are studied. The main applications of the fermionic formulae for 
unrestricted one dimensional sums considered in [HKKOTY] are related to 
the fermionic formulae for the branching functions and characters of some 
integrable representations of the affine Lie algebra $\wh{sl}(n)$.

Among applications considered in the present paper are the following ones:

\vskip 0.2cm
$\bullet$ New combinatorial formula for modified Hall--Littlewood 
polynomials $Q'_{\ld}(X_n;t)$ (Theorem~3.1). 
Let $\ld$ be a partition, $l(\ld )\le n$, then
$$Q'_{\ld}(X_n;t)=\sum_{\mu}{\cal P}_{\ld\mu}(t)m_{\mu}(X_n), \eqno (0.1)
$$
where $m_{\mu}(X_n)$ denotes the monomial symmetric function 
corresponding to partition $\mu$, $X_n=(x_1,\ldots ,x_n)$, and
$${\cal P}_{\ld\mu}(t)=\sum_{\{\nu \}}t^{c(\nu )}\prod_{k=1}^{n-1}
\prod_{i\ge 1}\left[\matrix{\nu_i^{(k+1)}-\nu_{i+1}^{(k)}\cr
\nu_i^{(k)}-\nu_{i+1}^{(k)}}\right]_t , \eqno (0.2)
$$
summed over all flags of partitions $\nu =\{ 
0=\nu^{(0)}\subset\nu^{(1)}\subset\cdots\subset\nu^{(n)}=\ld'\}$, such 
that $|\nu^{(k)}|=\mu_1+\cdots +\mu_k$, $1\le k\le n$, and
$c(\nu )=\ds\sum_{k=0}^{n-1}\sum_{i\ge 1}\pmatrix{\nu_i^{(k+1)}-
\nu_i^{(k)}\cr 2}$.

\vskip 0.2cm
$\bullet$ New combinatorial formula for the transition matrix $M(e,P)$ 
(Theorem~3.4). Let $\ld ,\mu$ be partitions, then
$$\eno{
&M(e,P)_{\ld\mu}=\sum_{\eta}K_{\eta\ld}K_{\eta'\mu}(t):=
{\cal R}_{\mu\ld}(t),~~~{\rm where}\cr
&{\cal R}_{\ld\mu}(t)=\sum_{\{\nu\}}\prod_{k=1}^{r-1}\prod_{i\ge 1}
\left[\matrix{\nu_i^{(k+1)}-\nu_{i+1}^{(k+1)}\cr \nu_i^{(k)}-
\nu_{i+1}^{(k+1)}}\right]_t, & (0.3)}
$$
summed over all flags of partitions $\{\nu\} =\{ 0=\nu^{(0)}\subset\nu^{(1)}
\subset\cdots\subset\nu^{(r)}=\ld'\}$ such that $\nu^{(k)}/\nu^{(k-1)}$ 
is a horizontal $\mu_k$--strip, $1\le k\le r$, $r=l(\mu )$.

\vskip 0.2cm
$\bullet$ New combinatorial formula for the number of rational points
${\cal F}_{\mu}^{\ld}({\bf F}_q)$ 
over the finite field ${\bf F}_q$ of the unipotent partial flag variety 
${\cal F}_{\mu}^{\ld}$ (Section~1.4):
$${\cal F}_{\mu}^{\ld}({\bf F}_q)=q^{n(\ld )}{\cal P}_{\ld\mu}(q^{-1}), 
\eqno (0.4)
$$
where polynomial ${\cal P}_{\ld\mu}(t)$ is given by (0.2).

\vskip 0.2cm
$\bullet$ New interpretation of the number $\al_{\ld}(S;p)$ of chains of 
subgroups
$$\{ e\}\subseteq H^{(1)}\subseteq\cdots\subseteq H^{(m)}\subseteq G
$$
in a finite abelian $p$--group of type $\ld$ such that each subgroup 
$H^{(i)}$ has order $p^{a_i}$, $1\le i\le m$ (Subsection~1.6):
$$p^{n(\ld )}\al_{\ld}(S;p^{-1})={\cal P}_{\ld\mu}(p). \eqno (0.5)
$$
Here $S=\{ a_1<a_2<\cdots <a_m\}$ is a subset of the set $[1,|\ld |-1]$, 
and $\mu :=\mu(S)=(a_1,a_2-a_1,\ldots ,a_m-a_{m-1},|\ld |-a_m)$.

\vskip 0.2cm
$\bullet$ New interpretation of the Schilling--Warnaar $t$--supernomial 
coefficients $\left[\matrix{{\bf L}\cr a}\right]_t$ and $T(L,a)$, [ScW], 
and Example~1, Subsection~3.1. Let $\ld =(\ld_1,\ldots ,\ld_k)$ be a 
partition, then
$$\left[\matrix{{\bf L}\cr a}\right]_t=\sum_{\eta}K_{\eta\mu}{\wt 
K}_{\eta\ld}(t)=t^{n(\ld )}{\cal P}_{\ld\mu}(t^{-1}), \eqno (0.6)
$$
where $L_i=\ld_i'-\ld_{i+1}'$, $1\le i\le k$, $\ld_{k+1}=0$, 
$\mu =\ds\left({1\over 2}|\ld |-a,~{1\over 2}|\ld |+a\right)$, and
$$\left[\matrix{{\bf L}\cr a}\right]_q=\sum_{j_1+\cdots +j_k=a+{|\ld 
|\over2}}t^{\ds\sum_{l=2}^kj_{l-1}(L_l+\cdots +L_k-j_l)}\left[\matrix{
L_k\cr j_k}\right]\left[\matrix{L_{k-1}+j_k\cr j_{k-1}}\right]\ldots
\left[\matrix{L_1+j_2\cr j_1}\right].
$$
$$T({\bf L},a):=t^{{1\over 4}L^tC_k^{-1}L-{a^2\over k}}\left[
\matrix{{\bf L}\cr a}\right]_{1/t}=t^{-B}{\cal P}_{\ld\mu}(t), 
\eqno (0.7)
$$
where $B=\ds{{1\over 2}n(\ld )+{1\over k}n(\mu')-{1\over 4k}\left(|\ld |^2
-(k+2)|\ld |\right) ={1\over 4}L^tC_k^{-1}L+{1\over k}\left( n(\mu')+
{|\mu |\over 2}\right)} $; 

\hskip -0.7cm$\left( C_k^{-1}\right)_{ij}:=\min (i,j)-
\ds{ij\over k}$, $1\le i,j\le k-1$, stands for the inverse of the Cartan 
matrix $C_k$ of the Lie algebra of type $A_{k-1}$.

We introduce also the $SU(n)$--analogue of $t$--multinomial coefficients 
(0.6) and (0.7) (Definition~3.2).

\vskip 0.2cm
$\bullet$ Definition and examples of the generalized mahonian statistics 
on the set of transport matrices ${\cal P}_{\ld\mu}$ (Section~2). This is 
a natural generalization of notion of {\it mahonian statistic} on the set 
of {\it words}
introduced and studied by D.~Foata in particular case $\ld =(1^n)$, [F],  
see also [Ma], [An], [ZB], [FZ], [GaW].

\vskip 0.2cm
$\bullet$ Connection between the rigged configurations polynomials 
$RC_{\ld R}(t)$ for a sequence of rectangular partitions $R=(R_1,\ldots 
,R_p)$, cf. [Ki1], and the classically restricted one dimensional sums 
$f_R^{\rm cl}(b_{T_{\max}};\mu )$ corresponding to the tensor product of 
"rectangular" crystals $B_{R_1}\otimes\cdots\otimes B_{R_p}$ (Section~7).

\vskip 0.2cm
$\bullet$ Definition, examples and properties of the two parameter 
deformation $B_{\ld\mu}(q,t)$ of the unrestricted one dimensional sum 
${\cal P}_{\ld\mu}(t)$ (Section~8).

\vskip 0.3cm
The paper is organized as follows:

In Section~1 we recall the definition of modified Hall--Littlewood 
polynomials, and explain a connection between the character of level 1 
basic representation of the affine Lie algebra $\wh{sl}(n)$, and the 
limit $N\to\infty$ of the modified Hall--Littlewood function 
corresponding to partition $\mu =(1^N)$, see [Ki2]. This result was 
extended to more general cases in [Ki2], [NY] and [HKKOTY]. 
In Subsections~1.4 and 1.5 we explain a connection between  the modified 
Hall--Littlewood polynomials and the unipotent partial flag varieties  
[HS], [LLT], [Sh], the Demazure characters [Ka2], [HKMOTY1,2], and the 
number of chains of subgroups in a finite abelian $p$--group, [Bu1,2,3], 
[F], [St].

In Section~2 we introduce the generalized mahonian and dual mahonian 
statistics on the set of transport matrices and on the set of 
(0,1)--transport matrices, respectively, and give few examples of
such statistics.

In Section~3 we state the fermionic formulae for polynomials ${\cal 
P}_{\ld\mu}(t)=\ds\sum_{\eta}K_{\eta\mu}K_{\eta\ld}(t)$ (Theorem~3.1) 
and ${\cal R}_{\ld\mu}(t)=\ds\sum_{\eta}K_{\eta\mu}K_{\eta'\ld}(t)$ 
(Theorem~3.4), and study their special cases. In particular, we show 
that polynomials ${\cal P}_{\ld\mu}(t)$ and 
$t^{n(\ld )}{\cal P}_{\ld\mu}(t^{-1})$ give a natural generalization
of generalized $p$--binomial coefficients introduced and studied by 
F.~Regonati [R], L.~Butler [Bu1], S.~Fishel [F], $\ldots$, and supernomial 
and multinomial coefficients introduced by A.~Schilling 
and S.O.~Warnaar, [Sc], [ScW], [W], and A.N.~Kirillov [Ki2]. 

In Section~4 we give algebraic proofs of main results, formulated in
Section~3, namely proofs of Theorems~3.1 and 3.4. A combinatorial 
proof of these theorems will appear elsewhere.

The main purpose of Section~5 is to show frequent apparitions of the one 
dimensional sums related to the tensor product of crystals 
$B_{R_1}\otimes\cdots\otimes B_{R_p}$ in different branches of 
Mathematics such as: representation theory, combinatorics, algebraic 
geometry, theory of finite abelian groups, and integrable systems. 
In our opinion, the fundamental role 
played by one dimensional sums in Mathematics and Mathematical Physics 
may be explained by the fact that one dimensional sums can be 
considered as a natural $q$--analog of the tensor product multiplicities. 
In the literature there exist at least 4 or 5 ways to define a 
$q$--analog of the Littlewood--Richardson numbers, see, e.g., [GoW], 
[BKMW], [CL], [LLT], [KLLT], [KS], $\ldots$.
 
In Section~6 we overview several known ways to define the $q$--analogues 
of the tensor product multiplicities, and formulate conjectures 
(Conjectures~6.4, 6.5, 6.8 and 6.9) 
which relate the classically restricted one dimensional sums $f_R^{\rm 
cl}(b_{T_{\min}};\mu ):=CK_{\mu R}(t)$, the ribbon Kostka polynomials 
$K_{\ld\mu}^{(p)}(t)$, introduced by A.~Lascoux, B.~Leclerc and 
J.-Y.~Thibon, [LLT], and the generalized Kostka polynomials $K_{\ld 
R}(t)$, introduced by M.~Shimozono and J.~Weyman. We expect that only in 
the case of dominant sequence of rectangular partitions $R$ the crystal 
$CK_{\mu R}(t)$, the ribbon $K_{\Lambda (R)\mu}^{(p)}(t)$ and the 
generalized Kostka polynomials $K_{\mu R}(t)$ give the equivalent 
$q$--analogues of the tensor product multiplicities. We also formulate 
some unsolved problems.

In Section~7 we formulate few conjectures about multinomial fermionic
formulae for homogeneous unrestricted one dimensional sums, and 
generalized Kostka--Foulkes polynomials corresponding to a sequence of
rectangles.

In Section~8 we suggest two parameter deformations of polynomials 
${\cal P}_{\ld\mu}(t)$ and one dimensional sums.

\vskip 0.3cm
{\bf Acknowledgments.} This paper presents an extended version of my talk 
"Fermionic formulae for the branching functions of the affine Lie algebra 
$\wh{sl}(n)$" delivered at PhD Centennial Conference, Department of 
Mathematics, University of Wisconsin--Madison, May~22--24, 1997. I would 
like to thank Richard Askey and Georgia Benkart for invitation to this 
conference. I am thankful to Goro Hatayama, Atsuo Kuniba, Franklin 
M.~Maley, Masato Okado, 
Mark Shimozono, Taichiro Takagi, Jean--Yves Thibon and Yasuhiko Yamada for 
very fruitful discussions and suggestions. I wish to thank L.~Vinet for 
hospitality at the CRM, University of Montreal, where this work was 
completed.

\vskip 0.5cm
%\vfil\eject
{\bf \S 1. Modified Hall-Littlewood polynomials.}
\vskip 0.3cm

{\bf 1.1.} Definition.
\vskip 0.2cm

Let $\ld$ be a partition, $l(\ld )\le n$, $Q_{\ld}(X_n;q)$ and 
$P_{\ld}(X_n;q)$ be the Hall--Littlewood polynomials corresponding to 
$\ld$, see e.g. [M], Chapter~III.

\vskip 0.2cm
{\bf Definition 1.1.} {\it A modified Hall--Littlewood polynomial 
$Q'_{\ld}(X_n;q)$ is defined to be
$$Q'_{\ld}(X_n;q)=Q_{\ld}(X_n/(1-q);q):=Q_{\ld}(X'_n;q), \eqno (1.1)
$$
where the variables $X'_n$ are the products $xq^{j-1}$, $j\ge 1$, 
$x\in X_n:=(x_1,\ldots ,x_n)$.}
\vskip 0.2cm

%{\bf Remark.} We will use also for short an abbriviation 
%"mHL--polynomials" instead of "modified Hall--Littlewood polynomials".
%\vskip 0.2cm 

The $Q'_{\ld}(X;q)$ serve to interpolate between the Schur functions 
$s_{\ld}$ and the complete homogeneous symmetric functions $h_{\ld}$, 
because
$$Q'_{\ld}(X;0)=s_{\ld}(X)
$$
as it is clear from (1.1), and
$$Q'_{\ld}(X;1)=h_{\ld}(X)
$$
as it is clear from Cauchy's identity (1.3) below.

\vskip 0.2cm
{\bf Proposition 1.2.} 
$$Q_{\ld}'(X;q)=\sum_{\mu}\left(\sum_{\eta}K_{\eta\mu}K_{\eta\ld}(q)
\right) m_{\mu}(X)=\sum_{\eta}s_{\eta}(X)K_{\eta\ld}(q).\eqno (1.2)
$$
\vskip 0.2cm

{\it Proof.} Let us remind that the Hall-Littlewood polynomials $Q_{\ld}$ 
and $P_{\ld}$ satisfy the following orthogonality condition (see, e.g. 
[M], Chapter~III, (4.4))
$$\sum_{\ld}Q_{\ld}(X;q)P_{\ld}(Y;q)=\prod_{x\in X, \ y\in Y}
{1-qxy\over 1-xy}.
$$
Hence, (cf. [M], Example~7a on p.234)
$$\eno{\sum_{\ld}Q'_{\ld}(X;q)P_{\ld}(Y;q)&=\prod_{k\ge 0}
\prod_{x\in X, \ y\in Y}{1-q^{k+1}xy\over 1-q^kxy}=\prod_{x\in X, \ y\in Y}
(1-xy)^{-1}\cr \cr
&=\sum_{\ld}s_{\ld}(x)s_{\ld}(y). & (1.3)}
$$
\vbox{Here we have used the Cauchy identity for Schur functions, [M], Chapter~I, 
(4.3). It remains to remind the definition of the Kostka--Foulkes 
polynomials:
$$s_{\ld}(Y)=\sum_{\mu}K_{\ld\mu}(q)P_{\mu}(Y;q).\eqno (1.4)
$$
\qed}

{\bf 1.2.} Modified Hall--Littlewood polynomials for partition $\ld =(1^N)$.
\vskip 0.2cm

{\bf Corollary 1.3.} {\it Let $\mu =(\mu_1,\ldots ,\mu_n)$ be a 
composition, $|\mu |=N$. Then
$$\sum_{\ld}K_{\ld\mu}K_{\ld (1^N)}(q)=q^{n(\mu')}\left[\matrix{N\cr
\mu_1,\ldots ,\mu_n}\right]_q,\eqno (1.5)
$$
where $\ds\left[\matrix{N\cr \mu_1,\ldots ,\mu_n}\right]_q={(q;q)_N\over
(q;q)_{\mu_1}\ldots (q;q)_{\mu_n}}$ is the $q$--analog of gaussian 
multinomial coefficient, and $(a;q)_n:=\ds\prod_{j=0}^{n-1}(1-aq^j)$.}

\vskip 0.2cm
{\it Proof.} First of all we have to compute $Q'_{(1^N)}(X;q)$. For this 
goal, let us remark that ([M], Chapter~III, (2.8))
$$Q_{(1^N)}(X;q)=(q;q)_Ne_N(X), \eqno (1.6)
$$
where $e_m(X)$ is the elementary symmetric function of degree $m$ in the 
variables $X$. Hence
$$\sum_{\ld}s_{\ld}(X)K_{\ld (1^N)}(q)=Q_{(1^N)}(X/((1-q);q)=
(q;q)_Ne_N(X/(1-q)).
$$
Now let us put $E(X)=\ds\sum_{N\ge 0}e_N(X)t^N=\prod_{x\in X}(1+tx)$. 
Then we have
$$E\left( X/(1-q)\right) =\prod_{x\in X}(-tx;q)_{\infty},
$$
and using the Euler identity
$$(x;q)_{\infty}=\sum_{n=0}^{\infty}{x^nq^{n(n-1)\over 2}\over (q;q)_n},
\eqno (1.7)
$$
we obtain the following result
$$E\left( X/(1-q)\right) =\sum_{N\ge 0}\left(\sum_{\mu\vdash N}
{q^{n(\mu')}\over (q;q)_{\mu}}m_{\mu}(X)\right)t^N, \eqno (1.8)
$$
where for a composition $\mu =(\mu_1,\mu_2,\ldots ,\mu_n)$ we set
$(q;q)_{\mu}:=\ds\prod_{j=1}^n(q;q)_{\mu_j}$. Finally, from (1.6) and (1.8) we 
obtain immediately that
$$Q'_{(1^N)}(X_n;q)=\sum_{\ld}s_{\ld}(X_n)K_{\ld (1^N)}(q)=
\sum_{\mu\vdash N}q^{n(\mu')}
\left[\matrix{N\cr \mu_1,\ldots ,\mu_n}\right] m_{\mu}(X_n). \eqno (1.9)
$$
\qed

%\vskip 0.2cm
%\vfil\eject
{\bf 1.3.} Hall-Littlewood polynomials and characters of the affine Lie 
algebra $\wh{sl}(n)$.
\vskip 0.2cm

We consider the identity (1.9) as the finitization of the 
Weyl--Kac--Peterson character formula (WKR--formula for short, see, e.g. 
[Kac], (12.7.12)) for the level 1 
basic representation $L(\Lambda_0)$ of the affine Lie algebra 
$\wh{sl}(n)$. Indeed, the WKP--formula for the character ch$L(\Lambda_0)$ 
may be recovered as an appropriate limit of (1.9). More exactly, let us 
consider the following form of (1.9):
$$\eno{
&q^{-{(N^2-N)n\over 2}}\sum_{\ld}{s_{\ld}(x_1,\ldots ,x_n)\over
(x_1\ldots x_n)^N}K_{\ld (1^{nN})}(q) = \cr 
&~~\sum_{\matrix{k\in{\bf Z}^n\cr |k|=0,\ k_i\ge -N,\ \forall i}}
x_1^{k_1}\cdots x_n^{k_n}q^{{1\over 2}k_i^2}\left[\matrix{nN\cr
k_1+N,\ldots ,k_n+N}\right]_q. & (1.10)} 
$$
First of all, $\lim_{N\to\infty}{\rm RHS}(1.10)=$
$$
{1\over (q;q)^{n-1}_{\infty}}\sum_{m\in{\bf Z}^n, \ |m|=0}x^mq^{{1\over 2}
\sum m_i^2}={1\over (q;q)^{n-1}_{\infty}}\sum_{k\in{\bf Z}^{n-1}}
z_1^{k_1}\cdots z_{n-1}^{k_n-1}q^{{1\over 2}kA_{n-1}k^t},
$$
where $z_i=\ds{x_i\over x_{i-1}}$, $1\le i\le n-1$, $x_0:=x_n$. On the 
other hand,
$$\lim_{N\to\infty}{\rm LHS}(1.10)=\sum_{\ld 
=(\ld_1\ge\ld_2\ge\cdots\ge\ld_n)\in{\bf Z}^n, \ |\ld |=0}
s_{\ld}(x_1,\ldots ,x_n)b_{\ld}(q), \eqno (1.11)
$$
where $b_{\ld}(q)$ is defined to be
$$b_{\ld}(q):=\lim_{N\to\infty}q^{-{n(N^2-N)\over 
2}}K_{\ld_N(1^{|\ld_N|})}(q)={q^{n(\ld')}\over (q;q)_{\infty}^{n-1}}
\prod_{1\le i\le j\le n}(1-q^{\ld_i-\ld_j-i+j}), \eqno (1.12)
$$
and for a given weight $\ld$ we set $\ld_N:=\ld +(N^n)$.
The last equality in (1.12) follows from the hook--formula (see, e.g. [M], 
Example~2 on p.243):
$$K_{\ld (1^N)}(q)={q^{n(\ld')}(q;q)_N\over\prod_{x\in\ld}(1-q^{h(x)})},
$$
where $h(x):=\ld_i+\ld_j'-i-j+1$ is the hook--length corresponding to the 
box $x=(i,j)\in\ld$.

Finally, it follows from (1.10)--(1.12) that
$$\sum_{\ld}s_{\ld}(x_1,\ldots ,x_n)b_{\ld}(q)={\Theta (x)\over 
(q;q)_{\infty}^{n-1}}, \eqno (1.13)
$$
summed over all partitions $\ld$ such that $l(\ld )\le n$, and $|\ld 
|\equiv 0(\mod n)$, and
where 
$$\Theta (x)=\ds\sum_{m=(m_1,\ldots ,m_n)\in{\bf Z}^n, \ |m|=0}
x^mq^{{1\over 2}(m_1^2+\cdots +m_n^2)}
$$ 
is the theta--function 
corresponding to the basic representation $L(\Lambda_0)$ of $\wh{sl}(n)$, 
[Kac], \S 12.7. It is well-known that the RHS(1.13) is equal to the character 
of the level 1 basic representation $L(\Lambda_0)$. Hence, $b_{\ld}(q)$ coincide 
with the branching functions ([Kac], \S 12.2) for the level 1 basic 
representation $L(\Lambda_0)$ of the affine Lie algebra $\wh{sl}(n)$ (cf. [Ki2]).

For further results concerning a connection between modified 
Hall--Littlewood functions and characters, and branching functions of the 
affine Lie algebra $\wh{sl}(n)$, see [Ki2], [NY], and [HKKOTY].
\vskip 0.3cm

{\bf 1.4.} Modified Hall--Littlewood polynomials and unipotent flag 
varieties.
\vskip 0.2cm

Polynomials ${\cal P}_{\ld \mu}(q):=\ds\sum_{\eta}K_{\eta\mu}K_{\eta\ld}(q)$ 
have the following geometric interpretation due to [HS] and [Sh]. Let 
$V$ be an $n$--dimensional vector space over an  algebraically closed 
field $k$, and let $\mu$, $l(\mu )=r$, be a composition of $n$. A 
$\mu$--flag in $V$ is a sequence $F=\{ V_1,\ldots ,V_r\}$ of subspaces of 
$V$ such that $V_1\subset V_2\subset\cdots\subset V_r=V$, and 
dim$V_i=\mu_1+\cdots +\mu_i$, $1\le i\le r$. Let ${\cal F}_{\mu}$ denote 
the set of all $\mu$--flags in $V$. The group $G:=GL(V)$ acts 
transitively on ${\cal F}_{\mu}$, so that ${\cal F}_{\mu}$ may be 
identified with $G/P$, where $P$ is the subgroup which fixes a given 
flag, and therefore ${\cal F}_{\mu}$ is a non--singular projective 
algebraic variety, the partial flag variety of $V$. 

Now let $u\in G$ be a unipotent endomorphism of $V$ of type $\ld$, so 
that $\ld$ is a partition of $n$ which describes the Jordan canonical 
form of $u$, and let ${\cal F}_{\mu}^{\ld}\subset{\cal F}_{\mu}$ be the 
set of all $\mu$--flags $F\in{\cal F}_{\mu}$ fixed by $u$. The set ${\cal 
F}_{\mu}^{\ld}$ is a closed subvariety of ${\cal F}_{\mu}$.

It has been shown by N.~Shimomura ([Sh], see also [HS]), that

$\bullet$ if $k={\bf C}$ is the field of complex numbers, the variety 
${\cal F}_{\mu}^{\ld}$ %/{\bf C}$ 
admits a cell decomposition, involving only cells of even real dimensions, 
and
$$t^{2n(\ld )}{\cal P}_{\ld\mu}(t^{-2}):=
\sum_{\eta}K_{\eta\mu}\wt K_{\eta\ld}(t^2)=
\sum_it^{2i}{\rm dim}H_{2i}({\cal F}_{\mu}^{\ld},{\bf Z}) \eqno (1.14)
$$
is the Poincare polynomial of ${\cal F}_{\mu}^{\ld}/{\bf C}$, where
$\wt K_{\eta\ld}(t):=t^{n(\ld )}K_{\eta\ld}(t^{-1})$;

$\bullet$ if $k$ contains the finite field of $q$ elements, ${\bf F}_q$, 
the number ${\cal F}_{\mu}^{\ld}(q)$ of ${\bf F}_q$--rational 
points of ${\cal F}_{\mu}^{\ld}$ is equal to 
$q^{n(\ld )}{\cal P}_{\ld\mu}(q^{-1})$.

\vskip 0.3cm
{\bf 1.5.} Modified Hall--Littlewood polynomials and Demazure characters.
\vskip 0.2cm

Let $\g$ be a symmetrizable Kac--Moody algebra. Recall that for every 
dominant integral weight $\Lambda$, there exists a unique (up to 
isomorphism) irreducible module $V=V(\Lambda )$ of highest weight 
$\Lambda$. The character of $V$, denoted by ${\rm ch}V$, is the formal sum
$${\rm ch}V=\sum_{\Lambda'}\left(\dim V_{\Lambda'}\right)e^{\Lambda'},
$$
summed over all weights $\Lambda'$, where $V_{\Lambda'}$ is the weight 
subspace of $V$ of weight $\Lambda'$, and where $e^{\Lambda'}$ is a 
formal exponential. This sum makes sense because each $V_{\Lambda'}$ is 
finite--dimensional. For definitions and further details, see [Kac].

Let $\b$ be the Borel subalgebra of $\g$  and let $w$ be an element of 
the Weyl group $W$. The $\b$--module generated by the one dimensional 
extremal weight subspace $V_{w(\Lambda )}$ is denoted by $V_w(\Lambda )$ 
and called a Demazure module. They are finite--dimensional subspaces 
which form a filtration of $V$ which is compatible with the Bruhat order 
of $W$, i.e. $V_w(\Lambda )\subseteq V_{w'}(\Lambda )$ whenever $w\le w'$ with 
respect to the Bruhat order, $w,w'\in W$, and $\ds\bigcup_{w\in 
W}V_w(\Lambda )=V(\Lambda )$, see, e.g., [Ka2].

From now let us assume that $\g=\wh{sl}(n)$. Let $\Lambda_i$ and $r_i$, 
$0\le i\le n-1$, denote the fundamental weight and simple reflection with 
respect to the simple root $\al_i$, of $\wh{sl}(n)$. It is convenient to 
define $\Lambda_i$ and $r_i$ for all $i\in{\bf Z}$ using the agreement 
$\Lambda_i=\Lambda_{i+n}$, $r_i=r_{i+n}$.

Now we are ready to explain an interpretation of the modified 
Hall--Littlewood polynomial $Q'_{(l^L)}(X_n;q)$ corresponding to a 
rectangular partition $(l^L)$ as the character of certain Demazure's 
module. This result is due to [KMOTU2]:

Let $L\ge 1$ be an integer, and $w:=w_{L,n}=r_{Ln-1}r_{Ln-2}\cdots 
r_{L+2}r_{L+1}r_L$ be an element of the affine Weyl group 
$W(A_{n-1}^{(1)})$ of type $A_{n-1}^{(1)}$. Then
$$e^{-l\Lambda_0}{\rm ch}V_w(l\Lambda_L)=q^{-E_0}Q'_{(l^L)}(X_n;q), \eqno (1.15)
$$
where $E_0=l\ds\left[{L\over n}\right]\left( L-{n\over 2}\left(\left[{L\over 
n}\right] +1\right)\right)$.

\vskip 0.3cm
{\bf 1.6.} Modified Hall--Littlewood polynomials and chains of subgroups 
in a finite abelian 

\hskip 0.8cm $p$--group.
\vskip 0.2cm

Let $p$ be a prime number. It is well--known (see, e.g. [H])
that any abelian group $G$ 
of order $p^n$ is isomorphic to a direct product of cyclic groups
$$G\approx{\bf Z}/p^{\ld_1}{\bf Z}\times\cdots\times{\bf Z}/p^{\ld_l}{\bf 
Z}
$$
where $\ld_1\ge\ld_2\ge\cdots\ge\ld_l>0$, $\ld_1+\cdots +\ld_l=n$. 
The partition $\ld$ is called the type of $G$.

For any partition $\nu\subseteq\ld$, let us denote by $\al_{\ld}(\nu ;p)$ 
the number of subgroups of type 
$\nu$ in a finite abelian $p$--group of type $\ld$. 

More generally, for any flag of partitions $\{\nu\} =\{\nu^{(1)}\subseteq
\nu^{(2)}\subseteq\cdots\subseteq\nu^{(m)}\subset\ld\}$ 
denote by $\al_{\ld}(\nu^{(1)},\ldots ,\nu^{(m)};p)$ (or 
$\al_{\ld}(\{\nu\} ;p)$ for short) the number of chains of subgroups
$$\{ e\}\subseteq H^{(1)}\subseteq H^{(2)}\subseteq\cdots\subseteq 
H^{(m)}\subseteq G
$$
in a finite abelian $p$--group $G$ of type $\ld$ such that the type of 
$H^{(i)}$ is $\nu^{(i)}$.

The problem of counting the number of subgroups of type $\nu$ in a finite 
abelian $p$--group of type $\ld$ has a long history and goes back at 
least to the beginning of 1900's, see e.g., papers by G.A.~Miller [Mi] and 
by H.~Hiller [Hi]. In 1934 G.~Birkhoff [Bi] has discovered an interesting 
connection between the set of subgroups of finite abelian $p$--group and 
that of so--called standard matrices of G.~Birkhoff. In 1948 three 
mathematicians, S.~Delsarte [De], P.~Dyubyuk [Dy], and 
Yenchien Yeh [Y] published formulae for the number $\al_{\ld}(\nu ;p)$ of 
subgroups of type $\nu$ in a finite abelian $p$--group of type $\ld$:
$$\al_{\ld}(\nu ;p)=\prod_{j\ge 1}p^{\nu_{j+1}'(\ld_j'-\nu_j')}
\left[\matrix{\ld_j'-\nu_{j+1}'\cr \nu_j'-\nu_{j+1}'}\right]_p,\eqno(1.16)
$$
where $\ld'$ is the conjugate of $\ld$, and $\nu'$ is the conjugate of 
$\nu$.

In order to explain a connection between the numbers $\al_{\ld}(\nu ;p)$ 
and $\al_{\ld}(\{\nu\} ;p)$ and unrestricted one dimensional sums ${\cal 
P}_{\ld\mu}(t)$, it is convenient to introduce the following  polynomials 
$p^{n(\ld )}\al_{\ld}(\nu ;p^{-1})$ and $p^{n(\ld )}\al_{\ld}(\{\nu\} 
;p^{-1})$.

\vskip 0.2cm
{\bf Proposition 1.4.} {\it i) For any partitions $\nu\subseteq\ld$,
$$p^{n(\ld )}\al_{\ld}(\nu ;p^{-1})=\ds\prod_{j\ge 
1}p^{\pmatrix{\ld_j'-\nu_j'\cr 2}+\pmatrix{\nu_j'\cr 2}}
\left[\matrix{\ld_j'-\nu_{j+1}'\cr 
\nu_j'-\nu_{j+1}'}\right]_p. \eqno (1.17)
$$
ii) Let $\{\nu\} =\{ 
0=\nu^{(0)}\subset\nu^{(1)}\subset\cdots\subset\nu^{(m)}\subset\nu^{(m+1)}
=\ld\}$ be a flag of partitions. Then
$$p^{n(\ld )}\al_{\ld}(\{\nu\} ;p^{-1})=p^{c(\nu )}\prod_{i=1}^m\prod_{j\ge 
1}\left[\matrix{(\nu^{(i+1)})_j'-(\nu^{(i)})_{j+1}'\cr 
(\nu^{(i)})_j'-(\nu^{(i)})_{j+1}'}\right]_p, \eqno(1.18)
$$ 
where $c(\nu )=\ds\sum_{i=0}^m\sum_{j\ge 
1}\pmatrix{(\nu^{(i+1)})_j'-(\nu^{(i)})_j'\cr 2}$.}
\vskip 0.2cm

Proofs of (1.17) and (1.18) easily follow from the formula (1.16).

\vskip 0.2cm
{\bf Definition 1.5} (see, e.g., [Bu1]). {\it Let $p$ be a prime number, 
$\ld$ be a partition of $n$, 
and $S=\{ 1\le a_1<\cdots <a_m<n\}$ be a subset of $[1,n-1]$. Let us denote 
by $\al_{\ld}(S;p)$ the number of chains of subgroups
$$\{ e\}\subseteq H^{(1)}\subseteq\cdots\subseteq H^{(m)}\subseteq G
$$
in a finite abelian $p$--group $G$ of type $\ld$, where each subgroup 
$H^{(i)}$ has order $p^{a_i}$.}
\vskip 0.2cm

It follows from Definition~1.5, that
$$\al_{\ld}(S;p)=\sum_{\{\nu\}}\al_{\ld}(\{\nu\} ;p), \eqno (1.19)
$$
summed over all flags of partitions $\{\nu\} =\{ 
0=\nu^{(0)}\subset\nu^{(1)}\subset\cdots\subset\nu^{(m)}\subset\nu^{(m+1)}
=\ld\}$ such that $|\nu^{(k)}|=a_k$, $1\le k\le m$.

\vskip 0.2cm
{\bf Proposition 1.6.} {\it Let $\ld$ be a partition, and $S=\{ 1\le 
a_1<\cdots <a_m<|\ld |\}$ be a subset of $[1,|\ld |-1]$. Then
$$p^{n(\ld )}\al_{\ld}(S;p^{-1})={\cal P}_{\ld\mu}(p), \eqno(1.20)
$$
where $\mu :=\mu (S)$ stands for the composition $\mu =(a_1,a_2-a_1,\ldots 
,a_m-a_{m-1},|\ld |-a_m)$.}
\vskip 0.2cm

Proof follows easily from (1.19) and (0.2).

\qed

\vskip 0.2cm
{\bf Corollary 1.7.} {\it Let $\ld$ and $S$ be as in Proposition~1.6. Then
$$\al_{\ld}(S;p)=\sum_{\eta}K_{\eta\mu (S)}{\wt K}_{\eta\ld}(p), 
\eqno(1.21)
$$
where $\mu (S)=(a_1,a_2-a_1,\ldots ,a_m-a_{m-1},|\ld |-a_m)$, and 
${\wt K}_{\eta\ld}(p)=p^{n(\ld )}K_{\eta\ld}(p^{-1})$.}
\vskip 0.2cm

Below we will give few examples of application of the formula (1.21).

$\bullet$ Follow to [R], [Bu1], [Fi], let us define the generalized 
$p$--binomial coefficient $\left[\matrix{\ld'\cr k}\right]$ to be the 
number of subgroups of order $p^k$ of a finite abelian group of type 
$\ld$. If $\ld =(1^n)$, then $\ld'=(n)$, and $\left[\matrix{\ld'\cr 
k}\right]$ coincides with $p$--binomial coefficient $\left[\matrix{n\cr 
k}\right]_p$. Note also that
$$\left[\matrix{\ld'\cr k}\right] =\sum_{\nu\vdash k}\al_{\ld}(\nu 
;p)=\sum_{\nu\vdash k}p^{{\overline c}(\nu )}\prod_{j\ge 1}\left[
\matrix{\ld'_j-\nu_{j+1}'\cr \nu_j'-\nu_{j+1}'}\right]_p, \eqno(1.22)
$$
where ${\overline c}(\nu )=\ds\sum_{j\ge 1}\nu_{j+1}'(\ld_j'-\nu_j')$.

It follows from formulae (1.22) and (3.3) that
$$\left[\matrix{\ld'\cr k}\right] =\left[\matrix{{\bf L}\cr a}\right]_p,
$$
where $L_i=\ld_i'-\ld_{i+1}'$, $1\le i\le l(\ld')$, ${\bf L}=(L_i)$, 
$a=\ds{n\over 2}-k$, and $\left[\matrix{{\bf L}\cr a}\right]_t$ stands 
for the Schilling--Warnaar $t$--supernomial coefficient (0.6).

It is easy to see from Corollary~1.7, that
$$\left[\matrix{\ld'\cr k}\right] =\al_{\ld}(\{ k\} 
;p)=\sum_{\eta}K_{\eta\mu (k)}{\wt K}_{\eta\ld}(p),
$$
where $\mu (k)=(k,n-k)$. On the other hand, it is clear that the Kostka--Foulkes 
number \hbox{$K_{\eta ,(k,n-k)}=0$} unless $\eta =(\eta_1,\eta_2)$ and $\eta_1\ge 
k$; in the later case $K_{\eta ,(k,n-k)}=1$. Thus,
$$\left[\matrix{\ld'\cr k}\right] =\sum_{l\ge k}{\wt K}_{(l,n-l),\ld}(p),
$$
and
$$\left[\matrix{\ld'\cr k}\right] -\left[\matrix{\ld'\cr k-1}\right] 
={\wt K}_{(k,n-k),\ld}(p).
$$
This result is due to Lynne Butler [Bu3].

$\bullet$ Let $\ld$ be a partition, $|\ld |=n$, and $S=\{ a_1,\ldots 
,a_m\}$ be a subset of $[1,n-1]$. Follow [St], [Bu1], [Bu2], consider the 
following polynomial
$$\beta_{\ld}(S;p)=\sum_{T\subseteq S}(-1)^{|S-T|}\al_{\ld}(T;p).
$$
It is known, [St], [Bu1], [Bu2], that $\beta_{\ld}(S;p)$ is equal to the top 
(and only non--vanishing) Betti number of a certain simplicial complex 
$\Delta_{\ld}(S)=\Delta_{\ld}(S;p)$; we refer the reader to [St] for 
definition of the simplicial complex $\Delta_{\ld}(S)$ and further details.

Let us show that polynomial $\beta_{\ld}(S;p)$ has nonnegative 
coefficients. Indeed,
$$\beta_{\ld}(S;p)=\sum_{T\subseteq S}(-1)^{|S-T|}\al_{\ld}(T;p)
=\sum_{\eta}\left(\sum_{T\subseteq S}(-1)^{|S-T|}K_{\eta\mu 
(T)}\right){\wt K}_{\eta\ld}(p).
$$

Our nearest aim is to show that the number $\left(\ds\sum_{T\subseteq 
S}(-1)^{|S-T|}K_{\eta\mu (T)}\right)\ge 0$. For this goal let us show 
that the latter number counts the number of Littlewood--Richardson tableaux 
of a certain skew shape $b(S)$ and weight $\eta$. More precisely, for a 
given subset $S=\{ a_1<a_2<\cdots <a_m\}$ of the set $[1,n-1]$, the skew 
shape $b(S)$ is the {\it border strip} with $a_1$ squares in row 1, $a_2-a_1$ 
squares in row $2,\ldots ,n-a_m$ squares in row $m+1$. For the reader's 
convenience, let us remind that

$\bullet$ a skew shape is called a border strip if consecutive rows 
overlap by exactly one square (see, e.g., [M], p.5);

$\bullet$ a skew tableau $T$ is called Littlewood--Richardson tableau, 
if the word $w(T)$ corresponding to the tableau $T$ is a lattice 
permutation (see, e.g. [M]. Chapter~I, \S 9).

\vskip 0.2cm
{\bf Proposition 1.8} (R.~Stanley [St]). {\it Let $S=\{ a_1<\cdots <a_m\}$ be 
a subset of the set $[1,n-1]$, and $b(S)$ stands for border strip with 
$a_1$, squares in row 1, $a_2-a_1$ squares in row $2,\ldots ,n-a_m$ 
squares in row $m+1$. Then
$$\beta_{\ld}(S;p)={\wt K}_{b(S),\ld}(p),
$$
where ${\wt K}_{b(S),\ld}(p)$ stands for the cocharge Kostka--Foulkes 
polynomial corresponding to the skew shape $b(S)$, see, e.g., [Ki1], 
[Bu1].}
\vskip 0.2cm

Indeed, 
$$\eno{
\beta_{\ld}(S;p)=\sum_{\eta}\left(\sum_{T\subseteq 
S}(-1)^{|S-T|}K_{\eta\mu (T)}\right){\wt K}_{\eta\ld}(p)
&=\sum_{\eta}\#|{\rm Tab}^0(b(S),\eta )|{\wt K}_{\eta\ld}(p)\cr
&={\wt K}_{b(S),\ld}(p).}
$$
To deduce the second equality we used the following formula 
$$\sum_{T\subseteq S}(-1)^{|S-T|}K_{\eta\mu (T)}=\#|{\rm Tab}^0(b(S),\eta 
)|, \eqno (1.23)
$$
where ${\rm Tab}^0(b(S),\eta )$ stands for the set of all 
Littlewood--Richardson tableaux of skew shape $b(S)$ and weight $\eta$. 
The formula (1.23) can be obtained using the results from [KKN].

\qed

Finally, let us describe (see [Bu1], Definition~1.3.1) the statistic {\it 
value}, denoted by $v$ on the set of tabloids. This statistic generates 
the generalized mahonian statistic VAL on the set of transport matrices, 
see Section~2.

\vskip 0.2cm
{\bf Definition 1.9} ([Bu1]). {\it Let $T$ be a tabloid of any shape and 
weight, and $x\in T$ be an entry of $T$. Then the {\it value} $v(x)$ of 
the entry $x$ in $T$ is the number of smaller entries in the same column 
and above $x$, or in the next column to the right and below $x$. The 
value $v(T)$ of $T$ is the sum of the values of the entries in 
$T:v(T)=\ds\sum_{x\in T}v(x)$.}
\vskip 0.2cm

{\bf Example.} Consider
$$T= \hbox{
          \normalbaselines\m@th\offinterlineskip
          \vtop{\hbox{\Hfourbox(1,2,1,2)}
          \vskip-0.4pt
          \hbox{\Hfourbox(2,2,1,3)}
          \vskip-0.4pt
          \hbox{{\Hthreebox(3,2,2)}}
          \vskip-0.4pt
          \hbox{\Hthreebox(3,2,3)}}}\ \  
\in \ T(4433,\ 374),
$$
Then $v(T)=(0+1+3+2)+(1+0+0+0)+(0+0+2+3)+(0+1)=13$. Note that the 
Shimomura statistic $d(T)$ (see Section~2) of the tabloid $T$ is equal to 
10, see Example in Subsection~2.3.

\vskip 0.2cm
{\bf Proposition 1.10} ([Bu1]). {\it Let $\ld$ be a partition, $|\ld |=n$, 
and $S=\{ a_1<a_2<\cdots <a_m\}$ be a subset of $[1,n-1]$. Then
$$\al_{\ld}(S;p)=\sum_Tp^{v(T)},
$$
summed over all tabloids $T$ of shape $\ld$ and weight $\mu :=\mu (S)$.}
%=(a_1,a_2-a_1,\ldots ,a_m-a_{m-1},n-a_{m-1})$.}

%\vfil\eject
\vskip 0.5cm
{\bf \S 2. Generalized mahonian statistics.}
\vskip 0.3cm

{\bf 2.1.} Mahonian statistics on the set $M(\mu )$.
\vskip 0.2cm

We start with recalling the definition of mahonian statistic on words, 
[F].  A word is a finite sequence of letters, $w=w_1\ldots w_N$, where 
each letter is in the set $\{ 1,\ldots ,n\}$. Let 
$\mu =(\mu_1,\ldots ,\mu_n)$ be a composition, $|\mu |=N$, denote 
by $M(\mu )$ the set of all words $w=w_1\ldots w_N$ of weight $\mu$, i.e. 
$\mu_i$ is the number of occurrences of $i$ in the word $w$.
It is well--known (see, e.g. [An]) that the cardinality of the set $M(\mu 
)$ is equal to the multinomial coefficient $\ds\pmatrix{N\cr 
\mu_1,\ldots ,\mu_n}$.

\vskip 0.2cm
{\bf Definition 2.1} ([F]). {\it A function $\varphi$ on the set $M(\mu )$ is 
called mahonian statistic, if}
$$\sum_{w\in M(\mu )} q^{\varphi (w)}=\left[\matrix{N\cr \mu_1,\ldots 
,\mu_n}\right]_q.
$$
\vskip 0.2cm

{\bf Examples.} $1^0$ (Inversion number, [Ma]). Let $w\in M(\mu )$ be a word, 
define the number of inversions for the word $w$ to be
$$INV(w)=\sum_{1\le i<j\le N}\chi (w_i>w_j).
$$
where $\chi (A)=1$ if $A$ is true and 0 otherwise.

$2^0$ (Major index, [Ma]). Define the major index of the word $w$ to be
$$MAJ(w)=\sum_{m=1}^{N-1}m\chi (w_m>w_{m+1}).
$$

$3^0$ (Modified major index, [Ki2]). Define the modified major index of the word 
$w$ to be
$$\wt{MAJ}(w)=\sum_{m=1}^{N-1}m\chi (w_m\ge w_{m+1})-n(\mu'), \ \ {\rm 
where} \ \ n(\mu')=\sum_{i=1}^n\pmatrix{\mu_i\cr 2}.
$$

$4^0$ (Zeilberger's index, [ZB]). For given $w\in M(\mu )$ let $w_{ij}$ be the 
subword of $w$ formed by deleting all letters $w_m$ such that $w_m\ne i$ 
or $j$. For example, if $w=2411213144321\in M(5323)$, then 
$w_{12}=21121121$, $w_{13}=1113131$, $w_{14}=41111441$, $w_{23}=22332$, 
$w_{24}=242442$, $w_{34}=43443$.

$\bullet$ Zeilberger's index, or $Z$--index, of a word $w$ is defined to be 
the sum of %the 
major indices of all 2--letter subwords $w_{ij}$ of $w$:
$$Z(w)=\sum_{1\le i<j\le n}MAJ(w_{ij}).
$$

\vbox{$\bullet$ Modified Zeilberger's index, or $\wt Z$--index, of a word $w$ is 
defined to be the sum of 
modified major indices of all 2--letter subwords $w_{ij}$ of $w$:
$$\wt Z(w)=\sum_{1\le i<j\le n}\wt{MAJ}(w_{ij}).
$$

Next example will require some definitions. First, for a word $w$ 
let denote $\overline w$  the non--decreasing rearrangement of the 
letters of $w$. Second, if $a$ and $b$ are positive integers, with $a\le 
n$, let
$$C[a,b]=\cases{[a+1,a+2,\ldots ,b],& if $a\le b$;\cr 
[1,2,\ldots ,b,a+1,a+2,\ldots ,n],& if $a>b$.}
$$}

$5^0$ (Denert's index, M.~Denert, see e.g., [FZ], [GaW]). 
Define the Denert index of a word $w$ to be
$$DEN(w)=\sum_{1\le i<j\le n}\chi (w_i\in C[w_j,\overline w_j]).
$$
For example, if $w=M(5323)$ as above, then $INV(w)=29$, $MAJ(w)=47$, 
$\wt{MAJ}(w)=42$, $Z(w)=46$, $\wt Z(w)=31$, $DEN(w)=46$.

\vskip 0.2cm
{\bf Theorem 2.2} ([Ma], [ZB], [Ki2], [FZ]). {\it The statistics INV, MAJ, 
$\wt{MAJ}$, $Z$, $\wt Z$, DEN are mahonian.}

\vskip 0.3cm
{\bf 2.2.} Dual mahonian statistics.
\vskip 0.2cm

Now we are going to extend the notion of mahonian statistic to the 
set of transport matrices. Let us denote by ${\cal 
P}_{\ld\mu}$ (respectively ${\cal R}_{\ld\mu}$) the set of all matrices 
of non--negative integers (respectively the set of all (0,1)--matrices) 
with row sums 
$\ld_i$ and column sums $\mu_j$. It is clear that if $\ld =(1^N)$ then 
the both sets ${\cal P}_{(1^N)\mu}$ and ${\cal R}_{(1^N)\mu}$ can be naturally 
identified with the set $M(\mu )$.

\vskip 0.2cm
{\bf Definition 2.3.} {\it A function $\psi$ on the set ${\cal 
R}_{\ld\mu}$ is called dual mahonian statistic if}
$$\sum_{m\in{\cal R}_{\ld\mu}}q^{\psi (m)}=\sum_{\eta}K_{\eta\mu}
K_{\eta'\ld}(q).
$$
\vskip 0.2cm

Let us give a few examples of the dual mahonian statistics.

$1^0$. The first example is due to A.~Zelevinsky, see [M], Chapter~III, \S 
6, Example~5, p.244. Let $\ld$ and $\mu$ be compositions of the same 
integer $n$, and $m=(m_{ij})\in{\cal R}_{\ld\mu}$. For each element 
$a=m_{ij}$ of the matrix $m$ we denote by $i(a):=i$, and $j(a):=j$ its 
first and second coordinates. We denote by supp$(m)=\{ m_{ij}\in 
m~|~m_{ij}\ne 0\}$ the set of all nonzero entries of $m$. If 
$a=m_{ij}\in{\rm supp}(m)$ we define the height of $a$ to be ${\rm 
ht}(a)=\ds\sum_{1\le k\le i}m_{kj}$. For each $a\in{\rm supp}(m)$ let us define
$$i^+(a)=\cases{i(b),& if $\exists b\in{\rm supp}(m)$ such that 
$j(a)=j(b)$ and ${\rm ht}(b)={\rm ht}(a)+1$;\cr +\infty ,& if such $b$ 
doesn't exist.}
$$
Follow A.~Zelevinsky $[ibid]$, for each $m\in{\cal R}_{\ld\mu}$ we 
define $\wt{ZEL}(m)=\ds\sum_{a\in{\rm supp}(m)}\wt z(a)$, where $\wt 
z(a)$ is equal to the number of $b\in{\rm supp}(m)$ such that 

$i)$ $j(b)<j(a)$, 

$ii)$ ht$(b)={\rm ht}(a)$, 

$iii)$ $i(a)<i(b)<i^+(a)$.

\vskip 0.2cm
{\bf Theorem 2.4} (A.~Zelevinsky).
$${\cal R}_{\ld\mu}(q)=\sum_{m\in{\cal R}_{\ld\mu}}q^{\wt{\rm ZEL}(m)},
$$
{\it in other words, the statistic $\wt{\rm ZEL}$ is dual mahonian.}
\vskip 0.2cm

There exists a bijection between the set ${\cal R}_{\ld\mu}$ and that of 
all column strict tabloids of shape $\ld'$ and weight $\mu$. Let $\nu$ 
and $\mu$ be 
compositions of the same integer $n$. A tabloid of shape $\nu$ and weight 
$\mu$ is a filling of the diagram of boxes with row lengths 
$\nu_1,\nu_2,\ldots ,\nu_r$, such that the number $i$ 
occurs $\mu_i$ times, 
and such that each column is nondecreasing. A tabloid of shape $\nu$ and 
weight $\mu$ is called a column strict if each column is strictly 
decreasing. For example, 
$$\matrix{\hbox{
          \normalbaselines\m@th\offinterlineskip
          \vtop{\hbox{\Htwobox(1,2)}
          \vskip-0.4pt
          \hbox{{\Hthreebox(1,3,1)}}
          \vskip-0.4pt
          \hbox{\Fsquare(0.4cm,2)}
          \vskip-0.4pt
          \hbox{\Htwobox(4,3)}}}
&& {\rm and} && \hbox{
          \normalbaselines\m@th\offinterlineskip
          \vtop{\hbox{\Htwobox(1,1)}
          \vskip-0.4pt
          \hbox{{\Hthreebox(2,2,1)}}
          \vskip-0.4pt
          \hbox{\Fsquare(0.4cm,3)}
          \vskip-0.4pt
          \hbox{\Htwobox(4,3)}}}}
$$
\vskip 0.2cm%
\hskip -0.7cm are tabloid and column strict tabloid of weight (3221) and shape 
(2312). We denote by $T(\nu ,\mu )$ (respectively, $\wt T(\nu ,\mu )$) 
the set of all tabloids (respectively, the set of all column strict 
tabloids) of shape $\nu$ and weight $\mu$. 

Now we are ready to describe a 
bijection ${\cal R}_{\ld\mu}\leftrightarrow\wt T(\ld',\mu )$ in the case 
when $\ld$ is a partition. Namely, consider a matrix 
$m=(m_{ij})\in{\cal R}_{\ld\mu}$. 
Let us fill the shape $\ld'$ by positive integers according to the 
following rule: if $m_{ij}\ne 0$, put the number $i$ in the box of the shape $\ld'$ with coordinates 
$(i,j)$. As a result we obtain the tabloid $T$ of shape $\ld'$ and weight 
$\mu$. For example, consider $\ld =(3221)$, $\mu =(2231)$, and
$$m=\left\{\matrix{1&0&1&0\cr 0&1&1&0\cr 1&1&0&1\cr 1&0&0&0}\right\}
\in{\cal R}_{\ld ,\mu}.
$$
The corresponding column strict tabloid is $\hbox{
          \normalbaselines\m@th\offinterlineskip
          \vtop{\hbox{\Hfourbox(1,2,1,3)}
          \vskip-0.4pt
          \hbox{{\Hthreebox(3,3,2)}}
          \vskip-0.4pt
          \hbox{\Fsquare(0.4cm,4)}}}$.

\vskip 0.3cm          
It is easy to see that the correspondence $m\to T$ defines a bijection. Let 
us continue and define a statistic $\wt{ZEL}$ on the set of all column 
strict tabloids. Namely, for any column strict tabloid $T$ of shape 
$\nu$, let $\wt{ZEL}(T)$ denote the number of pairs $(x,y)\in\nu\times\nu$ 
such that $y$ lies to the left from $x$ (in the same row) and 
$T(x)<T(y)<T(x\!\!\downarrow )$. We have used here the following notations: 
if a box 
$x$ has coordinates $(i,j)$, then $x\!\!\downarrow =(i+1,j)$ and 
$\vec{x}=(i,i+1)$; for any 
$x\in\nu$, $T(x)$ is the integer located in the box $x$ of the tabloid 
$T$; if $x\!\!\downarrow$ does not belong to the shape $\nu$, then we put
$T(x\!\!\downarrow )=+\infty$. It is clear that if a composition $\nu$ 
contains only one part, then $\wt T(\nu ,\mu )=M(\mu )$, and $\wt{ZEL}$ 
coincides with statistic $INV$.

\vskip 0.2cm
{\bf Theorem 2.5}  (A.~Zelevinsky).
$${\cal R}_{\ld\mu}(q)=\sum_{T\in\wt T(\ld',\mu )}q^{\wt{ZEL}(T)}.
$$

$2^0$. Let $\ld$ be a partition and $\mu$ be a composition of the same 
integer $n$, and $m\in{\cal R}_{\ld\mu}$. There is an explicit 
one--to--one correspondence, due to Knuth [Kn], between the set of 
(0,1)--matrices with row sums $\ld_i$ and column sums $\mu_j$, and pairs 
of semistandard tableaux of conjugate shapes and weights $\ld$, $\mu$, 
(Knuth's dual correspondence):
$$\matrix{{\cal R}_{\ld\mu}&\cong&\coprod_{\eta}{\rm SST}(\eta ,\mu )
\times{\rm SST}(\eta',\ld )\cr \cr m&\leftrightarrow &(P,Q).}
$$
Let us define the charge $CH$ of a matrix $m\in{\cal R}_{\ld\mu}$ to be 
the Lascoux--Sch\"utzenberger charge ([LS]) of the corresponding 
semistandard tableaux $Q$ of weight $\ld$:
$$CH(m)=c(Q).
$$
It follows from the results of Lascoux and Sch\"utzenberger [LS], and Knuth 
[Kn], that
$${\cal R}_{\ld\mu}(q)=\sum_{m\in{\cal R}_{\ld\mu}}q^{{CH}(m)}.
$$

It is an interesting problem to find a bijective proof that the 
statistics $\wt{ZEL}$ and $CH$ have the same distribution on the set 
${\cal R}_{\ld\mu}$.
\vskip 0.3cm

{\bf 2.3.} Generalized mahonian statistics.

\vskip 0.3cm
{\bf Definition 2.6.} {\it A function $\varphi$ on the set of transport 
matrices ${\cal P}_{\ld\mu}$ is called generalized mahonian statistic if
$$\sum_{m\in{\cal P}_{\ld\mu}}q^{\varphi (m)}=q^{E_0}\sum_{\eta}
K_{\eta\mu}K_{\eta\ld}(q),
$$
for a certain constant $E_0:=E_{0,\varphi}$.}
\vskip 0.2cm

There is a well--known bijection between sets ${\cal P}_{\nu\mu}$ and 
$T(\nu, \mu)$. To describe this bijection, let $m\in{\cal P}_{\nu\mu}$, 
and $D(\nu )$ be the diagram of the composition ${\bf\nu}$. To obtain a tabloid, 
let us fill the first $m_{1j}$ boxes of the $j$-th row of $D(\nu )$ by the 
number 1, the next $m_{2j}$ boxes of same row by the number 2, and so on. 
As a result we obtain the tabloid of shape $\nu$ and weight $\mu$. This 
construction defines the bijection under consideration. To go further, 
let us recall the Shimomura cells decomposition [Sh] of the fixed point 
variety ${\cal F}_{\mu}^{\ld}$ of a unipotent $u$ of type $\ld$ ($\ld$ is 
a partition) acting on the partial flag variety ${\cal F}_{\mu}$. The 
cells in Shimomura's decomposition are indexed by tabloids of shape $\ld$ 
and weight $\mu$. The dimension $d(T)$ of the cell $c_T$ indexed by 
$T\in T(\ld ,\mu )$ is computed by algorithm described below ([Sh], 
[LLT]), and defines the mahonian statistic $\wt d(T)=n(\ld )-d(T)$ 
(=codimension of the cell $c_T$) on the set ${\cal P}_{\ld\mu}$:
$${\cal P}_{\ld\mu}(t)=\sum_{T\in T(\ld ,\mu )}t^{\wt d(T)}.
$$

The dimensions $d(T)$ are given by the following algorithm ([LLT], 
Section~8.1).

1) If $T\in T(\ld ,(n))$ then $d(T)=0$;

2) If $\mu =(\mu_1,\mu_2)$ has exactly two parts, and $T\in T(\ld ,\mu 
)$, then $d(T)$ is computed as follows. A box $x$ of $T$ is called 
special if $T(x)$ is the lowest 1 of the column containing $x$. For a box 
$y$ such that $T(y)=1$, put $d(y)=0$; if $T(y)=2$, set $d(y)$ equals to 
the number of nonspecial 1's lying in the row of $y$, plus the number of 
special 1's lying in the same row, but from the right side of $y$. Then 
$d(T)=\sum d(y)$, summed over all $y\in T$ such that $T(y)=2$.

3) Let $\mu =(\mu_1,\ldots ,\mu_k)$ and $\mu^*=(\mu_1,\ldots ,\mu_{k-1})$. 
For $T\in T(\ld ,\mu )$, let $T_1$ be the tabloid obtained from $T$ by changing 
the entries $k$ into 2 and all the other ones by 1. Let $T_2$ be the 
tabloid  of weight $\mu^*$ obtained from $T$ by erasing all the entries $k$,
and rearranging the 
columns in the appropriate order. Then $d(T)=d(T_1)+d(T_2)$.

\vskip 0.2cm
{\bf Example.} Consider
$$T= \hbox{
          \normalbaselines\m@th\offinterlineskip
          \vtop{\hbox{\Hfourbox(1,2,1,2)}
          \vskip-0.4pt
          \hbox{\Hfourbox(2,2,1,3)}
          \vskip-0.4pt
          \hbox{{\Hthreebox(3,2,2)}}
          \vskip-0.4pt
          \hbox{\Hthreebox(3,2,3)}}}\ \  
\in \ T(4433,\ 374),
$$
then
$$ T_1= \hbox{
          \normalbaselines\m@th\offinterlineskip
          \vtop{\hbox{\Hfourbox(1,1,1,{\bf 1})}
          \vskip-0.4pt
          \hbox{\Hfourbox(1,1,1,2)}
          \vskip-0.4pt
          \hbox{{\Hthreebox(2,1,{\bf 1})}}
          \vskip-0.4pt
          \hbox{\Hthreebox(2,{\bf 1},2)}}}  \ \ \ \ \ 
T_2= \hbox{
          \normalbaselines\m@th\offinterlineskip
          \vtop{\hbox{\Hfourbox(2,1,{\bf 1},2)}
          \vskip-0.4pt
          \hbox{\Hthreebox(2,{\bf 1},2)}
          \vskip-0.4pt
          \hbox{{\Htwobox(2,2)}}
          \vskip-0.4pt
          \hbox{\Fsquare(0.4cm,2)}}}
$$
\vskip 0.3cm
\hskip -0.7cm where the special entries are printed in bold type. Thus, 
$d(T)=d(T_1)+d(T_2)=\hbox{(3+2+1)}+(2+1+1)=10$.
\vskip 0.2cm

There is a variant of this construction due to A.~Lascoux, B.~Leclerc and 
J.-Y.~Thibon [LLT], in which the shape $\ld$ is allowed to be an arbitrary 
composition. Such a variant has already been used by I.~Terada [T] in the 
case of complete flags (i.e. $\mu =(1^N))$. 

Let $\nu$ be a 
composition, and $T\in T(\nu ,\mu )$. Follow to [LLT], define an integer 
$e(T)$ by the following rules:

$i)$ for $T\in T(\nu ,(n))$, $e(T)=d(T)=0$;

$ii)$ for $T\in T(\nu ,(\mu_1,\mu_2))$, $e(T)=d(T)$;

$iii)$ otherwise $e(T)=e(T_1)+e(T_2)$ where $T_1$ is defined as above, but 
this time $T_2$ is obtained from $T$ by erasing the entries $k$, without 
reordering.

Let us define $\wt e(T)=n(\ld )-e(T)$.

\vskip 0.2cm
{\bf Proposition 2.7} ([LLT]). Let $\nu$ be a composition and $\ld =\nu^+$ be 
the corresponding partition. Then
$$\sum_{T\in T(\ld ,\mu )}t^{\wt d(T)}=\sum_{T\in T(\nu ,\mu )}t^{\wt e(T)}
={\cal P}_{\ld \mu}(t).
$$

\vskip 0.2cm
{\bf Example} (cf. [LLT], Example~8.4). Take $\ld =(321)$, $\mu =(42)$ 
and $\nu =(312)$. The set $T(\ld ,\mu )$ consists of the following tabloids
\vskip 0.1cm
$$\matrix{ T 
&& \hbox{
          \normalbaselines\m@th\offinterlineskip
          \vtop{\hbox{\Hthreebox(1,1,2)}
          \vskip-0.4pt
          \hbox{{\Htwobox(1,2)}}
          \vskip-0.4pt
          \hbox{\Fsquare(0.4cm,1)}}}
&& \hbox{
          \normalbaselines\m@th\offinterlineskip
          \vtop{\hbox{\Hthreebox(1,2,1)}
          \vskip-0.4pt
          \hbox{{\Htwobox(1,2)}}
          \vskip-0.4pt
          \hbox{\Fsquare(0.4cm,1)}}}
&& \hbox{
          \normalbaselines\m@th\offinterlineskip
          \vtop{\hbox{\Hthreebox(1,1,2)}
          \vskip-0.4pt
          \hbox{{\Htwobox(1,1)}}
          \vskip-0.4pt
          \hbox{\Fsquare(0.4cm,2)}}}                                
&& \hbox{
          \normalbaselines\m@th\offinterlineskip
          \vtop{\hbox{\Hthreebox(1,1,1)}
          \vskip-0.4pt
          \hbox{{\Htwobox(1,2)}}
          \vskip-0.4pt
          \hbox{\Fsquare(0.4cm,2)}}}
&& \hbox{
          \normalbaselines\m@th\offinterlineskip
          \vtop{\hbox{\Hthreebox(1,1,1)}
          \vskip-0.4pt
          \hbox{{\Htwobox(2,1)}}
          \vskip-0.4pt
          \hbox{\Fsquare(0.4cm,2)}}}\cr \cr \cr
\wt d(T) && 2~~~~~ && 1~~~~~ && 2~~~~ && 4~~~~ && 3~~~~}
$$
\vskip 0.2cm

The set of tabloids of shape $\nu$ and weights $\mu$ contains the 
following ones
\vskip 0.1cm
$$\matrix{T&&\hbox{
          \normalbaselines\m@th\offinterlineskip
          \vtop{\hbox{\Hthreebox(1,2,1)}
          \vskip-0.4pt
          \hbox{{\Fsquare(0.4cm,1)}\hskip 0.39cm\Fsquare(0.4cm,2)}
          \vskip-0.4pt
          \hbox{\Fsquare(0.4cm,1)}}}
&&\hbox{
          \normalbaselines\m@th\offinterlineskip
          \vtop{\hbox{\Hthreebox(1,1,2)}
          \vskip-0.4pt
          \hbox{{\Fsquare(0.4cm,1)}\hskip 0.39cm\Fsquare(0.4cm,2)}
          \vskip-0.4pt
          \hbox{\Fsquare(0.4cm,1)}}}
&&\hbox{
          \normalbaselines\m@th\offinterlineskip
          \vtop{\hbox{\Hthreebox(1,2,1)}
          \vskip-0.4pt
          \hbox{{\Fsquare(0.4cm,1)}\hskip 0.39cm\Fsquare(0.4cm,1)}
          \vskip-0.4pt
          \hbox{\Fsquare(0.4cm,2)}}}
&&\hbox{
          \normalbaselines\m@th\offinterlineskip
          \vtop{\hbox{\Hthreebox(1,1,1)}
          \vskip-0.4pt
          \hbox{{\Fsquare(0.4cm,1)}\hskip 0.39cm\Fsquare(0.4cm,2)}
          \vskip-0.4pt
          \hbox{\Fsquare(0.4cm,2)}}}
&&\hbox{
          \normalbaselines\m@th\offinterlineskip
          \vtop{\hbox{\Hthreebox(1,1,1)}
          \vskip-0.4pt
          \hbox{{\Fsquare(0.4cm,2)}\hskip 0.39cm\Fsquare(0.4cm,1)}
          \vskip-0.4pt
          \hbox{\Fsquare(0.4cm,2)}}}\cr \cr \cr
\wt e(T)&&1~~~~&&2~~~~&&2~~~~&&4~~~~&&3~~~~}
$$

Let us give yet another example of generalized mahonian statistic 
denoted by VAL. This example is due essentially to Lynne Butler [Bu1].

Let $\ld$ be a partition and $\mu$ be a composition, $|\ld |=|\mu |$. On 
the set $T(\ld ,\mu )$ of tabloids of shape $\ld$ and weight $\mu$ one 
can define the statistic {\it value} $v$, see [Bu1], Definition~1.3.1, or 
Subsection~1.6,

\vskip 0.2cm
{\bf Definition 2.8.} {\it Let us define ${\rm VAL}(T)=n(\ld )-v(T)$.}
\vskip 0.2cm

{\bf Example.} Take $\ld =(321)$ and $\mu =(42)$, Consider the set of 
tabloids $T(\ld ,\mu)$ in the same order as in the previous Example. Then 
the values of statistic VAL on the set $T(\ld ,\mu)$ are the following 
3,4,2,1,2, and
$$\sum_{T\in T(\ld ,\mu )}t^{{\rm VAL}(T)}={\cal P}_{\ld\mu}(t).
$$

\vskip 0.2cm
{\bf Proposition 2.9} ([Bu1]). {\it Let $\ld$ be a partition and $\mu$ be 
a composition. Then}
$$\sum_{T\in T(\ld ,\mu)}t^{{\rm VAL}(T)}={\cal P}_{\ld\mu}(t).
$$
\vskip 0.2cm

{\bf Problem 1.} Find a bijective proof that if $\ld$ is a partition, 
then the Shimomura statistic $\wt d$, LLT--statistic $\wt e$, statistic 
VAL, and the energy function $E$ are equidistribute on the set of 
transport matrices ${\cal P}_{\ld\mu}$.          

\vskip 0.5cm
{\bf \S 3. Main results.}
\vskip 0.2cm

{\bf 3.1.} Combinatorial formula for modified Hall-Littlewood polynomials.
\vskip 0.3cm

{\bf Theorem 3.1.} ([HKKOTY]) {\it Let $\ld$ be a partition, and $\mu$, 
$l(\mu )=r$, be a composition of the same integer $n$, then}
 
$${\cal P}_{\ld\mu}(t):=\sum_{\eta}K_{\eta\mu}K_{\eta\ld}(t)=
\sum_{\{\nu\}}t^{c(\nu )}\prod_{k=1}^{r-1}\prod_{i\ge 1}
\left[\matrix{\nu_i^{(k+1)}-\nu_{i+1}^{(k)}\cr 
\nu_i^{(k)}-\nu_{i+1}^{(k)}}\right]_t, \eqno (3.1)
$$
{\it summed over all flags of partitions $\nu =\{ 
0=\nu^{(0)}\subset\nu^{(1)}\subset\cdots\subset\nu^{(r)}=\ld'\}$, such 
that $|\nu^{(k)}|=\mu_1+\cdots +\mu_k$, $1\le k\le r$; and
$$c(\nu )=\sum_{k=0}^{r-1}\sum_{i\ge 
1}\left(\matrix{\nu_i^{(k+1)}-\nu_i^{(k)}\cr 2}\right) .
$$
where for any real number $\al$ we put $\pmatrix{\al\cr 2}:={\al (\al 
-1)}{2}$.}

Proof of Theorem 3.1 will be given in Subsection~4.2.
\vskip 0.2cm

{\bf Remark.} It is well--known ([Kn]; [M], Chapter~I, Section~6) that 
${\cal P}_{\ld\mu}(1)$ is equal to the number of matrices of 
non--negative integers with row sums $\ld_i$ and column sums $\mu_j$. 
This number is equal also to that of pairs of semistandard tableaux of 
the same shape and weights $\ld$ and $\mu$, [Kn].

\vskip 0.2cm
{\bf Examples.} $1^0$. Let us take a length two composition $\mu=(\mu_1,\mu_2)$, 
and a partition $\ld$. Let $\ld'=(\ld_1',\ldots ,\ld_k')$ be the 
conjugate partition. Then the identity (3.1) takes the following form
$$\sum_{\eta}K_{\eta\mu}K_{\eta\ld}(t)=\sum_{\nu\vdash\mu_1}
t^{c(\nu )}\prod_{i=1}^k\left[\matrix{\ld_i'-\nu_{i+1}\cr \nu_i-\nu_{i+1}}
\right]_t, \eqno (3.2)
$$
summed over all partitions $\nu$ of $\mu_1$, $l(\nu )=k$, and $c(\nu )=
\ds\sum_{i=1}^k\pmatrix{\ld_i'-\nu_i\cr 2}+
\sum_{i=1}^k\pmatrix{\nu_i\cr 2}$.

Let us put $L_i=\ld_i'-\ld_{i+1}'$ and $j_i=\ld_i'-\nu_i$, $1\le i\le k$, 
$j_{k+1}=0$. Then we have 
$$\ds\sum_{i=1}^kj_i=\mu_2,~~~~
\left[\matrix{\ld_i'-\nu_{i+1}\cr \nu_i-\nu_{i+1}}\right]_t=
\left[\matrix{L_i+j_{i+1}\cr j_i}\right]_t,
$$
and
$$\eno{
c(\nu )&=\sum_{i=1}^k\pmatrix{\ld_i'\cr 2}+\sum_{i=1}^kj_i(j_i-\ld_i')\cr
&=\sum_{i=1}^kj_i(j_i-L_i-L_{i+1}-\cdots -L_k)+{1\over 2}\sum_{1\le i,j\le k}
\min (i,j)L_iL_j-{1\over 2}|\mu |.}
$$
Thus, RHS(3.2)$=t^A\left[\matrix{{\bf L}\cr a}\right]_{1/t}$,
%$=t^B{\wt T}({\bf L},\mu_2)$, 
where
 
\hskip -0.7cm$A=\ds{1\over 2}\sum_{1\le i,j\le k}\min 
(i,j)L_iL_j-{1\over 2}|\mu |$, %$B=A-(L_1+\cdots +L_k)\mu_2$,
and $a=-\ds{\mu_1-\mu_2\over 2}$;
 
\hskip -0.7cm$\left[\matrix{{\bf L}\cr a}\right]_t$,  
%and ${\wt T}({\bf L}, a)$ are 
stands for the Schilling--Warnaar $t$--supernomial 
coefficients (0.5), see [ScW], (2.9). %and (3.2).
\vskip 0.2cm

It follows from the formulae above that
$$\left[\matrix{{\bf L}\cr a}\right]_t=
\sum_{\eta}K_{\eta\mu}{\wt K}_{\eta\ld}(t), \eqno (3.3)
$$
where~~ $\wt K_{\eta\ld}(t)=t^{n(\ld )}K_{\eta\ld}(t^{-1})$,~~~ 
$\mu=\ds\left({1\over 2}\left(\sum_{i=1}^kiL_i\right)-a,~ 
{1\over 2}\left(\sum_{i=1}^kiL_i\right)
+a\right)$, 

\hskip -0.7cm and ~~$\ld_i'=L_i+\cdots +L_k$, $1\le i\le k$.

Now let us assume additionally that $\ld_1'=\cdots =\ld_k'=N$, or 
equivalently, $\ld =(k^N)$. Then the RHS(3.2) can be rewritten in the 
following form
$$\sum_{\nu\vdash\mu_1}t^{c(\nu )}\left[\matrix{N\cr 
N-\nu_1,\nu_1-\nu_2,\ldots ,\nu_{k-1}-\nu_k,\nu_k}\right]_t, \eqno (3.4)
$$
where $c(\nu )=k\pmatrix{N\cr 2}-N\mu_1+
\ds\sum_{i=1}^k\nu_i^2$. The sum in (3.4) is taken over all 
partitions $\nu$ of $\mu_1$ such that $l(\nu )=k$.
%and ($N=m_1+\cdots +m_k$)
%$$\left[\matrix{N\cr m_1,\ldots ,m_k}\right]_q:={(q;q)_N\over 
%(q;q)_{m_1}\cdots (q;q)_{m_k}},
%$$
%is the $q$--multinomial coefficient.

If we put $m_i=k(\nu_1+\cdots +\nu_{k-i})-(k-i)\mu_1$, $1\le i\le k-1$, 
(and, consequently, $k\nu_i=m_{k-i}-m_{k-i+1}+\mu_1$, $1\le i\le k-1$), 
then the sum (3.4) coincides with the RHS(2.49), [Ki2], Theorem~14 (in 
[Ki2] we have used $q$ instead of $t$). 

Let us remark that the sum (3.4) is closely related to the special 
value $p=0$ of the Schilling and Warnaar 
$q$--multinomial coefficient $\left[\matrix{L\cr a}\right]_k^{(p)}$ 
([Sc], \S 2, and [W], 
Definition~1). More precisely, we state that sum 
(3.4) is equal to $t^{k\pmatrix{N\cr 2}}\left[\matrix{N\cr 
\mu_1}\right]_k^{(0)}(t^{-1})$. This statement is equivalent (cf. (3.3)) 
to the following one:
$$\left[\matrix{N\cr\mu_1}\right]_k^{(0)}=\sum_{\eta}K_{\eta\mu}
{\wt K}_{\eta ,(k^N)}(q). \eqno (3.5)
$$
Formulae (3.3) and (3.5) suggest the following definition:
\vskip 0.2cm

{\bf Definition 3.2.} {\it Let $\ld$ be a partition and $\mu$ be a 
composition, $|\ld |=|\mu |$. Define the $t$--multinomial coefficient 
$\left[\matrix{\ld\cr\mu}\right]^{(0)}$ to be}
$$\left[\matrix{\ld\cr\mu}\right]^{(0)}=\sum_{\eta}K_{\eta\mu}
{\wt K}_{\eta\ld}(t).
$$
\vskip 0.2cm

Thus, see Corollary~1.7, if $t=p$ is a prime number and $l(\mu )=m+1$, then 
the $t$--multinomial coefficient $\left[\matrix{\ld\cr\mu}\right]^{(0)}$ 
counts the number of chains of subgroups
$$\{ e\}\subseteq H^{(1)}\subseteq H^{(2)}\subseteq\cdots\subseteq 
H^{(m)}\subseteq G
$$
of a finite abelian $p$--group $G$ of type $\ld$ such that each subgroup 
$H^{(i)}$ has order $p^{\mu_1+\cdots +\mu_i}$.

It follows from (3.3) that if a composition $\mu =(\mu_1,\mu_2)$ consists 
of two parts then the $t$--multinomial coefficient 
$\left[\matrix{\ld\cr\mu}\right]^{(0)}$ coincides with the 
$t$--supernomial coefficient $\left[\matrix{{\bf L}\cr a}\right]_t$, 
where $a=-\ds{\mu_1-\mu_2\over 2}$, and if $\ld'=(\ld_1',\ldots ,\ld_k')$
is the conjugate partition, 
then ${\bf L}:=(L_1,\ldots ,L_k)$ with $L_i=\ld_i'-\ld_{i+1}'$, $1\le 
i\le k$.

More generally, let $B$ be a crystal (see, e.g., [Ka1], [KMOTU1], 
[HKKOTY]), and $b\in B$. Define the (unrestricted) $t$--multinomial 
coefficient $T^{(b)}(\ld ;\mu )$ to be 
$$T^{(b)}(\ld ;\mu )=t^{-E_{\min}}\sum_{p\in{\cal P}_{\mu}(b,\ld )}t^{E(p)},
$$
where ${\cal P}_{\mu}(b,\ld )$ is the set of paths $p=b\otimes 
b_1\otimes\cdots\otimes b_m\in B\otimes B_{(\mu_1)}\otimes\cdots\otimes 
B_{(\mu_m)}$ such that $wt(b_1)+\cdots +wt(b_m)=\ld$; $E(p)$ is the 
energy of a path $p$ (see, e.g., [HKKOTY]).

Similarly, one can define classically restricted and restricted 
$t$--multinomial coefficients. We intend to consider the properties 
(including recurrence relations, bosonic formulae, multinomial analogue 
of Bailey's lemma, and applications to polynomial identities and 
$q$--series) of these $t$--multinomial coefficients in a separate 
publication.

$2^0$. If $\mu =(1^n)$, then (3.1) coincides with 
%a particular case $t=0$ of 
the formula for modified Green's polynomials $X^{\ld}_{(1^n)}(t)$
from [M], Example 4 on p.249. 

Let us describe two generalized mahonian 
statistics on the set $M(\ld )$. The first one is the 
Lascoux--Sch\"utzenberger charge $c$ defined on the set of dominant 
weight words $w$, i.e. $w\in M(\ld )$, where $\ld$ is a partition, see 
[LS]; [M], Chapter~III, \S 6, p.242. The second one is the $LP$ statistic 
(see, e.g., [GaW]) which can be defined for arbitrary words. 

\vskip 0.2cm
{\bf Definition 3.3.} {\it Let 
$w$ be a word, define $lp_i(w)$ to be the number of distinct letters 
to the left of position $i$ and having the same multiplicity as the 
letter in position $i$ in the truncated word $w_1\ldots w_i$. Let 
$LP(w)=\ds\sum_{i\ge 2}lp_i(w)$.}
\vskip 0.2cm

 For example, 
$LP(3422231413)=0+2+0+0+1+1+0+2+0=6$. One can show that if $\nu$ is a 
composition, $\ld =\nu^+$ is the corresponding partition, then
$${\cal P}_{\ld (1^n)}(q)=\sum_{w\in M(\nu )}q^{LP(w)}=
\sum_{w\in M(\ld )}q^{c(w)}.
$$

$3^0$. If $\ld =(1^N)$, then the RHS(3.1) coincides with that of (1.5).

$4^0$. Let $\mu$ be a composition of length $n$, and $\ld =(2^N)$, so 
that $|\mu |=2N$. In this case we have

$\bullet$ $\nu^{(k)}=(\nu_1^{(k)},\nu_2^{(k)})$, 
$|\nu^{(k)}|=\mu_1+\cdots +\mu_k$, $1\le k\le n$;

$\bullet$ 
$0\le\nu_2^{(k-1)}\le\nu_2^{(k)}\le\nu_1^{(k)}\le\nu_1^{(k+1)}\le N$, if 
$1\le k\le n-1$, and $\nu^{(n)}=(N,N)$.

If we define $m_i=2\nu_1^{(i)}-\mu_1-\cdots -\mu_i\ge 0$, $0\le i\le 
n-1$, and
$$\beta_i={\mu_i+m_{i-1}-m_i\over 2}\in{\bf Z}_{\ge 0}, \ \ 
1\le i\le n, \ \ m_0=m_n=0,
$$
then the RHS(3.1) takes the following form
$$\sum_{m\in{\bf Z}_{\ge 0}^{n-1}}t^{c(m)}\left[\matrix{N\cr \beta_1,
\ldots ,\beta_n}\right]_t\prod_{k=1}^{n-1}\left[\matrix{\beta_{k+1}+
m_{k+1}\cr m_k}\right]_t, \eqno (3.6)
$$
summed over all sequences $m=(m_1,\ldots ,m_{n-1})\in{\bf Z}_{\ge 
0}^{n-1}$ such that $m_i+m_{i-1}+\mu_i\equiv 0$(\mod 2), $1\le i\le n$, 
$m_0=m_n=0$, and $c(m)=\ds\sum_{i=1}^n\pmatrix{\mu_i\cr 2}+{1\over 
4}mC_{n-1}m^t$, where $C_{n-1}$ is the Cartan matrix of type $A_{n-1}$. 

It is well--known that the LHS(3.1) does not depend on the permutations of 
components of the composition $\mu$. Hence,  the same is valid for the 
RHS(3.1) as well. This is not obvious at all because the number of terms 
in the 
right hand side sum (3.1) do depends on the composition $\mu$, but not 
only on the corresponding partition $\mu^+$. For example, let us take 
$\mu =(1221)$ and $\ld =(2^3)$. The summands in the RHS(3.1) correspond to 
the following flags of partitions $\nu =\{\nu^{(1)}\subset\nu^{(2)}
\subset\nu^{(3)}\subset\nu^{(4)}\}$:
\vskip 0.2cm
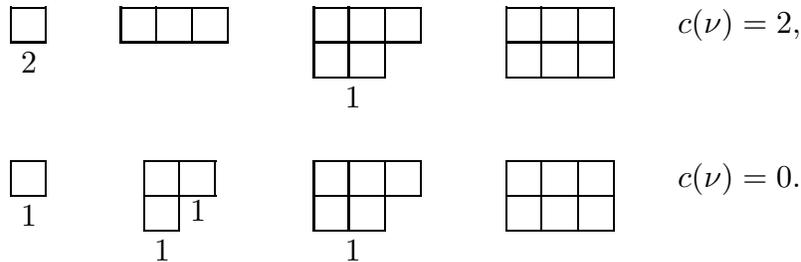
$$\matrix{\hbox{
          \normalbaselines\m@th\offinterlineskip
          \vtop{\hbox{{\Fsquare(0.4cm,)}}
          \vskip 2.8pt
          \hbox{~~2}}}
&&{\hbox{\Hthreebox(, ,)}}
&&\hbox{
          \normalbaselines\m@th\offinterlineskip
          \vtop{\hbox{{\Hthreebox( , ,)}}
          \vskip-0.4pt
          \hbox{{\Htwobox( , )}}
          \vskip 2.8pt
          \hbox{~~~~1}}}
&&\hbox{
          \normalbaselines\m@th\offinterlineskip
          \vtop{\hbox{{\Hthreebox( , , )}}
          \vskip-0.4pt
          \hbox{\Hthreebox( , , )}}}
&& c(\nu )=2,\cr \cr \cr
\hbox{
          \normalbaselines\m@th\offinterlineskip
          \vtop{\hbox{{\Fsquare(0.4cm,)}}
          \vskip 2.8pt
          \hbox{~~1}}}          
&&\hbox{
          \normalbaselines\m@th\offinterlineskip
          \vtop{\hbox{{\Htwobox( , )}}
          \vskip-0.4pt
          \hbox{{\Fsquare(0.4cm, )}~1}
          \vskip 2.8pt
          \hbox{~~1}}}
&&\hbox{
          \normalbaselines\m@th\offinterlineskip
          \vtop{\hbox{\Hthreebox( , ,)}
          \vskip-0.4pt
          \hbox{{\Htwobox( , )}}
          \vskip 2.8pt
          \hbox{~~~~1}}}
&&\hbox{
          \normalbaselines\m@th\offinterlineskip
          \vtop{\hbox{{\Hthreebox( , , )}}
          \vskip-0.4pt
          \hbox{\Hthreebox( , , )}}}          
&& c(\nu )=0.\cr \cr }
$$
Hence, the RHS(3.1)$=1+4t+7t^2+7t^3+4t^4+t^5+t^2(1+2t+3t^2+2t^3+t^4)=
1+4t+8t^2+9t^3+7t^4+3t^5+t^6$.

On the other hand, for the partition $\mu^+=(2211)$ the contribution to the 
RHS(3.1) is given by the following flags of partitions:
\vskip 0.2cm
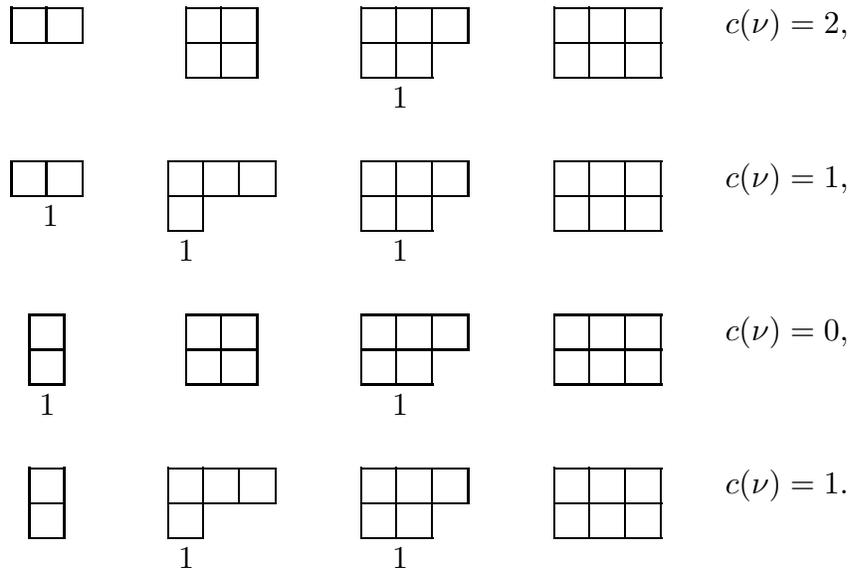
$$\matrix{{\hbox{
          \normalbaselines\m@th\offinterlineskip
          \Htwobox(, )}}
&&\hbox{
          \normalbaselines\m@th\offinterlineskip
          \vtop{\hbox{\Htwobox( ,)}
          \vskip-0.4pt
          \hbox{{\Htwobox( , )}}}}
&&\hbox{
          \normalbaselines\m@th\offinterlineskip
          \vtop{\hbox{\Hthreebox( , ,)}
          \vskip-0.4pt
          \hbox{{\Htwobox( , )}}
          \vskip 2.8pt
          \hbox{~~~~1}}}
&&\hbox{
          \normalbaselines\m@th\offinterlineskip
          \vtop{\hbox{{\Hthreebox( , , )}}
          \vskip-0.4pt
          \hbox{\Hthreebox( , , )}}}
&& c(\nu )=2,\cr \cr \cr
\hbox{
          \normalbaselines\m@th\offinterlineskip
          \vtop{\hbox{\Htwobox( , )}
          \vskip 2.8pt
          \hbox{~~~~1}}}          
&&\hbox{
          \normalbaselines\m@th\offinterlineskip
          \vtop{\hbox{{\Hthreebox( , , )}}
          \vskip-0.4pt
          \hbox{{\Fsquare(0.4cm, )}}
          \vskip 2.8pt
          \hbox{~~1}}}
&&\hbox{
          \normalbaselines\m@th\offinterlineskip
          \vtop{\hbox{\Hthreebox( , ,)}
          \vskip-0.4pt
          \hbox{{\Htwobox( , )}}
          \vskip 2.8pt
          \hbox{~~~~1}}}
&&\hbox{
          \normalbaselines\m@th\offinterlineskip
          \vtop{\hbox{{\Hthreebox( , , )}}
          \vskip-0.4pt
          \hbox{\Hthreebox( , , )}}}          
&& c(\nu )=1,\cr \cr \cr
\hbox{
          \normalbaselines\m@th\offinterlineskip
          \vtop{\hbox{\Vtwobox( , )}
          \vskip 2.8pt
          \hbox{~~1}}}  
&&\hbox{
          \normalbaselines\m@th\offinterlineskip
          \vtop{\hbox{\Htwobox( ,)}
          \vskip-0.4pt
          \hbox{{\Htwobox( , )}}}}
&&\hbox{
          \normalbaselines\m@th\offinterlineskip
          \vtop{\hbox{\Hthreebox( , ,)}
          \vskip-0.4pt
          \hbox{{\Htwobox( , )}}
          \vskip 2.8pt
          \hbox{~~~~1}}}
&&\hbox{
          \normalbaselines\m@th\offinterlineskip
          \vtop{\hbox{{\Hthreebox( , , )}}
          \vskip-0.4pt
          \hbox{\Hthreebox( , , )}}}
&& c(\nu )=0,\cr \cr \cr
\hbox{
          \normalbaselines\m@th\offinterlineskip
          \Vtwobox( , )}          
&&\hbox{
          \normalbaselines\m@th\offinterlineskip
          \vtop{\hbox{{\Hthreebox( , , )}}
          \vskip-0.4pt
          \hbox{{\Fsquare(0.4cm, )}}
          \vskip 2.8pt
          \hbox{~~1}}}
&&\hbox{
          \normalbaselines\m@th\offinterlineskip
          \vtop{\hbox{\Hthreebox( , ,)}
          \vskip-0.4pt
          \hbox{{\Htwobox( , )}}
          \vskip 2.8pt
          \hbox{~~~~1}}}
&&\hbox{
          \normalbaselines\m@th\offinterlineskip
          \vtop{\hbox{{\Hthreebox( , , )}}
          \vskip-0.4pt
          \hbox{\Hthreebox( , , )}}}          
&& c(\nu )=1.\cr \cr}
$$
Hence, the RHS(3.1)$=t^2(1+t+t^2)+t(1+3t+5t^2+5t^3+3t^4+t^5)+(1+2t+2t^2+t^3)
+t(1+2t+2t^2+t^3)=1+4t+8t^2+9t^3+7t^4+3t^5+t^6$.

We see that ${\cal P}_{\ld\mu}(t)={\cal P}_{\ld\mu^+}(t)$, but the 
corresponding sums of the products of $t$--binomial coefficients have 
different structures.

%\vfil\eject
\vskip 0.3cm
{\bf 3.2.} New combinatorial formula for the transition matrix $M(e,P)$.
\vskip 0.2cm

Now we are going to describe the fermionic formula for the following sum
$${\cal R}_{\ld\mu}(t)=\sum_{\eta}K_{\eta\mu}K_{\eta'\ld}(t).
$$
This sum is the $(\ld ,\mu)$--entry of the matrix transposed to the 
transition matrix between elementary and Hall--Littlewood polynomials, 
namely, if
$$e_{\ld}=\sum_{\mu}M(e,P)_{\ld\mu}P_{\mu}, \ \ {\rm then}
$$
$$M(e,P)_{\ld\mu}=\sum_{\nu}K_{\nu\ld}K_{\nu'\mu}(q)={\cal R}_{\mu\ld}(q).
$$

It is well--known ([Kn]) that ${\cal R}_{\ld\mu}(1)$ counts the number of 
(0,1)--matrices with row sums $\ld_i$ and column sums $\mu_j$. This 
number is equal also to the number of pairs of semistandard tableaux of 
conjugate shapes and weights $\ld$ and $\mu$, see, e.g. [M], Chapter~I, 
Section~6.

\vskip 0.2cm
{\bf Theorem 3.4.} ([HKKOTY]) {\it Let $\mu$ be a composition, $l(\mu )=r$. Then
$${\cal R}_{\ld\mu}(t)=\sum_{\{\nu\}}\prod_{k=1}^{r-1}\prod_{i\ge 1}
\left[\matrix{\nu_i^{(k+1)}-\nu_{i+1}^{(k+1)}\cr
\nu_i^{(k)}-\nu_{i+1}^{(k+1)}}\right]_t, \eqno (3.7)
$$
where the sum is taken over all flags of partitions $\nu=\{ 0=\nu^{(0)}
\subset\nu^{(1)}\subset\cdots\subset\nu^{(r)}=\ld'\}$ such that 
$\nu^{(k)}/\nu^{(k-1)}$ is a horizontal strip of length $\mu_k$, $1\le 
k\le r$.}

Proof of Theorem 3.4 will be given in Subsection~4.1.

\vskip 0.2cm
{\bf Remark.} The last condition on the flag $\nu$ means that $\nu$ 
defines a semistandard tableau of shape $\ld'$ and weight $\mu$. Thus, the 
number of terms in the RHS(3.7) is equal to that of semistandard tableaux 
of shape $\ld'$ and weight $\mu$.

\vskip 0.2cm
{\bf Examples.} $1^0$. It is clear that if $\mu =(1^n)$, then
$${\cal R}_{\ld\mu}(q)={\cal P}_{\ld\mu}(q)=X^{\ld}_{(1^n)}(q).
$$

$2^0$. If $\ld =(1^N)$, and $\mu =(\mu_1,\ldots ,\mu_n)$, $|\mu |=N$ then
$${\cal R}_{\ld\mu}(t)=\left[\matrix{N\cr \mu_1,\ldots ,\mu_n}\right]_t.
$$
Indeed, the RHS(3.7) contains only one product
$$\left[\matrix{\mu_1+\mu_2\cr \mu_1}\right]_t\left[\matrix{
\mu_1+\mu_2+\mu_3\cr \mu_1+\mu_2}\right]_t\cdots\left[\matrix{
N\cr \mu_1+\cdots +\mu_{n-1}}\right]_t.
$$

$3^0$. Let $\mu$ be a composition of length $n$, and 
$\ld =(2^{\ld_2}1^{\ld_1-\ld_2})$, so that $\ld'=(\ld_1,\ld_2)$. In this 
case the following partitions give the contribution to the RHS(3.7):

$\bullet$ $\nu^{(k)}=(\nu_1^{(k)},\nu_2^{(k)})$, 
$|\nu^{(k)}|=\mu_1+\cdots +\mu_k$, $1\le k\le n$;

$\bullet$ $0\le\nu_2^{(k)}\le\nu_2^{(k+1)}\le\nu_1^{(k)}\le\nu_1^{(k+1)}$,
$1\le k\le n-1$, $\nu_2^{(1)}=0$, $\nu^{(n)}=(\ld_1,\ld_2)$.

If we define $m_k=\nu_1^{(k+1)}-\nu_1^{(k)}$, $1\le k\le n-1$, 
$m_0=\mu_1$, then the RHS(3.7) takes the following form
$$\sum_{m\in{\bf Z}_{\ge 0}^{n-1}}\left[\matrix{\ld_2\cr 
\mu_1-m_1,\mu_2-m_2,\ldots ,\mu_n-m_{n-1}}\right]_t\prod_{k=1}^{n-1}
\left[\matrix{\ds\sum_{i=0}^k(2m_i-\mu_{i+1})\cr m_k}\right]_t,
$$
summed over all sequences $m\in{\bf Z}_{\ge 0}^{n-1}$, such that 
$m_1+\cdots +m_{n-1}=\ld_1-\mu_1$.

Finally, let us assume that $n=3$ and $\ld_1=\ld_2$. Then $m_2=0$, 
$m_1=\ld_1-\mu_1$ and the RHS(3.7) takes the form $(\mu_1+\mu_2+\mu_3=N)$
$$\left[\matrix{N\cr N-\mu_1,N-\mu_2, N-\mu_3}\right]_t.
$$
Hence, the number of (0,1)--matrices of size $N\times 3$ with row sums 
$\mu_i$, $i=1,2,3$, and column sums $\ld_i=2$, $1\le i\le N$, is equal to 
$\ds{N!\over (N-\mu_1)!(N-\mu_2)!(N-\mu_3)!}$.

$4^0$. Consider $\mu =(1221)$ and $\ld =(321)$. The summands in the 
RHS(3.7) correspond to the following flags of partitions $\nu 
=\{\nu^{1}\subset\nu^{(2)}\subset\nu^{(3)}\subset\nu^{(4)}\}$:
$$\matrix{\hbox{
          \normalbaselines\m@th\offinterlineskip
          \Fsquare(0.4cm,)}
&&\hbox{
          \normalbaselines\m@th\offinterlineskip
          \Hthreebox(, ,)}
&&\hbox{
          \normalbaselines\m@th\offinterlineskip
          \vtop{\hbox{{\Hthreebox( , ,)}}
          \vskip-0.4pt
          \hbox{\Htwobox( , )}}}
&&\hbox{
          \normalbaselines\m@th\offinterlineskip
          \vtop{\hbox{{\Hthreebox( , , )}}
          \vskip-0.4pt
          \hbox{\Htwobox( , )}
          \vskip-0.4pt
          \hbox{\Fsquare(0.4cm,)}}}
&&\matrix{\cr \cr 1&2&2\cr 3&3\cr 4}
&&\matrix{ \cr\cr \left[\matrix{ 3-0\cr 1-0}\right]_t,} \cr \cr
\hbox{
          \normalbaselines\m@th\offinterlineskip
          \Fsquare(0.4cm,)}          
&&\hbox{
          \normalbaselines\m@th\offinterlineskip
          \vtop{\hbox{\Htwobox( , )}
          \vskip-0.4pt
          \hbox{\Fsquare(0.4cm, )}}}
&&\hbox{
          \normalbaselines\m@th\offinterlineskip
          \vtop{\hbox{\Htwobox( , )}
          \vskip-0.4pt
          \hbox{\Htwobox( , )}
          \vskip-0.4pt
          \hbox{\Fsquare(0.4cm,)}}}
&&\hbox{
          \normalbaselines\m@th\offinterlineskip
          \vtop{\hbox{\Hthreebox( , , )}
          \vskip-0.4pt
          \hbox{\Htwobox( , )}
          \vskip-0.4pt
          \hbox{\Fsquare(0.4cm,)}}}
&&\matrix{\cr \cr 1&2&4\cr 2&3\cr 3}
&&\matrix{\cr 1,} \cr \cr 
\hbox{
          \normalbaselines\m@th\offinterlineskip
          \Fsquare(0.4cm,)}          
&&\hbox{
          \normalbaselines\m@th\offinterlineskip
          \vtop{\hbox{\Htwobox( , )}
          \vskip-0.4pt
          \hbox{\Fsquare(0.4cm, )}}}
&&\hbox{
          \normalbaselines\m@th\offinterlineskip
          \vtop{\hbox{\Hthreebox( , ,)}
          \vskip-0.4pt
          \hbox{\Fsquare(0.4cm,)}
          \vskip-0.4pt
          \hbox{\Fsquare(0.4cm,)}}}
&&\hbox{
          \normalbaselines\m@th\offinterlineskip
          \vtop{\hbox{\Hthreebox( , , )}
          \vskip-0.4pt
          \hbox{\Htwobox( , )}
          \vskip-0.4pt
          \hbox{\Fsquare(0.4cm,)}}}
&&\matrix{\cr\cr 1&2&3\cr 2&4\cr 3}
&&\matrix{\cr\cr \left[\matrix{3-1\cr 2-1}\right]_t,} \cr \cr           
\hbox{
          \normalbaselines\m@th\offinterlineskip
          \Fsquare(0.4cm,)}          
&&\hbox{
          \normalbaselines\m@th\offinterlineskip
          \vtop{\hbox{\Htwobox( , )}
          \vskip-0.4pt
          \hbox{\Fsquare(0.4cm, )}}}
&&\hbox{
          \normalbaselines\m@th\offinterlineskip
          \vtop{\hbox{\Hthreebox( , ,)}
          \vskip-0.4pt
          \hbox{\Htwobox( , )}}}
&&\hbox{
          \normalbaselines\m@th\offinterlineskip
          \vtop{\hbox{\Hthreebox( , , )}
          \vskip-0.4pt
          \hbox{\Htwobox( , )}
          \vskip-0.4pt
          \hbox{\Fsquare(0.4cm,)}}}
&&\matrix{\cr\cr 1&2&3\cr 2&3\cr 4}
&&\matrix{\cr\cr \left[\matrix{2-0\cr 1-0}\right]_t.}
\cr \cr \cr}
$$
Hence, the RHS(3.7)$=(1+t+t^2)+1+(1+t)+(1+t)=4+3t+t^2$.

Let us remark that RHS(3.7) does not depend on the permutations of 
components of the composition $\mu$. This is clear since the LHS(3.7) 
does. However, the number of summands in the RHS(3.7) do depends on the 
composition $\mu$, but not only on the corresponding partition $\mu^+$.
\vskip 0.5cm

{\bf \S 4. Proofs of Theorems 3.1 and 3.4.}
\vskip 0.2cm

Let $f_{\mu\nu}^{\ld}(t)$ be the structural constants for the 
Hall--Littlewood functions, i.e.
$$P_{\mu}(x;t)P_{\nu}(x;t)=\sum_{\ld}f_{\mu\nu}^{\ld}(t)P_{\ld}(x;t).
\eqno (4.1)
$$
It is well--known (see, e.g., [M], Chapter~III, \S 3, p.215, formula 
(3.2)) that
$$f_{\mu (1^m)}^{\ld}(t)=\prod_{i\ge 1}\left[\matrix{\ld'_i-\ld'_{i+1}\cr
\ld'_i-\mu'_i}\right]_t, \eqno (4.2)
$$
and therefore $f_{\mu (1^m)}^{\ld}(t)=0$ unless $\ld$--$\mu$ is a vertical 
$m$--strip.

Now let $T$ be a pure supertableau of shape $\ld$ and weight $\mu$, i.e. 
$T$ is a sequence of partitions 
$0=\ld^{(0)}\subset\ld^{(1)}\subset\cdots\subset\ld^{(r)}=\ld $, such 
that each skew diagram $\ld^{(i)}-\ld^{(i-1)}$ ($1\le i\le r$) is a 
vertical $\mu_i$--strip. For such tableau $T$, let us define
$$f_T(t)=\prod_{i\ge 1}f_{\ld^{(i-1)}(1^{\mu_i})}^{\ld^{(i)}}(t).
$$
Then the RHS(3.7) can be rewritten in the following form 
$\ds\sum_Tf_T(t)$, summed over all pure supertableaux of shape $\ld$ and 
weight $\mu$.
\vskip 0.3cm

{\bf 4.1.} Proof of Theorem 3.4.
\vskip 0.2cm

%We will use notations from the previous remark. 
It is well known that 
the Hall--Littlewood polynomial $P_{\ld}(X_n;t)$, when $\ld =(1^m)$, 
coincides with the $m$-th elementary symmetric function in the variables $X_n$:
$$P_{(1^m)}(X_n;t)=e_m(X_n),
$$
see e.g., [M], Chapter III, (2.8).

Using (4.1) and (4.2) we can write
$$e_m(x)P_{\nu}(x;t)=\sum_{\ld}f_{\nu (1^m)}^{\ld}(t)P_{\ld}(x;t), 
\eqno (4.3)
$$
and more generally using induction,
$$e_{\mu_1}(x)\ldots e_{\mu_r}(x)P_{\nu}(x;t)=\sum_{\ld}R_{\ld\mu}^{(\nu )}
P_{\ld}(x;t), \eqno (4.4)
$$
where
$$R_{\ld\mu}^{(\nu )}(t)=\sum f_T(t), \eqno (4.5)
$$
summed over all pure supertableaux $T$ of skew shape $\ld -\nu$ and 
weight $\mu$; in other words, the sum in (4.5) is taken over all sequences 
of partitions $\nu =\ld^{(0)}\subset\ld^{(1)}\subset\cdots\subset\ld^{(r)}
=\ld$, such that each skew diagram $\ld^{(i)}-\ld^{(i-1)}$ ($1\le i\le 
r$) is a vertical $\mu_i$--strip, and 
$$f_T(t)=\ds\prod_{i=1}^rf_{\ld^{(i-1)}(1^{\mu_i})}^{\ld^{(i)}}(t).
$$

To finish the proof of Theorem~3.4 we need the following formulae (see, 
e.g., [M], Table~1 on p.101 and Table on p.241):
$$\eno{
e_{\mu_1}(x)\cdots e_{\mu_r}(x)&=\sum_{\eta}K_{\eta'\mu}s_{\eta}(x),\cr
s_{\ld}(x)&=\sum_{\mu}K_{\ld\mu}(t)P_{\mu}(x;t). & (4.6)}
$$
Thus, we have
$$e_{\mu_1}(x)\cdots e_{\mu_r}(x)P_{\nu}(x;t)=
\sum_{\ld}\left(\sum_{\eta ,\beta}K_{\eta'\mu}K_{\eta\beta}(t)
f_{\nu\beta}^{\ld}(t)\right) P_{\ld}(x;t),
$$
and consequently,
$$R_{\ld\mu}^{(\nu )}(t)=\sum_{\eta ,\beta}K_{\eta'\mu}K_{\eta\beta}(t)
f_{\nu\beta}^{\ld}(t).\eqno (4.7)
$$
Finally, if we take $\nu =\emptyset$ in (4.7), then 
$f_{\emptyset\beta}^{\ld}(t)=\delta_{\ld\beta}$, and formula (3.5) follows.

\qed

{\bf 4.2.} Proof of Theorem 3.1.
\vskip 0.2cm

Proof of Theorem 3.1 is similar to that of Theorem~3.4 and based on the 
following
\vskip 0.2cm

{\bf Lemma 4.1.} {\it Let $\mu$ be a partition, $l(\mu )\le n$, and
$$h_k(X_n)P_{\mu}(X_n;t)=\sum_{\ld}g^{\ld}_{\mu}(t)P_{\ld}(X_n;t), \eqno 
(4.7)
$$
where $h_k(X_n)$ denotes the complete homogeneous symmetric function of 
degree $k$ in the variables $X_n=(x_1,\ldots ,x_n)$. Then
$$g_{\mu}^{\ld}(t)=t^{\ds\sum_{i\ge 1}\pmatrix{\ld'_i-\mu'_i\cr 2}}
\prod_{i\ge 1}\left[\matrix{\ld_i'-\mu'_{i+1}\cr \ld_i'-\mu_i'}\right]_t, 
\eqno (4.8)
$$
and therefore $g_{\mu}^{\ld}(t)=0$ unless $\mu\subset\ld$, $|\ld /\mu |=k$.}
\vskip 0.2cm

Let us postpone the proof of Lemma~4.1 to the end of this subsection and show 
first how using the formula (4.8) one can deduce the formula (3.1) from 
Theorem~3.1.

To do this we will need the formula (4.6) and the following one (see, e.g., 
[M], Table~1 on p.101):
$$h_{\mu_1}(x)\ldots h_{\mu_r}(x)=\sum_{\eta}K_{\eta\mu}s_{\eta}(x).
$$
Thus, we have
$$h_{\mu_1}(x)\ldots h_{\mu_r}(x)P_{\nu}(x;t)=\sum_{\ld}\left(\sum_{\eta 
,\beta}K_{\eta\mu}K_{\eta\beta}(t)f_{\nu\beta}^{\ld}(t)\right) P_{\ld}(x;t). 
\eqno (4.9)
$$
On the other hand, we can compute the LHS(4.9) using Lemma~4.1. Namely,
$${\rm LHS}(4.9)=\sum_{\ld}{\cal P}_{\ld\mu}^{(\nu )}(t)P_{\ld}(x;t), 
\eqno (4.10)
$$
where
$${\cal P}_{\ld\mu}^{(\nu )}(t)=\sum_{\pi}g_{\pi}(t), \eqno (4.11)
$$
summed over all reverse plain partitions $\pi$ of skew shape $\ld -\nu$ 
and weight $\mu$; in other words, the sum in (4.11) is taken over all 
sequences of partitions $\nu 
=\ld^{(0)}\subset\ld^{(1)}\subset\cdots\subset\ld^{(r)}=\ld$ such that 
$|\ld^{(i)}/\ld^{(i-1)}|=\mu_i$, $1\le i\le r$, and
$$g_{\pi}(t)=\prod_{i=1}^rg_{\ld^{(i-1)}}^{\ld^{(i)}}(t).
$$
Thus, it follows from (4.9)--(4.11) that 
$${\cal P}_{\ld\mu}^{(\nu )}(t)=\sum_{\eta 
,\beta}K_{\eta\mu}K_{\eta\beta}(t)f_{\nu\beta}^{\ld}(t). \eqno (4.12)
$$
Finally, if we take $\nu =\emptyset$ in (4.12), then 
$f_{\emptyset\beta}^{\ld}(t)=\delta_{\ld\beta}$, and formula (3.1) follows.
\vskip 0.2cm

{\it Proof of Lemma 4.1.} We will prove (4.8) by induction on the number 
$|\ld /\mu |$. It is clear that if $|\ld /\mu |=1$, then 
$g_{\mu}^{\ld}=f_{\mu ,(1)}^{\ld}={\rm RHS}(4.8)$. Because of the 
relation $\ds\sum_{r=0}^k(-1)^re_rh_{k-r}=\delta_{k,0}$, it is enough to 
prove that if $\nu\subset\ld$, $|\ld\setminus\nu |>0$, then
$$\sum_{\mu}(-1)^{|\mu -\nu |}f_{\nu ,(1^{|\mu -\nu |})}^{\mu}
g_{\mu}^{\ld}=0, \eqno (4.13)
$$
summed over all partitions $\mu$ such that $\nu\subset\mu\subset\ld$. 
Now, using (4.8) and (4.2), we can write
$$\eno{
{\rm RHS}(4.13)&=\sum_{\nu\subset\mu\subset\ld}(-1)^{\sum(\mu_i'-\nu_i')}
t^{\sum\pmatrix{\ld_i'-\mu_i'\cr 2}}\prod_{i\ge 1}\left[\matrix{\ld_i'-
\mu_{i+1}'\cr \ld_i'-\mu_i'}\right]_t\left[\matrix{\mu_i'-\mu_{i+1}'\cr
\mu_i'-\nu_i'}\right]_t\cr
&=\prod_i\Phi_i(t), & (4.14)}
$$
where
$$\Phi_i(t)=\sum_{\nu_i'\le\mu_i'\le\ld_i'}(-1)^{\mu_i'-\nu_i'}
t^{\pmatrix{\ld_i'-\mu_i'\cr 2}}{(t;t)_{\ld_{i-1}'-\mu_i'}\over
(t;t)_{\ld_i'-\mu_i'}(t;t)_{\mu_i'-\nu_i'}(t;t)_{\nu_{i-1}'-\mu_i'}},
$$
$\ld_0'=\nu_0'=0$, and by definition $(t;t)_m=0$, if $m<0$.

Consider at first $\Phi_1(t)$. We have
$$\eno{
(-1)^{\ld_1'-\nu_1'}(t)_{\ld_1'-\nu_1'}\Phi_1(t)&=
\sum_{\nu_1'\le\mu_1'\le\ld_1'}(-1)^{\ld_1'-\mu_1'}t^{\pmatrix{\ld_1'-\mu_1'\cr 
2}}\left[\matrix{\ld_1'-\nu_1'\cr \ld_1'-\mu_1'}\right]_t\cr
&=\sum_{m\ge 0}(-1)^mt^{\pmatrix{m\cr 2}}\left[\matrix{\ld_1'-\nu_1'\cr 
m}\right]_t=\delta_{\ld_1',\nu_1'}.}
$$
The last equality follows from the $q$--binomial theorem
$$\sum_{m=0}^N(-z)^mq^{\pmatrix{m\cr 2}}\left[\matrix{N\cr m}\right]_q=
(z,q)_N:=\prod_{i=1}^N(1-q^{i-1}z).
$$
Thus, if the product (4.14) does not equal to zero, then $\ld_1'=\nu_1'$, 
and
$$(-1)^{\ld_2'-\nu_2'}(t)_{\ld_2'-\nu_2'}\Phi_2(t)=\sum_{m\ge 0}(-1)^m
t^{\pmatrix{m\cr 2}}\left[\matrix{\ld_2'-\nu_2'\cr m}\right]_t=
\delta_{\ld_2',\nu_2'}.
$$
Repeating these arguments we see that the product (4.14) does not equal 
to zero only if $\ld =\nu$. But this is a contradiction with our 
assumption $|\ld /\nu |>0$. This proves (4.13) and 
(by induction) Lemma~4.1.

\qed

{\bf Remark.} The similar proofs of Theorems~3.1 and 3.4 can be found in 
[HKKOTY].

It seems the formula (4.8) is new. The formula (4.12) in the case $\nu 
=\emptyset$, probably, goes back to R.~Stanley, unpublished; see, e.g., 
[Bu2], Lemma~3.1.
\vskip 0.2cm

{\bf Corollary 4.2.} {\it Let $\ld$ and $\mu$ be partitions, $|\mu |=n$, 
and $f^{\ld}_{\nu\mu}(t)$ be the structural constants for the 
Hall--Littlewood functions, see [M], Chapter~III, or Section~4, (4.1). 
Then}
$$\sum_{\nu}t^{n(\nu )}f^{\ld}_{\nu\mu}(t)=t^{\sum_{i\ge 
1}\pmatrix{\ld_i'-\mu_i'\cr 2}}\prod_{i\ge 1}\left[\matrix{\ld_i'-
\mu_{i+1}'\cr\ld_i'-\mu_i'}\right]_t. \eqno (4.15)
$$
\vskip 0.2cm

{\it Proof.} It follows from Lemma~4.1 that the 
RHS(4.15)$=g_{\mu}^{\ld}(t)$. Hence,
$$\sum_{\ld}g_{\mu}^{\ld}(t)P_{\ld}=h_nP_{\mu}=\sum_{\nu}K_{(n)\nu}(t)
P_{\nu}P_{\mu}=\sum_{\ld}\left(\sum_{\nu}K_{(n)\nu}(t)
f_{\nu\mu}^{\ld}(t)\right)P_{\ld},
$$
and consequently, $g_{\mu}^{\ld}(t)=\ds\sum_{\nu}K_{(n)\nu}(t)
f_{\nu\mu}^{\ld}(t)$. The identity (4.15) follows from a simple 
observation that $K_{(n)\nu}(t)=t^{n(\nu )}$.

\qed

If $\mu =(1^n)$, then the RHS(4.15)$=t^{n(\ld )-n(\mu )}\left[\matrix{
\ld_1'\cr n}\right]_{t^{-1}}$, and identity (4.15) is reduced to that 
in [M], Chapter~III, Example~1.

\vskip 0.2cm
{\bf Exercise.} Let $\mu =(\mu_1,\ldots ,\mu_s)$ and $\nu =(\nu_1,\ldots 
,\nu_r)$ be compositions. For each partition $\eta$, $l(\eta )\le n$, 
denote by $K_{\eta ,\mu |\nu}$ the multiplicity of the highest weight 
%$\eta$ 
irreducible representation $V_{\eta}^{(n)}$ of the general linear 
group $\g l(n)$ in the tensor product
$$V_{\mu_1}^{(n)}\otimes\cdots\otimes V_{\mu_s}^{(n)}\otimes 
V_{(1^{\nu_1})}^{(n)}\otimes\cdots\otimes V_{(1^{\nu_r})}^{(n)}.
$$
Let $\ld$ be a partition, $l(\ld )\le n$. Find a fermionic formulae for 
the following sum
$$\sum_{\eta}K_{\eta, \mu |\nu}K_{\eta\ld}(q),
$$
which generalizes (3.1) and (3.7).

\vskip 0.2cm
{\bf Conjecture 4.3.} {\it Let $\ld ,\mu ,\nu$ be partitions. Define a 
family of polynomials $g_{\mu ;\nu}^{\ld}(t)$ via decomposition
$$s_{\nu}(x)P_{\mu}(x;t)=\sum_{\ld}g_{\mu ;\nu}^{\ld}(t)
P_{\ld}(x;t).
$$
Then $g_{\mu ;\nu}^{\ld}(t)$ is a polynomial with nonnegative integer 
coefficients.}
\vskip 0.2cm

{\bf Problem 2.} Find a combinatorial formula for polynomials $g_{\mu 
;\nu}^{\ld}(t)$.
\vskip 0.2cm

The answer on this problem is known when either $\nu =(1^N)$, see, e.g., 
[M], p.215, or $\nu =(n)$, see Lemma~4.1.

\vskip 0.5cm
{\bf \S 5. Polynomials ${\cal P}_{\ld\mu}(t)$ and their interpretations.}
\vskip 0.3cm

In this Section we summarize the known interpretations and some
properties of polynomials ${\cal P}_{\ld\mu}(t)$. The main reason
for this is the following: we suppose that all generalizations
of polynomials ${\cal P}_{\ld\mu}(t)$ considered in the coming
sections, should have properties similar to (5.2)-(5.10).

Polynomials ${\cal P}_{\ld\mu}(t)$ admit the following interpretations:

$\bullet$ Transition coefficients between modified Hall-Littlewood 
polynomials and monomial symmetric functions
$$Q'_{\ld}(X_n;t)=\sum_{\mu}{\cal P}_{\ld\mu}(t)m_{\mu}(X_n). \eqno (5.2)
$$

$\bullet$ Inhomogeneous unrestricted one dimensional sum with "special 
boundary conditions":
$${\cal P}_{\ld\mu}(t)=t^{n(\mu')}\sum_{m\in{\cal P}_{\ld\mu}}t^{E(m)},
\eqno (5.3)
$$
summed over the set ${\cal P}_{\ld\mu}$ of all transport matrices $m$ of 
type $(\ld ;\mu )$, i.e. the set of all matrices of non--negative 
integers with row sums $\ld_i$ and column sums $\mu_j$; $E(m)$ stands for 
the value of energy function $E(p)$ of the path $p$ which corresponds to the 
transport matrix $m$ under a natural identification of the set of paths 
${\cal P}_{\mu}(b_{\max},\ld )$ (see, e.g., [KMOTU2], or Subsection~3.1, 
Example~$1^0$) with that of transport matrices ${\cal 
P}_{\ld\mu}$.

\vskip 0.2cm
{\bf Problem 3.} Find a combinatorial rule for computation of the 
energy function $E(m)$ of a transport matrix $m\in{\cal P}_{\ld\mu}$.
\vskip 0.2cm

$\bullet$ Generating function of a generalized mahonian statistic $\varphi$ 
on the set of transport matrices ${\cal P}_{\ld\mu}$:
$${\cal P}_{\ld\mu}(t)=t^{n(\mu')}\sum_{m\in{\cal P}_{\ld\mu}}t^{\varphi 
(m)}.
$$
For examples of generalized mahonian statistics, see Section~2.6. 

\vskip 0.2cm
{\bf Problem 4.} It is natural to ask: are there exist combinatorial 
analogues of statistics 
$INV$, $MAJ$, $\wt{MAJ}$, $Z$, $\wt Z$ and $DEN$ (see Subsection~2.1), 
and $LP$ (see Definition~3.3) on the set of transport matrices ${\cal 
P}_{\ld\mu}$ with generating function ${\cal P}_{\ld\mu}(t)$?
\vskip 0.2cm

$\bullet$ The Poincare polynomial of the partial flag variety ${\cal 
F}_{\mu}^{\ld}/{\bf C}$:
$${\cal P}_{\ld\mu}(t)=\sum_{i\ge 0}t^{n(\ld )-i}\dim H_{2i}({\cal 
F}^{\ld}_{\mu};{\bf Z}). \eqno (5.4)
$$

$\bullet$ The number of ${\bf F}_q$--rational points of the partial flag 
variety ${\cal F}_{\mu}^{\ld}/{\bf F}_q$:
$$q^{n(\ld )}{\cal P}_{\ld\mu}(q^{-1})={\cal F}_{\mu}^{\ld}({\bf F}_q). 
\eqno (5.5)
$$

$\bullet$ The number of chains of subgroups
$$\{ e\}\subseteq H^{(1)}\subseteq H^{(2)}\subseteq\cdots\subseteq 
H^{(m)}\subseteq G
$$
in a finite abelian $p$--group $G$ of type $\ld$, such that each subgroup 
$H^{(i)}$ has order $p^{\mu_1+\cdots +\mu_i}$:
$$\al_{\ld}(S;p)=p^{n(\ld )}{\cal P}_{\ld\mu}(p^{-1}), \eqno(5.6)
$$
where $S:=S(\mu )=(\mu_1,\mu_1+\mu_2,\ldots ,\mu_1+\mu_2+\cdots +\mu_m)$, 
and $l(\mu )=m+1$. 

$\bullet$ String function of affine Demazure's module $V_w(l\Lambda_L)$ 
corresponding to the element $w=r_{Ln-1}r_{Ln-2}\ldots r_{L+2}r_{L+1}r_L$ 
of the affine Weyl group $W(A_{n-1}^{(1)})$:
$$t^{E_0}{\cal P}_{(l^L)\mu}(t)=\sum_{n\ge 0}\dim V_w(l\Lambda_L)_{\mu 
-n\delta}t^n, \eqno (5.7)
$$
for some known constant $E_0$; see [KMOTU2], or Subsection~1.6.

$\bullet$ Generalized $t$--supernomial and $t$--multinomial coefficients 
$\left[\matrix{\ld\cr\mu}\right]^{(0)}$ and $T^{(0)}(\ld ;\mu)$:
$$\eno{
&\left[\matrix{\ld\cr\mu}\right]^{(0)}=\sum_{\eta}K_{\eta\mu}\wt 
K_{\eta\ld}(t)=t^{n(\ld )}\sum_{\eta}K_{\eta\mu}K_{\eta\ld}(t^{-1}), & 
(5.8)\cr
&T^{(0)}(\ld ;\mu)=t^{-E_{\min}}{\cal P}_{\ld\mu}(t), & (5.9)}
$$
for some known constant $E_{\min}$.

As it was shown in Subsection~3.1, the coefficients (5.8) and (5.9) are a
natural generalization of those introduced by A.~Schilling and S.O.~Warnaar 
in the case $l(\mu )=2$, see [Ki2], [Sc], [ScW], [W].

$\bullet$ "Fermionic expression". Let $\ld$ be a partition and $\mu$ be a 
composition, $l(\mu )=n$, then
$${\cal 
P}_{\ld\mu}(t)=\sum_{\{\nu\}}t^{c(\{\nu\})}\prod_{k=1}^{n-1}\prod_{i\ge 
1}\left[\matrix{(\nu^{(k+1)})'_i-(\nu^{(k)})'_{i+1}\cr 
(\nu^{(k)})'_i-(\nu^{(k)})'_{i+1}}\right], \eqno (5.10)
$$
summed over all flags of partitions $\nu =\{ 
0=\nu^{(0)}\subset\nu^{(1)}\subset\cdots\subset\nu^{(n)}=\ld\}$, such 
that $|\nu^{(k)}|=\mu_1+\cdots +\mu_k$, $1\le k\le n$, and 
$$c(\{\nu\} )=\sum_{k=0}^{n-1}\sum_{i\ge 
1}\pmatrix{(\nu^{(k+1})'_i-(\nu^{(k)})'_i\cr 2}.
$$
See [HKKOTY] and Sections~3 and 4, where further details and applications of 
the fermionic formula (5.10) can be found.

$\bullet$ Truncated form, or finitization of characters and branching 
functions of (some) integrable representations of the affine Lie algebra 
of type $A_{n-1}^{(1)}$, and more generally, for Kac--Moody algebras, 
$W$--algebras, $\ldots$. 

The observation that certain special limits of polynomials
${\cal P}_{\ld\mu}(t)$ and Kostka--Foulkes polynomials may play an
important role in the representation theory of affine Lie algebras
originally was made in [Ki2].
It was observed in [Ki2] that the character formula for the 
level 1 vacuum representation $V(\Lambda_0)$ of the affine Lie algebra of 
type $A_{n-1}^{(1)}$ (see, e.g., [Kac], Chapter~13) can be obtained as an 
appropriate limit $N\to\infty$ of the modified Hall--Littlewood 
polynomials $Q'_{(1^N)}(X_n;q)$. The proof was based on the following 
formula
$${\cal P}_{(1^N)\mu}(q)=q^{n(\mu')}\left[\matrix{N\cr \mu_1,\ldots 
,\mu_n}\right]_q, \eqno (5.11)
$$
see [Ki2], (2.28), or Subsection~1.2, (1.5).

The latter observation about a connection between the character 
ch$(V(\Lambda_0))$ and modified Hall-Littlewood polynomials 
$Q'_{(1^N)}(X_n;q)$ immediately implies that the level 1 branching 
functions $b_{\ld}^{\Lambda_0}(q)$ can be obtained as an appropriate 
limit $\ld_N\to\infty$ of the "normalized" Kostka--Foulkes polynomials 
$q^{-A_N}K_{\ld_N,(1^N)}(q)$. We refer the reader to [Kac], Chapter~12, 
for definitions and basic properties of branching functions 
$b_{\ld}^{\Lambda}(q)$ corresponding to an integrable representation 
$V(\Lambda )$ of affine Lie algebra. 

It was conjectured in [Ki2], 
Conjecture~4, that the similar result should be valid for the branching 
functions $b_{\ld}^{\Lambda}(q)$ corresponding to the integrable highest 
weight $\Lambda$ irreducible representation $V(\Lambda )$ of the affine 
Lie algebra $\wh{sl}(n)$. This conjecture has been proved in [Ki2] in the 
following cases: $\wh{sl}(n)$ and $\Lambda =\Lambda_0$, $\wh{sl}(2)$ and 
$\Lambda =l\Lambda_0$, and $\wh{sl}(n)$ and $\Lambda =2\Lambda_0$. It had 
not been long before A.~Nakayashiki and Y.~Yamada [NY] proved this 
conjecture in the case $\wh{sl}(n)$ and $\Lambda =l\Lambda_i$, $0\le i\le 
n-1$. See also [KKN] for another proof of the result of Nakayashiki and 
Yamada in the case $i=0$. The general case has been investigated in [HKKOTY]. 
It happened that in general the so--called thermodynamical Bethe ansatz 
limit of Kostka--Foulkes polynomials gives the branching function of a 
certain reducible integrable representation of $\wh{sl}(n)$, see details 
in [HKKOTY]. 

\vskip 0.2cm
{\bf Problem 5.} Find an interpretation of 
the branching functions $b_{\ld}^{\Lambda}(q)$ of the integrable highest 
weight $\Lambda$ irreducible representation $V(\Lambda )$ of the affine 
Lie algebra $\wh{sl}(n)$ as the thermodynamical Bethe ansatz type limit of 
a certain family of the Kostka--Foulkes type polynomials.

\vskip 0.5cm

{\bf \S 6. Generalizations of polynomials ${\cal P}_{\ld\mu}(t)$
and $K_{\ld\mu}(t)$.}
\vskip 0.3cm

In this Section we summarize possible generalizations of polynomials 
${\cal P}_{\ld\mu}(t)$ and Kostka--Foulkes polynomials $K_{\ld\mu}(t)$, 
their properties, and some special cases. Let us remind that
$${\cal P}_{\ld\mu}(t)=\sum_{\eta}K_{\eta\mu}K_{\eta\ld}(t), \eqno (6.1)
$$
where $\ld$ is a partition, $\mu$ is a composition; summation in (6.1) 
runs over all partitions $\eta$; $K_{\eta\ld}(t)$ is the Kostka--Foulkes 
polynomial (see, e.g., 
[M], Chapter~III, Section~6), and $K_{\eta\mu}:=K_{\eta\mu}(1)$ is the 
Kostka number which is equal to the number of semistandard Young tableaux 
of shape $\eta$ and content $\mu$.

\vskip 0.3cm
{\bf 6.1.} Crystal Kostka polynomials.
\vskip 0.2cm

First let us recall the result of A.~Nakaya\-shi\-ki and 
Y.~Yamada [NY] that the Kostka--Foulkes 
polynomial $K_{\ld\mu}(t)$ coincides with the classically restricted one 
dimensional sum with special boundary conditions. For another proof, see, 
e.g., [KKN]; cf. [KMOTU2], [HKKOTY]. 

Let $n\ge 2$ be a natural integer which is fixed throughout this 
subsection.

\vskip 0.2cm
{\bf Definition 6.1.} {\it Let $R=\{ R_1,\ldots ,R_p\}$ be a sequence of 
partitions, $\mu$ be a partition such that $|\mu |=|R_1|+\cdots +|R_p|$. 
Define the polynomial $C{\cal P}_{R\mu }(t)$ to be the weight $\mu$ 
unrestricted one 
dimensional sum corresponding to the tensor product of crystals 
$B_{R_1}\otimes\cdots\otimes B_{R_p}$, and boundary condition 
$b_{T_{\min}}$, $T_{\min}\in STY(\ld ,\ld )$, where $B_{R_i}$ is the crystal 
(see, e.g., [Ka1]) corresponding to the irreducible  highest weight $R_i$
representation $V_{R_i}$ of the Lie algebra $sl(n)$.}
\vskip 0.2cm

{\bf Definition 6.2.} {\it The  crystal 
Kostka polynomial $CK_{\ld R}(q)$ corresponding to a set of 
partitions $R=\{ R_1,\ldots ,R_p\}$ is defined to be the weight $\ld$ 
classically restricted one dimensional sum corresponding to the tensor 
product of crystals $B_{R_1}\otimes\cdots\otimes B_{R_p}$, and boundary 
condition $b_{T_{\min}}$.}
\vskip 0.2cm

We refer the reader to [LS], [DLT] and [Ki1], where definition 
and basic properties of Kostka--Foulkes polynomials can be found, and
to [HKMOTU2] and [HKKOTY] for definitions of 
unrestricted, classically restricted and restricted one dimensional sums.  

Let us remark that
$$CK_{\ld R}(1)={\rm Mult}[V_{\ld}~:~V_{R_1}\otimes\cdots\otimes 
V_{R_p}], \eqno (6.2)
$$
i.e. $CK_{\ld R}(1)$ is equal to the multiplicity of the highest weight 
$\ld$ irreducible representation $V_{\ld}$ of $sl(n)$ in the tensor product 
$V_{R_1}\otimes\cdots\otimes V_{R_p}$. Thus, the crystal Kostka 
polynomial $CK_{\ld R}(q)$ may be considered as a $q$--analog of the tensor 
product multiplicity (6.2).

\vskip 0.3cm
{\bf 6.2.} Fusion Kostka polynomials.
\vskip 0.2cm

The problem of finding a ``natural'' $q$--analog of the tensor product
multiplicities has a long story. To our knowledge, there exists at least 
three natural algebraic ways to define a 
$q$--analog of the tensor product multiplicity (6.2). The first one is based 
on the so--called fusion rules for the tensor product of "restricted" 
representations of the quantized universal enveloping algebra $U_q(sl(n))$ 
when $q$ is a root of unity, see, e.g., [GoW]; [Kac], Exercises 
13.34-13.36; and [BKMW], where a combinatorial description of the fusion 
rules for representations of $sl(3)$ and $sl(4)$ are given. We denote by  
$FK_{\ld R}(q)$, and call it {\it fusion} Kostka polynomial, a $q$--analog 
of the tensor product multiplicity (6.2) which corresponds to the fusion rules.

Let us explain informally the meaning of the fusion Kostka polynomials 
$FK_{\ld R}(q)$. Let ${\cal F}_r(n)$ be the fusion algebra corresponding 
to the quantized universal enveloping algebra $U_q(sl(n))$, when $q=\exp 
\left(2\pi i/r+n\right)$. Each finite dimensional $sl(n)$--module $V$ 
defines an element $[V]$ of the fusion algebra ${\cal F}_r(n)$. This 
algebra is generated by the so--called "restricted" representations 
$V_{\ld}$, which correspond to partitions $\ld =(\ld_1,\ldots ,\ld_n)$ 
such that $\ld_1-\ld_n\le r$. It is well--known that the fusion algebra 
is commutative and associative.  
We denote by $\wh\otimes$ the product in the algebra ${\cal F}_r(n)$. 
This product depends on $r$ and $n$. Let $R=(R_1,\ldots ,R_p)$ be a 
sequence of partitions, denote by 
Mult$^{(r)}(V_{\ld}~:~V_{R_1}\wh\otimes\cdots\wh\otimes V_{R_p})$ the 
coefficient of $[V_{\ld}]$ in the decomposition of the product 
$[V_{R_1}]\wh\otimes\cdots\wh\otimes [V_{R_p}]$ in the fusion algebra 
${\cal F}_r(n)$:
$$[V_{R_1}]\wh\otimes\cdots\wh\otimes [V_{R_p}]=\sum_{\ld}
{\rm Mult}^{(r)}\left( V_{\ld}~:~V_{R_1}\wh\otimes\cdots\wh\otimes 
V_{R_p}\right) [V_{\ld}].
$$

\vskip 0.2cm
{\bf Definition 6.3.} {\it The fusion Kostka polynomial $FK_{\ld R}(q)$ 
is defined to be}
$$FK_{\ld R}(q)=\sum_{r\ge 0}\left({\rm Mult}^{(r+1)}(V_{\ld}~:~V_{R_1}
\wh\otimes\cdots\otimes V_{R_p})-{\rm Mult}^{(r)}(V_{\ld}~:~V_{R_1}
\wh\otimes\cdots\wh\otimes V_{R_p})\right) q^r. \eqno (6.3)
$$
\vskip 0.2cm

It is well--known that $FK_{\ld R}(q)$ is a polynomial with nonnegative 
integer coefficients, and
$$FK_{\ld R}(1)={\rm Mult}(V_{\ld}~:~V_{R_1}\otimes\cdots\otimes V_{R_p}).
$$
Thus, if all partitions $R_i$ have only one part $\mu_i$, i.e. 
$R_i=(\mu_i)$, then $FK_{\ld R}(1)=K_{\ld\mu}(1)$ is equal to the number 
$STY(\ld ,\mu)$ of semistandard Young tableaux of shape $\ld$ and weight 
$\mu$.

\vskip 0.2cm
{\bf Problem 6.} Give a combinatorial definition of a statistic on the 
set $STY(\ld ,\mu)$ which has the generating function $FK_{\ld R}(t)$.
\vskip 0.2cm
 
Let us give few illustrative examples of the fusion Kostka polynomials 
for Lie algebras $sl(5)$ and $sl(3)$.

$\bullet$ Algebra $sl(3)$:

\vskip 0.2cm
i) if $\ld =(433)$, $R=\{ (1)^{\otimes 10}\}$, then
$$FK_{\ld R}(q)=q+54q^2+115q^3+40q^4;
$$

ii) if $\ld =(422)$, $R=\{ (1)^{\otimes 8}\}$, then
$$FK_{\ld R}(q)=13q^2+30q^3+13q^4.
$$

$\bullet$ Algebra $sl(5)$:

\vskip 0.2cm
i) if $\ld =(65432)$,  $R=\{ (4321),(4321)\}$, then
$$FK_{\ld R}(q)=4q^4+10q^5+2q^6;
$$

ii) if $\ld =(98653)$, $R=\{ (6531),(6532)\}$, then
$$FK_{\ld R}(q)=4q^6+16q^7+13q^8+2q^9.
$$

{\bf Problem 7.} Let us introduce the {\it fusion} modified 
Hall--Littlewood polynomials
$$FQ_R'(X_n;t)=\sum_{\eta}FK_{\eta R}(t)s_{\eta}(X_n). \eqno (6.4)
$$
where $s_{\eta}(X_n)$ stands for the Schur function corresponding to a 
partition $\eta$.

Find algebraic, combinatorial, and geometric interpretations of 
the fusion modified Hall--Littlewood polynomials $FQ_n'(X_n;t)$.

\vskip 0.3cm
%\vfil\eject
{\bf 6.3.} Ribbon Kostka polynomials.
\vskip 0.2cm

The second way to define a $q$--analog of the tensor product multiplicity 
(6.2) is due to A.~Lascoux, B.~Leclerc and J.-Y.~Thibon, [LLT], and 
based on the using of ribbon tableaux. We refer the reader to [LLT], 
Sections~4 and 6, for definitions of a $p$--ribbon tableau $T$, spin 
$s(T)$ of a $p$--ribbon tableau, and "$p$--ribbon version" $\wt 
Q_{\ld}^{(p)}(X_n;t)$ of modified Hall--Littlewood polynomials. Here we 
are only reminding that if $\ld$ is a partition with empty $p$--core, 
then by definition
$$\wt Q_{\ld}^{(p)}(X_n;t)=\sum_{T\in{\rm Tab}_p(\ld ,\le n)}
t^{\bar s(T)}x^{w(T)}, \eqno (6.5)
$$
summed over the set Tab$_p(\ld ,\le n)$ of all $p$--ribbon tableaux of 
shape $\ld$ filled by numbers not exceeding $n$; $\bar s(T)=s(T)-\min\{ 
s(T)~|~T\in{\rm Tab}_p(\ld ,\le n)\}$ is a normalized spin of 
the $p$--ribbon tableau $T$, cf. [LLT], (25). It is known, [LLT], Theorem~6.1, 
that $\wt Q_{\ld}^{(p)}(X_n;t)$ is a symmetric polynomial. Let us define 
the ribbon polynomials ${\cal P}_{\ld\mu}^{(p)}(t)$ and the ribbon Kostka 
polynomials $K_{\ld\mu}^{(p)}(t)$ via decompositions (cf. [LLT]):
$$\eno{
&\wt Q_{\ld}^{(p)}(X_n;t)=\sum_{\mu}{\cal P}_{\ld\mu}^{(p)}(t)m_{\mu}(X_n),
& (6.6)\cr
&\wt Q_{\ld}^{(p)}(X_n;t)=\sum_{\mu}K_{\ld\mu}^{(p)}(t)s_{\mu}(X_n). & 
(6.7)}
$$

{\bf Remark.} The functions $\wt Q_{\ld}^{(p)}(X_n;t)$ were introduced 
and studied by A.~Lascoux, B.~Leclerc and J.-Y.~Thibon in [LLT], and 
denoted in [LLT] by $G_{\ld}^{(p)}(X_n;t)$. We denote these functions by 
$\wt Q_{\ld}^{(p)}(X_n;t)$, and call the {\it ribbon} modified 
Hall--Littlewood polynomials in order to underline a certain similarity with 
modified Hall--Littlewood polynomials $Q'_{\ld}(X_n;t)$. In fact, it was 
proved in [LLT], Theorem~6.6, that if $\ld$ is a partition, and $L\ge 
l(\ld )$, then
$$\wt Q_{L\ld}^{(L)}(X_n;t)=Q'_{\ld}(X_n;t),
$$
where $L\ld =(L\ld_1,L\ld_2,\ldots ,\ld L_n)$.
\vskip 0.2cm

It is well--known and goes back to D.~Littlewood, cf. [SW], that
$$K_{\ld\mu}^{(p)}(1)={\rm Mult}[V_{\mu}~:~V_{\ld^{(1)}}\otimes\cdots
\otimes V_{\ld^{(p)}}], \eqno (6.8)
$$
and ${\cal P}_{\ld\mu}^{(p)}(1)$ is equal to the number of weight $\mu$ 
unrestricted paths corresponding to the tensor product of crystals 
$B_{\ld^{(1)}}\otimes\cdots\otimes B_{\ld^{(p)}}$, where 
$(\ld^{(1)},\ldots ,\ld^{(p)})$ is the $p$--quotient of partition $\ld$ 
(see, e.g., [M], Chapter~I, Example~8, for definitions of $p$--core and 
$p$--quotient of a partition $\ld$).

Now we are going to formulate two conjectures about connections between 
polynomials $C{\cal P}_{R\mu}(t)$ and $CK_{\ld R}(t)$, see Definitions~6.1 
and 6.2, and the ribbon polynomials ${\cal P}_{\ld\mu}^{(p)}(t)$ and 
$K_{\ld\mu}^{(p)}(t)$. Namely, let $R=\{ R_1\ldots ,R_p\}$ be a sequence 
of partitions. According to the result of D.~Littlewood there exists the 
unique partition $\Lambda$ with the following properties (see, e.g., [M], 
Chapter~I, Example~8, and [SW]):
$$\eno{
i)&~~~p{\rm -core}(\Lambda )=\emptyset; & (6.9)\cr
ii)&~~~p{\rm -quotient}(\Lambda )=(R_1,\ldots ,R_p). ~~~~~~~~~~~~~~~~
~~~~~~~~~~~~~~~~~~~~~~~~~~~~~~~~~~~~~~~~~~~~~&(6.10)}
$$

{\bf Conjecture 6.4.} {\it Let $C{\cal P}_{R\mu}(t)$ be the weight $\mu$ 
unrestricted one dimensional sum corresponding to the tensor product of 
crystals $B_{R_1}\otimes\cdots\otimes B_{R_p}$, and boundary condition 
$b_{T_{\min}}$; let $\Lambda$ be the 
unique partition which satisfies the conditions (6.9) and (6.10). Then 
$$C{\cal P}_{R\mu}(t)=t^{E_0}{\cal P}_{\Lambda\mu}^{(p)}(t),
$$ 
for a certain constant $E_0$.}
\vskip 0.2cm

{\bf Conjecture 6.5.} {\it Let $CK_{\ld R}(t)$ be the weight $\ld$ 
classically restricted one dimensional sum corresponding to the tensor 
product of crystals $B_{R_1}\otimes\cdots\otimes B_{R_p}$ and boundary 
condition $b_{T_{\min}}$; let $\Lambda$ 
be the unique partition which satisfies the conditions (6.9) and (6.10). 
Then 
$$CK_{\ld R}(t)=q^{E_0}K_{\Lambda\ld}^{(p)}(t),
$$ 
for a certain constant $E_0$.}
\vskip 0.3cm

{\bf 6.4.} Generalized Kostka polynomials.
\vskip 0.2cm

The third way to define a $q$--analog, denoted by $K_{\ld R}(q)$, of 
the tensor product multiplicity (6.2), in the case $R$ is a sequence of 
rectangular 
partitions, is due to M.~Shimozono and J.~Weyman, see, e.g., [KS]. By 
definition the polynomials $K_{\ld R}(q)$ are the 
Poincare polynomials of isotypic components of Euler 
characteristics of certain ${\bf C}[\g l_n]$--modules supported in 
nilpotent conjugacy class closures.

To give precise definitions, we need little more notations. Our 
exposition follows to [KS]. Let $\eta =(\eta_1,\eta_2,\ldots ,\eta_p)$ be 
a sequence of positive integers that sum to $n$. Denote by ${\rm 
Root}s_{\eta}$ the set of ordered pairs $(i,j)$ such that $1\le 
i\le\eta_1+\cdots +\eta_r<j\le n$ for some $r$. For example, if $\eta 
=(1^n)$, then Root$s_{\eta}=\{ (i,j)~|~1\le i<j\le n\}$.

Let $X_n=(x_1,\ldots ,x_n)$ be the set of independent variables. For any 
sequence of integer numbers $\gamma =(\gamma_1,\ldots ,\gamma_n)$ we put 
$x^{\gamma}=x_1^{\gamma_1}\cdots x_n^{\gamma_n}$. The symmetric group 
$S_n$ acts on polynomials in $X_n=(x_1,\ldots ,x_n)$ by permuting 
variables. Define the operators $J$ and $\pi$ by
$$\eno{
J(f)&=\sum_{w\in S_n}(-1)^{l(w)}w(x^{\delta }f), & (6.11)\cr
\pi (f)&=J(1)^{-1}J(f), & (6.12)}
$$ 
where $J(1)=\ds\prod_{i<j}(x_i-x_j)$ is the Vandermond determinant, 
$\delta =(n-1,n-2,\ldots ,1,0)$.

For the dominant (weakly decreasing) integral weight $\ld 
=(\ld_1\ge\ld_2\ge\cdots\ge\ld_n)$, the character $s_{\ld}(X_n)$ of the 
highest weight $\ld$ $\g l(n)$ module $V_{\ld}$ is given by the Laurent 
polynomial $s_{\ld}(X_n)=\pi (x^{\ld})$. When $\ld$ is a partition (that 
is $\ld_n\ge 0$), $s_{\ld}$ is the Schur function.

Let $B_{\eta}(X_n;q)$, $H_{\gamma\eta}(X_n;q)$, and 
$K_{\ld ,\gamma ,\eta}(q)$ be the formal power series defined by
$$\eno{
B_{\eta}(X_n;q)&=\prod_{(i,j)\in{\rm Root}s_{\eta}}(1-qx_i/x_j)^{-1}, & 
(6.13)\cr
H_{\gamma\eta}(X_n;q)&=\pi\left( x^{\gamma}B_{\eta}(X_n;q)\right)=
\sum_{\ld}s_{\ld}(X_n)K_{\ld ,\gamma ,\eta}(q), & (6.14)}
$$
where $\ld$ runs over the dominant integral weights in ${\bf Z}^n$. It is 
known (M.~Shimozono and J.~Weyman) that the coefficients 
$K_{\ld ,\gamma ,\eta}(q)$ are in fact polynomials with integer coefficients. 
It is not true in general that the polynomials $K_{\ld ,\gamma ,\eta}(q)$ 
have nonnegative coefficients.

Now we are going to introduce the generalized Kostka polynomials $K_{\ld 
R}(q)$. Namely, let $R=(R_1,\ldots ,R_p)$ be a sequence of partitions. 
Denote by $\eta =(\eta_1,\ldots ,\eta_p)$ the sequence of lengths 
$\eta_i=l(R_i)$ of partitions $R_i$. Let $n=|\eta |$, and $\gamma 
(R)\in{\bf Z}^n_{\ge 0}$ denotes the composition obtained by 
concatenating the parts of the $R_i$ in order.

\vskip 0.2cm
{\bf Definition 6.6.} {\it The generalized Kostka polynomial $K_{\ld 
R}(q)$ corresponding to a partition $\ld$ and sequence of partitions $R$ 
is defined by the following formula}
$$K_{\ld R}(q)=K_{\ld ,\gamma (R),\eta}(q). \eqno (6.15)
$$
\vskip 0.2cm

It is known (M.~Shimozono and J.~Weyman) that
$$K_{\ld R}(1)={\rm Mult}(V_{\ld}~:~V_{R_1}\otimes\cdots\otimes V_{R_p}),
$$
i.e. $K_{\ld R}(1)$ is equal to the multiplicity of the highest weight 
$\ld$ irreducible representation $V_{\ld}$ of the Lie algebra $\g l(n)$ 
in the tensor product $V_{R_1}\otimes\cdots\otimes V_{R_p}$. The 
generalized Kostka polynomials are a far generalization of the 
Kostka--Foulkes polynomials $K_{\ld\mu}(q)$, the 
generalized exponents polynomials $F_q(V_{\ld})$ introduced by B.~Kostant 
and studied by R.~Gupta, W.~Hessenlink, S.~Kato, J.~Weyman, A.Broer, 
$\ldots$, see, e.g., [G], [DLT]. More precisely:

1. Let $R_i$ be the single row $(\mu_i)$ for all $i$, where $\mu 
=(\mu_1,\mu_2,\ldots )$ is a partition of length at most $n$. Then
$$K_{\ld R}(q)=K_{\ld\mu}(q), \eqno (6.16)
$$
where $K_{\ld\mu}(q)$ is the Kostka--Foulkes polynomial. The proof of 
(6.16) follows from the following well--known identity:
$$\pi\left(x^{\mu}\prod_{1\le i<j\le n}\left(1-q{x_i\over x_j}
\right)^{-1}\right) =\sum_{k\ge 0}e_n(X_n)^{-k}Q'_{\mu +(k^n)}(X_n;q).
\eqno (6.17)
$$
When $\mu =0$, the LHS(6.17)$=\ds\prod_{1\le i,j\le 
n}\left(1-qx_i/x_j\right)^{-1} =\sum_{k\ge 0}q^k{\rm ch}({\cal H}^k)$, where 
${\cal H}=\ds\bigoplus_{k\ge 0}{\cal H}^k$ is the graded module of harmonic 
polynomials. Follow [Gu], the generalized exponents polynomial $F_q(V)$ 
of a finite--dimensional $\g l(n)$--module $V$ is defined to be 
$F_q(V)=\ds\sum_{k\ge 0}\langle V,{\cal H}^k\rangle q^k$. It is follows 
immediately from (6.17) with $\mu =0$, that $F_q(V_{\ld})=0$, if $|\ld 
|\not\equiv 0(\mod n)$, and $F_q(V_{\ld })=K_{\ld (l^n)}(q)$, if $|\ld 
|=ln$. The last equality originally was proved by W.~Hessenlink, and 
"elementary" algebraic proof may be found in [DLT].

2. Let $R_i$ be the single column $(1^{\eta_i})$ for all $i$. Then
$$K_{\ld R}(q)=\wt K_{\ld'\eta^+}(q),
$$
is the cocharge Kostka--Foulkes polynomial, where $\ld'$ is the conjugate 
of the partition $\ld$ and $\eta^+$ is the partition obtaining by sorting 
the parts of $\eta$ into weakly decreasing order.

3. (M.~Shimozono and J.~Weyman). Let $k$ be a positive integer and $R_i$ be
the rectangle with $k$ columns and $\eta_i$ rows, $1\le i\le n$. Then 
$K_{\ld R}(q)$ is 
the Poincare polynomial of the isotypic component of the irreducible 
$GL(n)$--module of highest weight $(\ld_1-k,\ld_2-k,\ldots ,\ld_n-k)$ in 
the coordinate ring of the Zariski closure of the nilpotent conjugacy 
class which corresponds to the set of nilpotent matrices with the Jordan 
canonical form of type $(\eta^+)'$.

As it was mentioned, the generalized Kostka polynomials $K_{\ld R}(q)$ 
may have negative coefficients for general $\ld$ and $R$. Nevertheless, 
for the so--called dominant sequence of partitions $R$, one expect

\vskip 0.2cm
{\bf Conjecture 6.7} (A.~Broer, [KS]) {\it Let $R$ be a dominant sequence 
of partitions. Then}
$$K_{\ld R}(t)\in{\bf N}[t].
$$

Recall that a sequence of partitions $R=(R_1,\ldots ,R_p)$ is called {\it 
dominant}, if for all $1\le i\le p$, the last part of $R_i$ is at least 
as large as the first part of $R_{i+1}$. 

\vskip 0.3cm
{\bf 6.5.} Summary.
\vskip 0.2cm

In the previous Subsections we gave definitions of four families of 
polynomials which may be considered as the "natural" $q$--analogues of 
the tensor product multiplicities, namely,

$\bullet$ fusion Kostka polynomials $FK_{\ld R}(t)$,

$\bullet$ crystal Kostka polynomials $CK_{\ld R}(t)$,

$\bullet$ ribbon Kostka polynomials $K_{\Lambda\mu}^{(p)}(t)$,

$\bullet$ generalized Kostka polynomials $K_{\ld R}(t)$, 

\hskip -0.7cm where $R$ is a 
sequence of partitions, $\ld$ and $\mu$ are partitions, and $\Lambda$ is 
a partition without $p$--core.

It is natural to ask: what are the relations between these four families 
of polynomials?

First of all, for each sequence of partitions $R=(R_1,\ldots ,R_p)$ 
denote by $\Lambda :=\Lambda (R)$ the unique partition $\Lambda$ which 
has no $p$--core, and has $R$ as its $p$--quotient. It is known that

$\bullet$ $CK_{\ld R}(1)=FK_{\ld R}(1)=K_{\Lambda (R)\ld}^{(p)}(1)=
K_{\ld R}(1)=$RHS(6.2),

$\bullet$ $CK_{\ld R}(t)$, $FK_{\ld R}(t)$ are polynomials with 
nonnegative coefficients by definition.

It was conjectured in [LLT] that the ribbon Kostka polynomials 
$K_{\Lambda\mu}^{(p)}(t)$ have nonnegative coefficients. This conjecture 
was proved in [CL] in the case $p=2$. As for the generalized Kostka 
polynomials $K_{\ld R}(t)$, they do may have negative coefficients in 
general. For example, take $\ld =(2,2)$ and $R=((1),(3))$, then $K_{\ld 
R}(t)=t-1$.

It seems a very difficult problem to characterize all sequences of 
partitions $R=(R_1,\ldots ,R_p)$ such that $K_{\ld R}(t)\in{\bf N}[t]$ 
for all partitions $\ld$. But even if it happens that 
the generalized Kostka polynomial $K_{\ld R}(t)$ do has nonnegative 
coefficients for some $\ld$ and $R$, even in this case, 
$K_{\ld R}(t)\ne K_{\Lambda (R)\ld}$ in 
general. For example, take $\ld =(521)$ and $R=((31),(1),(1),(2))$. 
In this case $\Lambda (R)=(32111)$, and $K_{\ld R}(t)=t^5+3t^6+2t^7+t^8$, 
but $K_{\Lambda (R)\ld}(t)=2q^3+3q^4+2q^5$. 

Summarizing, it seems that there are no 
simple connection between the ribbon and generalized Kostka polynomials 
in general. Nevertheless, for the so--called dominant sequences of 
rectangular
partitions $R$, one can conjectured (see, e.g., [KS]) that the generalized 
and ribbon Kostka polynomials coincide. Recall that a sequence of 
partitions $R=(R_1,\ldots ,R_p)$ is called {\it dominant}, if for all 
$1\le i\le p-1$, the last part of $R_i$ is at least as large as the first 
part of $R_{i+1}$.

\vskip 0.2cm
{\bf Conjecture 6.8} ([KS]). {\it Let $R=(R_1,\ldots ,R_p)$ be a dominant 
sequence of rectangular partitions, and $\Lambda =\Lambda (R)$ be the unique 
partition with empty $p$--core and $p$--quotient $(R_1,\ldots ,R_p)$. Then}
$$K_{\ld R}(t)=K_{\Lambda\ld}^{(p)}(t). \eqno (6.18)
$$
\vskip 0.2cm

More generally, let $\ld$ be a partition with empty $p$--core and 
$p$--quotient $(\ld^{(1)},\ldots ,\ld^{(p)})$. Partition $\ld$ is called 
{\it $p$--dominant}, if there exists a permutation $s\in S_p$ such that\break 
$R=(\ld^{(s(1))},\ldots ,\ld^{(s(p))})$ is the dominant sequence of 
partitions.

\vskip 0.2cm
{\bf Conjecture 6.9.} {\it Let $\ld$ be a $p$--dominant partition, and 
$R:=R(\ld )$ be the dominant sequence of partitions obtained by 
rearrangement of the $p$--quotient of $\ld$. Assume that all partitions 
in the sequence $R$ have rectangular form, then for any partition $\mu$}
$$K_{\ld\mu}^{(p)}(t)=K_{\mu R}(t). \eqno (6.19)
$$
\vskip 0.2cm

{\bf Problem 8.} Let $\ld$ be a $p$--dominant partition, and $R:=R(\ld )$ 
be the dominant rearrangement of the $p$--quotient of $\ld$. For which 
partition $\mu$, the $p$--ribbon Kostka polynomial $K_{\ld\mu}^{(p)}(t)$ 
coincides with generalized Kostka polynomial $K_{\mu R}(t)$?
\vskip 0.2cm

As for the fusion Kostka polynomials $FK_{\ld R}(t)$, their connection 
with the corresponding crystal, ribbon or generalized Kostka polynomials 
is unclear.

Finally, let us consider few examples which illustrate the difference 
between the ribbon, fusion and generalized Kostka polynomials.

\vskip 0.2cm
{\bf Examples.} $i)$ Let $R=(R_1,R_2)$ be a dominant sequence of 
partitions. One can show that in this case
$$K_{\ld R}(q)={\rm Mult}(V_{\ld}:V_{R_1}\otimes V_{R_2})q^{E_0},\eqno 
(6.20)
$$
for a certain constant $E_0:=E(\ld R)$. However, the corresponding fusion 
and ribbon Kostka polynomials contain "in general"  more than one term.

$ii)$ Take $p=4$ and $\ld =(8,8,8,4,4)$. Then 4--quotient 
$\ld =((2,1),(1),(2),(2))$ and 4--core $(\ld )=\emptyset$. Hence, we 
see that $\ld$ is the 4--dominant partition, and $R:=R(\ld 
)=((2),(2),(2,1),(1))$. One can check that if $\mu =(4211)$, then
$$K_{\ld\mu}^{(4)}(q)=K_{\mu R(\ld )}(q)=q^2+2q^3+3q^4+q^5.
$$
If we take the same $p$ and $\ld$, but take $\mu =(4,2,2)$, then
$$K_{\ld\mu}^{(4)}(q)=K_{\mu R(\ld )}(q)=2q^3+2q^4+2q^5.
$$
However, if we take $p=4$, $\ld =(8,8,8,4,4)$, $\mu =(4,2,2)$, but take 
$R=((2,1),(1),(2),(2))$ $=p$--quotient $(\ld )$, then $K_{\mu 
R}(q)=q^6+3q^7+q^8+q^9\ne K_{\ld\mu}^{(4)}(q)$.

$iii)$ Take $p=2$ and $\ld =(11,9,9,7,7,5,2)$, then 2--core$(\ld 
)=\emptyset$ and 2--quotient of $\ld$ is equal to $((4,3,2,1),(6,5,4))$; 
hence, $\ld$ is the 2--dominant partition, and $R:=R(\ld 
)=((6,5,4),(4,3,2,1))$. Now let us take $\mu =(8,6,6,3,2)$, then
$$\eno{
K_{\mu R}(q)&=7q^5,\cr
FK_{\mu R}(q)&=2q^6+4q^7+q^8,\cr
K_{\ld\mu}^{(2)}(q)&=3q^4+4q^5=K_{\Lambda (R)\mu}^{(2)}(q).}
$$

These examples show that for a general $p$--dominant partition $\ld$ with 
the dominant rearrangement of the $p$--quotient $R:=R(\ld )$, the 
ribbon and generalized Kostka polynomials $K_{\ld\mu}^{(p)}(t)$ and 
$K_{\mu R}(t)$ give unequivalent $q$--analogues of the tensor product 
multiplicities.

\vskip 0.5cm
{\bf \S 7. Fermionic formulae.}
\vskip 0.3cm

By fermionic formulae for polynomial (or series) $f(t)\in{\bf N}[t]$ we 
roughly mean such expression for $f(t)$ which is free of signs, admits 
a quasi--particle interpretation, has an origin in the Bethe ansatz, etc. 
Thanks to the absence of signs, fermionic formulae are suitable for 
studying the limiting behavior and serve as a key to establish various 
formulae for the characters related to the affine Lie algebras and 
Virasoro algebra, see [Ki2], [HKKOTY] for examples illustrating this 
thesis.

%\vfil\eject
\vskip 0.3cm
{\bf 7.1.} Multinomial fermionic formulae for one dimensional sums.
\vskip 0.2cm

The starting point of our investigation is a simple observation that the 
number of transport 
matrices ${\cal P}_{\ld\mu}(1)$ of type $(\ld ;\mu )$ is equal to the 
coefficient of $x^{\mu}$ in the product 
$$h_{\ld_1}(X_n)\ldots h_{\ld_p}(X_n),
$$ 
where $\ld =(\ld_1,\ldots ,\ld_p)$, $p\le n$, $l(\mu 
)\le n$, and $h_k(X_n)$ denotes the complete homogeneous symmetric 
function of degree $k$ in the variables $X_n=(x_1,\ldots ,x_n)$. 

More 
generally, the number ${\cal P}_{R\mu}(1)$ of weight $\mu$ unrestricted 
paths corresponding to the tensor product of crystals 
$B_{R_1}\otimes\cdots\otimes B_{R_p}$ is equal to the coefficient of 
$x^{\mu}$ in the product of Schur functions $s_{R_1}\cdots s_{R_p}$. This 
is clear. On the other hand it is well--known and goes back to 
D.~Littlewood (see, e.g., [CL], [LLT], [SW]) that the latter coefficient 
is equal also 
to the number $|{\rm Tab}_p(\Lambda ,\mu )|$ of $p$--ribbon tableaux of 
shape $\Lambda$ and weight $\mu$, where $\Lambda$ is the unique partition 
which satisfies the conditions (6.10) and (6.11). 

\vskip 0.2cm
{\bf Problem 9} (cf. Conjecture~6.4). To construct a bijection $\psi$ 
between the sets ${\cal P}_{R\mu}(1)$ and ${\rm Tab}_p(\Lambda ,\mu )$ 
which transforms the energy function $E$ on the set ${\cal P}_{R\mu}(1)$ 
to the modified spin function $\bar s$ on that ${\rm Tab}_p(\Lambda ,\mu )$.

\vskip 0.2cm
{\bf Problem 10.} Let ${\cal P}_{R\mu}^0(1)$ be the set of weight $\mu$ 
classically restricted paths corresponding to the tensor product of 
crystals $B_{R_1}\otimes\cdots\otimes B_{R_p}$, see, e.g., [KMOTU2], or 
[HKKOTY]. It is clear that 
${\cal P}_{R\mu}^0(1)\subset{\cal P}_{R\mu}(1)$. Find characterization 
of  the subset 
${\rm Tab}_p^0(\Lambda ,\mu )\subset{\rm Tab}_p(\Lambda ,\mu )$ which 
corresponds to that ${\cal P}_{R\mu}^0(1)$ under the above bijection 
$\psi$.
\vskip 0.2cm

In the case $p=2$ the set ${\rm Tab}_2^0(\Lambda ,\mu)$ was characterized 
by C.~Carre and B.~Leclerc [CL] as the set of Yamanouchi domino tableaux. 
A weight preserving bijection between the set of domino tableaux of a 
fixed shape and that of ordinary tableaux of a related fixed shape, which 
maps Yamanouchi domino tableaux to ordinary Yamanouchi tableaux, was 
constructed by M.~van~Leeuwen [Le]. The question whether or not the 
bijection constructed by Leeuwen transforms the spin of a domino tableau 
to the value of the energy function for corresponding path is still open.

Let us continue and note that there exists yet another way to describe 
the coefficient of $x^{\mu}$ in the product of Schur functions 
$s_{R_1}\cdots s_{R_p}$ which is based on a combinatorial formula for 
Schur functions, see, e.g., [M], Chapter~I, (5.11). We consider here only the 
so--called "homogeneous case" $R_1=\cdots =R_p:=\ld$. The general case 
can be treated similarly. The starting point for obtaining the
multinomial fermionic formulae is the following combinatorial formula for 
the Schur functions mentioned above: 

let $\ld$ be a partition, $l(\ld )\le n$, then 
$$s_{\ld}(X_n)=\sum_Tx^{w(T)},
$$ 
where the sum runs over the set $STY(\ld ,\le n)$ of all semistandard 
Young tableaux of shape $\ld$ filled by 
numbers not exceeding $n$, and $w(T)$ is the weight of a tableau $T$. 

Let us define the multinomial coefficient $\pmatrix{L\cr \mu}_{\ld}$ via 
decomposition
$$\left( s_{\ld}(X_n)\right)^L=\sum_{\mu}\pmatrix{L\cr\mu}_{\ld}x^{\mu},
\eqno (7.1)
$$
where the sum is taken over the set of all compositions $\mu$ such that 
$l(\mu )\le n$ and $|\mu |=L|\ld |$;
$$\pmatrix{L\cr\mu}_{\ld}:=\sum_{\{ k_T\}}\pmatrix{L\cr \{ k_T\}},
\eqno (7.2)
$$
summed over all sequences of nonnegative integers $\{ k_T\}$ 
parameterized by the set of semistandard Young tableaux $STY(\ld ,\le n)$, 
such that $\ds\sum_{T\in STY(\ld ,\le n)}w(T)k_T=\mu$;

\hskip -0.7cm$\pmatrix{L\cr \{ k_T\}}:=\ds{L!\over \prod_T(k_T)!}$ 
stands for the gaussian multinomial coefficient.

The multinomial coefficients $\pmatrix{L\cr\mu}_{\ld}$ can be characterized 
by the following properties:

$\bullet$ $\pmatrix{0\cr\mu}_{\ld}=\delta_{0,\mu}$ (initial data);

$\bullet$ $\pmatrix{L+1\cr\mu}_{\ld}=\ds\sum_{T\in STY(\ld ,\le n)}
\pmatrix{L\cr\mu-w(T)}_{\ld}$ (recurrence relations), \hskip 2cm (7.3)

\hskip -0.7cm where we assume that the gaussian multinomial coefficient 
$\pmatrix{L\cr m_1,\ldots ,m_n}$ is equal to 0, if $m_i<0$ for some $i$.

It is natural to ask: what is a $q$--analog of the multinomial coefficient 
$\pmatrix{L\cr\mu}_{\ld}$, and what are the $q$--analogues of relations 
(7.1), (7.2), and (7.3)? The answers on these questions are either 
well--known or conjectured.
 
More precisely, for each semistandard tableau $T\in STY(\ld ,\le n)$ let
us denote by $\left[\matrix{L\cr\mu}\right]_{\ld}^{(T)}$ the weight $\mu$ 
unrestricted one dimensional sum with boundary condition $b_T\in B_{\ld}$.
Let $H:B_{\ld}\times B_{\ld}\to{\bf Z}$ stands for the local energy 
function corresponding to the crystal $B_{\ld}$, see, e.g., [Ka1], [Ka2].
In the sequel we will identify the sets $B_{\ld}$ and $STY(\ld ,\le n)$.

It is well--known, see, e.g., [KMOTU2], that one dimensional sums 
$\left[\matrix{L\cr\mu}\right]_{\ld}^{(T)}$ satisfy the following 
conditions:

$\bullet$ $\left[\matrix{L\cr\mu}\right]_{\ld}^{(T)}\vert_{q=1}=
\pmatrix{L\cr\mu -w(T)}_{\ld}$, 

$\bullet$ $\left[\matrix{L\cr\mu}\right]_{\ld}^{(T)}=
\delta_{0\mu}$ (initial datum),

$\bullet$ let $T_0\in STY(\ld ,\le n)$, then
$$\left[\matrix{L+1\cr\mu}\right]_{\ld}^{(T_0)}=
\sum_{T\in STY(\ld ,\le n)}q^{H(T_0,T)}\left[\matrix{L\cr\mu}
\right]_{\ld}^{(T)}. \eqno (7.4)
$$
(recurrence relations)

For $q$--analogue of (7.2) and (7.3), we are making the following 
conjectures:

\vskip 0.2cm
{\bf Conjecture 7.1.} {\it There exists a quadratic form 
$Q:B_{\ld}\times B_{\ld}\to{\bf Z}$, and a set of linear forms 
$l_T:B_{\ld}\to{\bf Z}$, $T\in B_{\ld}$, such that
$$\left[\matrix{L\cr\mu}\right]_{\ld}^{(T_0)}=\sum_{\{ k_T\}}
q^{\sum_{T,T'}Q(T,T')k_Tk_{T'}+\sum_Tl_{T_0}(T)k_T}\left[\matrix{L\cr
\{ k_T\}}\right]_q, \eqno (7.5)
$$
summed over all sequences of nonnegative integers $\{ k_T\}$, 
$T\in STY(\ld ,\le n)$, such that $\ds\sum_Tw(T)k_T=\mu$;
$\left[\matrix{L\cr\{ k_T\}}\right]_q=\ds{(q;q)_L\over\prod_T(q;q)_{k_T}}$ 
stands for a $q$--analog of the gaussian multinomial coefficient.}

\vskip 0.2cm
{\bf Remark.} The answer to this conjecture is known or conjectured in 
the case when partition $\ld =(l)$ consists of one part and for some 
special values of $T\in STY((l),\le n)$.
\vskip 0.2cm

{\bf Conjecture 7.2.} {\it Let $\ld =(\ld_1,\ldots ,\ld_s)$ be a partition.
For each integer $k\ge 1$, denote by $k\ld$ the following partition 
$(k\ld_1,\ldots ,k\ld_s)$. Then
$$q^{Ln(\ld')}{\cal P}_{L\ld ,\mu}^{(L)}(q)=\left[\matrix{L\cr\mu}
\right]_{\ld}^{(T_{\max})}, \eqno (7.6)
$$
where $T_{\max}$ denotes the unique maximal with respect to the
lexicographic order element in the set $STY(\ld ,\le n)$.}
\vskip 0.2cm

In other words, let $H_{\ld}^{(L)}(X_n;t)$ be the $H$--function defined 
in [LLT], Section~6; see also [CL] and [KLLT]. Then
$$H_{\ld}^{(L)}(X_n;t)=q^{-Ln(\ld')}\sum_{\mu}\left[\matrix{L\cr\mu}
\right]_{\ld}^{(T_{\max})}m_{\mu}(X_n). \eqno (7.7)
$$

{\bf Remark.} If partition $\ld =(l)$ consists of one part, the multinomial 
coefficients $\left[\matrix{L\cr\mu}\right]_{\ld}^{(T)}$ coincide with 
those introduced by A.~Schilling and S.O.~Warnaar after changing $q$ to 
$q^{-1}$ and multiplication on some power of $q$, [Ki2], [Sc], [ScW], [W].

%\vfil\eject
\vskip 0.3cm
{\bf 7.2.} Rigged configurations polynomials.
\vskip 0.2cm

In the previous subsection we explained an origin of a 
(conjectural) multinomial fermionic formulae for polynomials ${\cal 
P}_{L\ld ,\mu}^{(L)}(t)$. In this subsection we are going to to present 
yet another example of a (conjectural) 
fermionic formula for generalized Kostka polynomials corresponding to a 
collection of rectangles $R=(R_1,\ldots ,R_p)$, where 
$R_a=(\eta_a^{\mu_a})$ for all $1\le a\le p$. Our approach is using 
the so--called rigged configurations polynomials. 

\vskip 0.2cm
{\bf Definition 7.3.} {\it Let $\ld$ be a partition such 
that $|\ld |=\ds\sum_a\mu_a\eta_a$. We define the rigged configurations 
polynomial $RC_{\ld R}(q)$ to be
$$RC_{\ld R}(q)=\sum_{\{\nu\}}q^{c(\{\nu\})}\prod_{k\ge 1}\prod_{i\ge 1}
\left[\matrix{P_i^{(k)}(\nu )+m_i(\nu^{(k)})\cr m_i(\nu^{(k)})}\right]_q, 
\eqno (7.8)
$$
where sum runs over all sequences of partitions $\nu 
=\{\nu^{(1)},\nu^{(2)},\ldots\}$ such that}

$\bullet$ $|\nu^{(k)}|=\ds\sum_{j\ge k+1}\ld_j-\sum_{a=1}^p\mu_a\theta 
(\eta_a-k)$;

$\bullet$ $P_i^{(k)}(\nu ):=\ds\sum_{a=1}^p\min 
(i,\mu_a)\delta_{k,\eta_a}+Q_i(\nu^{(k-1)})-2Q_i(\nu^{(k)})+
Q_i(\nu^{(k+1)})\ge 0$ for $i,k\ge 1$.
\vskip 0.2cm

We have used the following notations:

\vskip 0.2cm
i) for any partition $\ld$, $Q_j(\ld )=\ds\sum_{i\le j}\ld_i'=\sum_{i\ge 
1}\min (j,\ld_i)$;

ii) if $x\in{\bf R}$, then $\theta (x)=1$, if $x\ge 0$, and $\theta (x)=0$, 
if $x<0$;

iii) $\nu^{(0)}=\emptyset$;

iv) $m_i(\ld )=\ld_i'-\ld_{i+1}'$ is the number of parts equal to $i$ of 
the partition $\ld$;

v) $c(\{\nu\} )=\ds\sum_{k\ge 1}\sum_{i\ge 1}\pmatrix{A_{ik}\cr 2}$, where 
$A_{ik}=(\nu^{(k-1)})'_i-(\nu^{(k)})'_i+\ds\sum_{a=1}^p\theta (\eta_a-k)
\theta (\mu_a-i).$

\vskip 0.2cm
{\bf Conjecture 7.4} (A.N.~Kirillov, M.~Shimozono, [KS]). {\it Let $\ld$ 
and $R$ be as above, and $\Lambda$ be the unique partition which 
satisfies the conditions (6.10) and (6.11). Assume that 
$\mu_1\ge\mu_2\cdots\ge\mu_p$, then}
$$K_{\Lambda\ld}^{(p)}(q)=RC_{\ld R}(q).
$$
\vskip 0.1cm

It is known (A.N.~Kirillov) that
$$RC_{\ld R}(1)={\rm Mult}[V_{\ld}:V_{R_1}\otimes\cdots\otimes V_{R_p}],
\eqno (7.9)
$$
where for each $1\le a\le p$, $R_a=\left(\eta_a^{\mu_a}\right)$ is a 
rectangular partition. A combinatorial proof of (7.9) is based on the 
construction of rigged configurations bijection (A.N.~Kirillov). 

Let us illustrate the formula (7.8) by simple example:

\vskip 0.1cm
{\bf Example.} Take $\ld =(44332)$ and 
$R=\left((2^3),(2^2),(2^2),(1),(1)\right)$. Then $|\nu^{(1)}|=4$, 
$|\nu^{(2)}|=6$, $|\nu^{(3)}|=5$, and $|\nu^{(4)}|=2$. It is not hard to 
check that there exist 6 configurations. 
%\vfil\eject
They are:
\vskip 0.2cm
$$\matrix{ \{\nu\} :&&
\hbox{
          \normalbaselines\m@th\offinterlineskip
	         \vtop{\hbox{{\Hthreebox(,,)}~0}
	         \vskip-0.4pt
	         \hbox{{\fsquare(0.4cm,)}~0}}}
&&\hbox{
          \normalbaselines\m@th\offinterlineskip
          \vtop{\hbox{{\Hthreebox( , , )}~1}
          \vskip-0.4pt
          \hbox{\Hthreebox( , , )}}}	      
&&\hbox{
          \normalbaselines\m@th\offinterlineskip
          \vtop{\hbox{{\Hthreebox( , ,)}~0}
          \vskip-0.4pt
          \hbox{{\Htwobox( , )}~0}}}
&&\hbox{{\Htwobox(, )}~0} 
&& c(\{\nu\} )=10,\cr \cr \cr
 \{\nu\} :&&
\hbox{
          \normalbaselines\m@th\offinterlineskip
	         \vtop{\hbox{{\Hthreebox(,,)}~0}
	         \vskip-0.4pt
	         \hbox{{\fsquare(0.4cm,)}~1}}}
&&\hbox{
          \normalbaselines\m@th\offinterlineskip
          \vtop{\hbox{{\Hthreebox( , , )}~1}
          \vskip-0.4pt
          \hbox{{\Htwobox( , )}~1}
          \vskip-0.4pt
          \hbox{{\fsquare(0.4cm,)}~0}}}      
&&\hbox{
          \normalbaselines\m@th\offinterlineskip
          \vtop{\hbox{{\Hthreebox( , ,)}~0}
          \vskip-0.4pt
          \hbox{{\Htwobox( , )}~1}}}
&&\hbox{{\Htwobox(, )}~0} 
&& c(\{\nu\} )=8,\cr \cr \cr
\{\nu\} :&&
\hbox{
          \normalbaselines\m@th\offinterlineskip
	         \vtop{\hbox{{\Htwobox(,)}~0}
	         \vskip-0.4pt
	         \hbox{{\Htwobox(,)}}}}
&&\hbox{
          \normalbaselines\m@th\offinterlineskip
          \vtop{\hbox{{\Htwobox(,)}~0}
          \vskip-0.4pt
          \hbox{{\Htwobox(,)}}
          \vskip-0.4pt
          \hbox{{\Htwobox(,)}}}}	      
&&\hbox{
          \normalbaselines\m@th\offinterlineskip
          \vtop{\hbox{{\Hthreebox( , ,)}~0}
          \vskip-0.4pt
          \hbox{{\Htwobox( , )}~2}}}
&&\hbox{{\Htwobox(, )}~0} 
&& c(\{\nu\} )=8,\cr \cr \cr
\{\nu\} :&&
\hbox{{\Hfourbox(,,,)}~0}
&&\hbox{
          \normalbaselines\m@th\offinterlineskip
          \vtop{\hbox{{\Hthreebox( , , )}~0}
          \vskip-0.4pt
          \hbox{{\Hthreebox( , , )}}}}	      
&&\hbox{
          \normalbaselines\m@th\offinterlineskip
          \vtop{\hbox{{\Hthreebox( , ,)}~0}
          \vskip-0.4pt
          \hbox{{\Htwobox( , )}~0}}}
&&\hbox{{\Htwobox(, )}~0} 
&& c(\{\nu\} )=12,\cr \cr \cr
\{\nu\} :&&
\hbox{
          \normalbaselines\m@th\offinterlineskip
	         \vtop{\hbox{{\Htwobox(,)}~0}
	         \vskip-0.4pt
	         \hbox{{\fsquare(0.4cm,)}~0}
	         \vskip-0.4pt
	         \hbox{{\fsquare(0.4cm,)}}}}
&&\hbox{
          \normalbaselines\m@th\offinterlineskip
          \vtop{\hbox{{\Htwobox( , )}~1}
          \vskip-0.4pt
          \hbox{{\Htwobox( , )}}
          \vskip-0.4pt
	         \hbox{{\fsquare(0.4cm,)}~0}
	         \vskip-0.4pt
	         \hbox{{\fsquare(0.4cm,)}}}}	      
&&\hbox{
          \normalbaselines\m@th\offinterlineskip
          \vtop{\hbox{{\Htwobox( ,) }~0}
          \vskip-0.4pt
          \hbox{{\Htwobox(,)}}
          \vskip-0.4pt
          \hbox{{\fsquare(0.4cm,)}~0}}}
&&\hbox{{\Htwobox(, )}~1} 
&& c(\{\nu\} )=6,\cr \cr \cr
\{\nu\} :&&
\hbox{
          \normalbaselines\m@th\offinterlineskip
	         \vtop{\hbox{{\Htwobox(,)}~0}
	         \vskip-0.4pt
	         \hbox{{\fsquare(0.4cm,)}~0}
	         \vskip-0.4pt
	         \hbox{{\fsquare(0.4cm,)}}}}
&&\hbox{
          \normalbaselines\m@th\offinterlineskip
          \vtop{\hbox{{\Htwobox( , )}~0}
          \vskip-0.4pt
          \hbox{{\Htwobox( , )}}
          \vskip-0.4pt
	         \hbox{{\fsquare(0.4cm,)}~0}
	         \vskip-0.4pt
	         \hbox{{\fsquare(0.4cm,)}}}}	      
&&\hbox{
          \normalbaselines\m@th\offinterlineskip
          \vtop{\hbox{{\Hthreebox( , ,) }~0}
          \vskip-0.4pt
          \hbox{{\fsquare(0.4cm,)}~0}
          \vskip-0.4pt
          \hbox{{\fsquare(0.4cm,)}}}}
&&\hbox{{\Htwobox(, )}~0} 
&& c(\{\nu\} )=8,
}         	            
$$
\vskip 0.1cm
Thus, the rigged configurations polynomial $RC_{\ld R}(q)$ is equal to
$$\eno{
&q^{10}\left[\matrix{3\cr 1}\right] +q^8\left[\matrix{2\cr 1}\right]
\left[\matrix{2\cr 1}\right]\left[\matrix{2\cr 1}\right]
\left[\matrix{2\cr 1}\right] +q^8\left[\matrix{3\cr 1}\right] +q^{12}
+q^6\left[\matrix{2\cr 1}\right]\left[\matrix{3\cr 1}\right] +q^8\cr \cr
&=q^6+2q^7+5q^8+6q^9+8q^{10}+5q^{11}+3q^{12}.}
$$

%\vfil\eject
\vskip 0.5cm
{\bf \S 8. Two parameter deformation of one dimensional sums.}
\vskip 0.3cm

Following [GH], we define the modified Macdonald polynomials in {\it infinite
number of variables} by 
$$\wt P_{\ld}(x;q,t)=P_{\ld}\left({x\over 1-t}; q,t\right),~~~~
\wt J_{\ld}(x;q,t)=J_{\ld}\left({x\over 1-t};q,t\right) \eqno (8.1)
$$
in the $\ld$--ring notation. Let us explain briefly the $\ld$--ring 
notation in the context of symmetric functions. Given a symmetric 
function $f(x)=f(x_1,x_2,\ldots )$ in infinite set of variables 
$x=(x_1,x_2,\ldots )$, the symbol 
$\ds f\left({x\over 1-t}\right)$ in the $\ld$--ring notation stands 
for the symmetric function $f(\wt x)$ obtained by the transformation 
of variables $\wt x=(x_it^j)_{i\ge 1, j\ge 0}$. In infinite number of variables, 
the symmetric function $f(x)$ can be written uniquely in the form 
$f(x)=\varphi (p_1(x),p_2(x),\ldots )$ as a polynomial of the power 
series $p_k(x)=\ds\sum_{j=1}^{\infty}x_j^k$, $k=1,2,\ldots$. Then the 
symbol $\ds f\left({x\over 1-t}\right)$ represents the symmetric function
$$f\left({x\over 1-t}\right)=\varphi\left({p_1(x)\over 1-t},
{p_2(x)\over 1-t^2},\ldots\right),
$$
obtained by the transformation $p_k(x)\to p_k(x)/(1-t^k)$, $k=1,2,\ldots$.
When we consider the modified Macdonald polynomials in $n$ variables 
$X_n=(x_1,\ldots ,x_n)$, each function $\wt P_{\ld}(X_n;q,t)$ and 
$\wt J_{\ld}(X_n;q,t)$ should be understood as the one obtained from the
corresponding symmetric function in infinite number of variables by setting 
$x_{n+1}=x_{n+2}=\cdots =0$.

An advantage of modified Macdonald polynomials is that they have nice 
transformation coefficients with classical Schur functions $s_{\ld}(x)$:
$$\wt J_{\ld}(x;q,t)=\sum_{\ld}K_{\ld\mu}(q,t)s_{\ld}(x), \eqno (8.2)
$$
where $K_{\ld\mu}(q,t)$ are the double Kostka coefficients. It is 
well--known that $K_{\ld\mu}[q,t]\in{\bf Z}[q,t]$ for 
all $\ld$ and $\mu$, see, e.g., [KN].

Our next aim is to construct two parameter deformation of polynomials
${\cal P}_{\ld\mu}(t)$ using the modified Macdonald polynomials 
$\wt J_{\ld}(X_n;q,t)$ instead of modified Hall--Littlewood polynomials
$Q_{\ld}'(X_n;t)$. To this end, let us consider, follow [KN], a family 
of polynomials $B_{\ld\mu}(q,t)$ via decomposition
$$\wt J_{\ld}(x;q,t)=\sum_{\mu}B_{\ld\mu}(q,t)m_{\mu}(x).
$$
It is clear that
$$B_{\ld\mu}(q,t)=\sum_{\eta}K_{\eta\mu}K_{\eta\ld}(q,t).
$$
Let us formulate some basic properties of polynomials $B_{\ld\mu}(q,t)$.
For further details and proofs, see [M], Chapter~VI, and [KN], Section~8.

Let $\ld$ and $\mu$ be a partitions of a given natural number $N$, then

\vskip 0.2cm
$\bullet$ $B_{\ld\mu}(1,1)=\pmatrix{N\cr\mu_1,\mu_2,\ldots}=
{\cal P}_{1^N\mu}(1)$;
\vskip 0.2cm

$\bullet$ $B_{\ld\mu}(0,t)={\cal P}_{\ld\mu}(t)$;
\vskip 0.2cm

$\bullet$ $B_{(N)\mu}(q,t)=q^{n(\mu')}\left[\matrix{N\cr\mu_1,\mu_2,
\ldots}\right]_q={\cal P}_{1^N\mu}(q)$;
\vskip 0.2cm

$\bullet$ $B_{\ld'\mu}(q,t)=q^{n(\ld')}t^{n(\ld )}B_{\ld\mu}
(t^{-1},q^{-1})$ (duality).
\vskip 0.2cm

It follows from duality that
$$B_{\ld\mu}(q,t)=q^{n(\ld')}(\wt{\cal R}_{\ld\mu}(t)+o(q^{-1})),
$$
where $\wt{\cal R}_{\ld\mu}(t)=\ds\sum_{\eta}K_{\eta\mu}
\wt K_{\eta'\ld}(t)$, and $\wt K_{\eta'\ld}(t)=t^{n(\ld 
)}K_{\eta'\ld}(t^{-1})$.

The properties of polynomials $B_{\ld\mu}(q,t)$ mentioned above show 
that they can be considered as a natural two parameter deformation
of the gaussian multinomial coefficients.

\vskip 0.2cm
{\bf Problem 11.} Find a ``path realization'' of polynomials 
$B_{\ld\mu}(q,t)$.
\vskip 0.2cm

{\bf Problem 12} (``Parabolic modified Macdonald polynomials''). Find two 
parameter deformation of polynomials $\wt Q_{\ld}^{(p)}(X_n;t)$ with
nice combinatorial, algebraic and geometric properties.
\vskip 0.2cm

At the end of this Section we give an example of polynomials 
$B_{\ld\mu}(q,t)$. 

\vskip 0.2cm
{\bf Example.} Take $\ld =(2^3)$, $\mu =(2^21^2)$, then
$$\eno{
B_{\ld\mu}(q,t)&=\sum_{\eta}K_{\eta\mu}K_{\eta\ld}(q,t)=
1+4t+8t^2+9t^3+7t^4+3t^5+t^6\cr
&+q\left[\matrix{3\cr 1}\right]_t
(1+5t+9t^2+7t^3+3t^4)+q^2\left[\matrix{3\cr 2}\right]_t(2t+6t^2+7t^3+4t^4)\cr
&+q^3(t^2+3t^3+5t^4+4t^5+2t^6).}
$$
Using the fermionic formulae (3.1) and (3.7), one can check that
$${\cal P}_{\ld\mu}(t)=B_{\ld\mu}(0,t)=1+4t+8t^2+9t^3+7t^4+3t^5+t^6,
$$
and
$${\cal R}_{\ld\mu}(t)=q^3t^6B_{\ld\mu}(q^{-1},t^{-1})|_{q=0}=
\sum_{\eta}K_{\eta\mu}K_{\eta'\ld}(t)=2+4t+5t^2+3t^3+t^4.
$$

\vfil\eject
%\vskip 0.5cm
%\vskip 1cm
{\bf References.}
\vskip 0.5cm

\item {[An]} Andrews G., {\it The theory of partitions,} Addison--Wesley 
Publishing Company, 1976.

\item {[BKMW]} Begin L., Kirillov A.N., Mathieu P. and Walton M.A., 
{\it Berenstein--Zelevinski triangles, elementary couplings and fusion 
rules,} Lett. in Math. Phys., 1993, v.28, p.257-268.

\item{[Bi]} Birkhoff G., {\it Subgroups of abelian groups,} Proc. London 
Math. Soc. (2), 1934-5, v.38, p.385-401.

\item{[Bu1]}  Butler L., {\it Subgroup lattices and symmetric functions,} 
Memoirs of AMS, 1994, v. 112, n. 539.

\item{[Bu2]} Butler L., {\it Generalized flags in finite abelian 
$p$--groups,} Discrete Appl. Math., 1991, v.112, p.67-81.

\item{[Bu3]} Butler L., {\it A unimodality result in the enumeration of 
subgroups of a finite abelian group,} Proc. Amer. Math. Soc., 1987, 
v.101(4), p.771-775.

\item {[CL]} Carr\'e C. and Leclerc B., {\it Splitting the square of a
Schur function into its symmetric and antisymmetric parts,}
J. of Alg. Combin., 1995, v.4, p.201-231.
      
\item{[DLT]} D\'esarmenien J., Leclerc B. and Thibon J.-Y., {\it 
Hall--Littlewood functions and Kostka--Foulkes polynomials in 
representation theory,} S\"eminaire Lotharingien de Combinatoire, 1994, 
v.32.

\item{[De]} Delsarte S., {\it Fonctions de M\"obius sur les groupes 
abelian finis,} Annals of Math., 1948, v.49, p.600-609.

\item{[Dy]} Dyubyuk P., {\it On the number of subgroups of a finite 
abelian group,} Izv. Akad. Nauk USSR, Ser. Mat., 1948, v.12, p.371-328.

\item{[Fi]} Fishel S., {\it Nonnegativity results for generalized 
$q$--binomial coefficients,} PhD thesis, the University of Minnesota, 1993.

\item {[F]} Foata D., {\it Distribution eul\'eriennes et mahoniennes sur 
le group des permutations,} in {\it Higher Combinatorics,} M. Aigner ed., 
Amsterdam, D.~Reidel, 1977, p.27-49.

\item {[FZ]} Foata D. and Zeilberger D., {\it Denert's permutation 
statistic is indeed Euler--Mahonian,} Studies in Applied Math., 1990, 
v.83, p.31-59.

\item{[GaW]} Galovich J. and White D., {\it Recursive statistics on 
words,} Discrete Math., 1996, v.157, p.169-191.

\item {[GH]} Garsia A. and Haiman M., {\it A graded representation model 
for Macdonald's polynomials}, Proc. Nat. Acad. Sci. USA, 1993,
v.90, p.3607-3610.

\item{[GoW]} Goodman F. and Wenzl H., {\it Littlewood--Richardson 
coefficients for Hecke algebras at roots of unity,} Adv. Math., 1990, 
v.82, p.244-265.

\item{[Gu]} Gupta R., {\it Generalized exponent via Hall--Littlewood 
symmetric functions,} Bull. Amer. Math. Soc., 1987, v.16, p.287-291.

\item{[HKKOTY]} Hatayama G, Kirillov A.N., Kuniba A., Okado M., Takagi 
T. and Yamada Y., {\it Character formulae of $\wh{sl_n}$--modules and 
inhomogeneous paths}, Preprint math.QA/9802085, 1998, 42p.

\item{[Hi]} Hilton H., {\it On subgroups of a finite abelian group,} 
Proc. london Math. Soc (2), 1907, v.5, p.1-5.

\item{[HS]} Hotta R. and Shimomura N., {\it The fixed point subvarieties 
of unipotent transformations on generalized flag varieties and Green 
functions,} Math. Ann., 1979, v.241, p.193-208.

\item{[H]} Huppert B., Endliche Gruppen, vol.1, Spring--Verlag, 1967.

\item{[Kac]} Kac V.G., {\it Infinite dimensional Lie algebras,} Cambridge 
University Press, 1990.

\item{[Ka1]} Kashiwara M., {\it On crystal bases of the $q$--analogue of 
universal enveloping algebras,} Duke Math. J., 1991, v.63, p.465-516.

\item{[Ka2]} Kashiwara M., {\it The crystal base and Littlemann's refined 
Demazure character formula,} Duke Math. J., 1993, v.71, p.839-858.

\item{[Ki1]} Kirillov A.N., {\it On the Kostka--Green--Foulkes polynomials 
and the Clebsch--Gordan numbers,} 
Journ. Geom. and Phys., 1988, v.5, n.3, p.365-389.

\item{[Ki2]} Kirillov A.N., {\it Dilogarithm identities,} Progress of Theor. 
Phys. Suppl., 1995, v.118, p.61-142.

\item{[KKN]} Kirillov A.N., Kuniba A. and Nakanishi T., {\it Skew Young 
diagram method in spectral decomposition of integrable
      lattice models II: Higher levels,} 
      Preprint q-alg/9711009, 1997, 27p.

\item {[KLLT]} Kirillov A,N., Lascoux A., Leclerc B. and Thibon J.-Y.,
{\it S\'eries g\'en\'eratrices pour les tableaux de dominos,} 
      C.R. Acad. Sci. Paris, 1994, t.318, Serie I, p.395-400.

\item {[KN]} Kirillov A.N. and Noumi M., {\it Affine Hecke algebras and 
raising operators for Macdonald polynomials,} to appear in Duke Math. 
Journ.; q-alg/9605004, 1996, 35p. 

\item {[KS]} Kirillov A.N. and Shimozono M., {\it  A generalization of the 
Kostka--Foulkes polynomials,} Preprint math.QA/9803062, 1998, 37p.

\item{[Kn]} Knuth D.E., {\it Permutations, matrices and generalized Young 
tableaux,} Pacific J. Math., 1970, v.34, p.709-727.

\item{[KMOTU1]} Kuniba A., Misra K.C., Okado M., Takagi T. and Uchiyama J., 
{\it Crystals for Demazure modules of classical affine Lie algebras,} 
Preprint q-alg/9707014, 29p.

\item{[KMOTU2]} Kuniba A., Misra K.C., Okado M., Takagi T. and Uchiyama J.,
{\it Characters of Demazure modules and solvable lattice models,} 
Nuclear Phys., 1998, v.B510[PM], p.555-576.

\item{[LLT]} Lascoux A., Leclerc B. and Thibon J.-Y., {\it Ribbon 
tableaux, Hall--Littlewood functions, quantum affine algebras, and 
unipotent varieties,} J. Math. Phys., 1997, v.38, p.1041-1068.

\item{[LS]} Lascoux A. and Sch\"utzenberger M.-P., {\it Sur une 
conjecture de H.O.~Foulkes,} C.R. Acad. Sci. Paris, 1978, t.286A, 
p.323-324.

\item {[Le]} Leeuwen M. van, {\it Some bijective correspondence
involving domino tableaux,} Preprint MAS--R9708, 1997, 28p.

\item{[M]} Macdonald I., {\it Symmetric functions and Hall polynomials,} 
2nd ed., Oxford, 1995.

\item{[Ma]} MacMahon P.A., {\it Combinatorial Analysis,} I, II, Cambridge 
University Press, 1915, 1916 (reprinted by Chelsea, New York, 1960).

\item{[Mi]} Miller G., {\it On subgroups of an abelian group,} Annals of 
Math., 1904, v.6, p.1-6.

\item{[NY]} Nakayashiki A. and Yamada Y., {\it Kostka polynomials and 
energy functions in solvable lattice models,} Preprint q-alg/9512027; to 
appear in Selecta Mathematica.

\item{[R]} Regonati F., {\it Sui numeri dei sottogruppi di dato ordine 
dei $p$--gruppi abeliani finiti,} Instit. Lombardo Rend. Sci., 1988, 
v.A122, p.369-380.

\item{[Sc]} Schilling A., {\it Multinomials and polynomial bosonic forms 
for the branching functions of the 
$\wh{su}(2)_M\times\wh{su}(2)_N/\wh{su}(2)_{M+N}$ conformal coset 
models,} Nucl. Phys. B, 1996, v.467, p.247-271.

\item {[ScW]} Schilling A. and Warnaar S.O., {\it Supernomial coefficients, 
polynomial identities and $q$--series,} Preprint q-alg/9701007, 34p.

\item{[Sh]} Shimomura N., {\it A theorem of the fixed point set of a 
unipotent transformation of the flag manifold,} J. Math. Soc. Japan, 
1980, v.32, p.55-64.

\item{[St]} Stanley R., {\it Supersolvable lattices,} Algebra Universalis, 
1972, v.2, p.197-217.

\item {[SW]} Stanton D. and White D., {\it A Schensted algorithm 
for rim hook tableaux,} J. Comb. Theory, Ser. A, 1985, v.40, p.211-247.

\item{[T]} Terada I., {\it A generalization of the length -- {\rm Maj} 
symmetry and the variety of $N$--stable flags,} Preprint, 1993.

\item{[W]} Warnaar S.O., {\it The Andrews--Gordon identities and 
$q$--multinomial coefficients,} Comm. Math. Phys., 1997, v.184, p.203-232.

\item{[Y]} Yeh Y., {\it On prime power abelian groups,} Bull. Amer. Math. 
Soc., 1948, v.54, p.323-327.

\item{[ZB]} Zeilberger D. and Bressoud D., {\it A proof of Andrew's 
$q$--Dyson conjecture,} Discrete Math., 1985, v.54, p.201-224.

\end